\documentclass[11pt]{amsart}
\usepackage[margin=1in]{geometry}
\usepackage{amsmath,amsfonts,amssymb,amsthm,epsfig,epstopdf,array,mathrsfs,dsfont,enumerate, fancyhdr, graphicx, proof, comment, multicol,tabularx,float,tikz,xypic,hyperref,appendix,setspace,wrapfig,tikz-cd,ytableau}
\let\C\relax
\usepackage[color,matrix,arrow]{xy}
\newcommand{\ol}{\overline}
\newcommand{\ul}{\underline}
\newcommand{\cO}{\mathcal{O}}
\newcommand{\C}{\mathbb{C}}
\newcommand{\Z}{\mathbb{Z}}
\renewcommand{\cR}{\mathcal{R}}
\newcommand{\Q}{\mathbb{Q}}

\newcommand{\Fcal}{\mathcal{F}}
\newcommand{\zfrak}{\mathfrak{z}}
\newcommand{\Vcal}{\mathcal{V}}
\newcommand{\Ecal}{\mathcal{E}}
\newcommand{\Ocal}{\mathcal{O}}
\newcommand{\Bcal}{\mathcal{B}}
\newcommand{\Dcal}{\mathcal{D}}

\newcommand{\Hrm}{\mathrm{H}}

\newcommand{\fu}{\mathfrak{u}}
\newcommand{\Lrm}{\mathrm{L}}

\newcommand{\Hom}{\mathrm{Hom}}
\newcommand{\IC}{\mathrm{IC}}

\newcommand{\fq}{\mathfrak{q}}

\newcommand{\st}{\mathrm{st}}
\newcommand{\gr}{\mathrm{gr}}

\newcommand{\Ind}{\mathrm{Ind}}
\newcommand{\Frob}{\mathrm{Frob}}
\newcommand{\bfrak}{\mathfrak{b}}
\newcommand{\nfrak}{\mathfrak{n}}
\newcommand{\scrA}{\mathscr{A}}
\newcommand{\bR}{\mathbf{R}}
\newcommand{\hfrak}{\mathfrak{h}}
\newcommand{\gfrak}{\mathfrak{g}}

\newcommand{\tfrak}{\mathfrak{t}}
\newcommand{\ufrak}{\mathfrak{u}}
\newcommand{\mfrak}{\mathfrak{m}}

\newcommand{\lfrak}{\mathfrak{l}}
\newcommand{\Xfrak}{\mathfrak{X}}

\newcommand{\Pfrak}{\mathfrak{P}}

\newcommand{\Cfrak}{\mathfrak{C}}

\newcommand{\hrm}{\mathrm{h}}

\newcommand{\Nbb}{\mathbb{N}}
\newcommand{\Mbb}{\mathbb{M}}

\newcommand{\Ad}{\mathrm{Ad}}

\newcommand{\Scal}{\mathcal{S}}
\newcommand{\Dbb}{\mathbb{D}}
\newcommand{\Sym}{\mathrm{Sym}}

\newcommand{\Cscr}{\mathscr{C}}

\newcommand{\Ucal}{\mathcal{U}}

\newcommand{\Mcal}{\mathcal{M}}

\newcommand{\Abb}{\mathbb{A}}
\newcommand{\Arm}{\mathrm{A}}

\newcommand{\Pcal}{\mathcal{P}}

\newcommand{\End}{\mathrm{End}}
\newcommand{\GL}{\mathrm{GL}}
\newcommand{\SL}{\mathrm{SL}}

\newcommand{\Lcal}{\mathcal{L}}

\newcommand{\gl}{\mathfrak{gl}}
\newcommand{\Sp}{\mathrm{Sp}}
\newcommand{\bC}{\mathbf{C}}
\newcommand{\bT}{\mathbf{T}}
\renewcommand{\cL}{\mathcal{L}}
\newcommand{\quash}[1]{}

\newcommand{\fin}{\mathrm{fin}}
\newcommand{\Acal}{\mathcal{A}}
\newcommand{\wt}{\widetilde}
\newcommand{\e}{\mathrm{e}}

\newcommand{\Mod}{\mathrm{-mod}}
\newcommand{\Supp}{\mathrm{Supp}}

\newcommand{\Hilb}{\mathrm{Hilb}}

\newcommand{\IHilb}{\mathrm{IHilb}}

\newcommand{\Fl}{\mathcal{F}l}

\newcommand{\red}{\mathrm{red}}
\newcommand{\reg}{\mathrm{reg}}

\newcommand{\Coh}{\mathrm{Coh}}
\newcommand{\ch}{\mathrm{ch}}

\newcommand{\aff}{\mathrm{aff}}

\newcommand{\Ngen}{\mathrm{Ngen}}

\newcommand{\area}{\mathrm{area}}

\newcommand{\lcodinv}{\mathrm{lcodinv}}
\newcommand{\pdinv}{\mathrm{pdinv}}
\newcommand{\codinv}{\mathrm{codinv}}
\newcommand{\rs}{\mathrm{rs}}
\newcommand{\vla}{{\check{\lambda}}}

\newcommand{\fg}{\mathfrak{g}}
\newcommand{\des}{\mathrm{des}}
\newcommand{\fc}{\mathfrak{c}}
\newcommand{\git}{/\!\!/}
\newcommand{\scrL}{\mathscr{L}}
\newcommand{\fm}{\mathfrak{m}}

\newcommand{\cT}{\mathcal{T}}

\newcommand{\Lie}{\mathrm{Lie}}
\newcommand{\Aut}{\mathrm{Aut}}

\newcommand{\ybf}{\mathbf{y}}
\newcommand{\fa}{\mathfrak{a}}
\newcommand{\ft}{\mathfrak{t}}

\newcommand{\hess}{\mathrm{Hess}}
\newcommand{\cell}{\mathrm{cell}}

\renewcommand{\cH}{\mathcal{H}}
\newcommand{\SO}{\mathrm{SO}}

\newcommand{\fgl}{\mathfrak{gl}}
\newcommand{\fh}{\mathfrak{h}}

\newcommand{\bG}{\mathbb{G}}

\theoremstyle{plain}
\newtheorem{thm}{Theorem}[section]
\newtheorem{prop}[thm]{Proposition}
\newtheorem{cor}[thm]{Corollary}
\newtheorem{lem}[thm]{Lemma}
\newtheorem{conj}[thm]{Conjecture}
\theoremstyle{definition}
\newtheorem{defn}[thm]{Definition}
\theoremstyle{remark}
\newtheorem{rem}[thm]{Remark}
\newtheorem{ex}[thm]{Example}

\usepackage{nicematrix,microtype}
\usetikzlibrary{decorations.pathreplacing}

\counterwithout{figure}{section}
\counterwithout{figure}{subsection}

\definecolor{color1}{HTML}{318CE7}
\hypersetup{
	colorlinks=true,	
	linkcolor=color1,	
	urlcolor=black,
	linktoc=all,
	citecolor=color1
}
\title[Affine Springer fibers, Cherednik algebra modules and Hilbert schemes]{Affine Springer fibers, Cherednik algebra modules and Hilbert schemes: torus link case}
\author{Pablo Boixeda Alvarez}
\author{Thomas Hameister}
\author{Xinchun Ma}

\newcommand{\Ffrak}{\mathfrak{F}}
\newcommand{\loc}{\mathrm{loc}}
\newcommand{\pol}{\mathrm{pol}}
\newcommand{\Rep}{\mathrm{Rep}}
\usetikzlibrary{decorations.pathmorphing}

\begin{document}
\begin{abstract}
    For a pair of positive integers $M,N$ with GCD $d$, we relate the bigraded $S_N\times S_d$-Frobenius characters of the homology of the positive part of a homogeneous affine Springer fiber for $\GL_N$ of slope $M/N$ to a direct sum of the associated graded of minimal supported modules of type $A_{N-1}$ rational Cherednik algebra at the same slope, and to the equivariant K-theory of a coherent sheaf on the Hilbert scheme of points on the plane.  For the first object, we use the Hessenberg paving of Goresky-Kottwitz-Macpherson; for the second and third objects, we resort to an explicit Hodge module on $\mathfrak{sl}_N$. The bridge between them is a non-coprime shuffle theorem of Kim-Oh. 
\end{abstract}
\maketitle
\tableofcontents

\setcounter{tocdepth}{1}

  \section{Introduction}

Fix positive integers $M,N$, put $d=\gcd(M,N)$, and write $M=dm$, $N=dn$, and $c=M/N=m/n$. To $(M,N)$ we associate the following objects, each of a completely different flavor:
    
\begin{enumerate}
    \item The bigraded Frobenius character of the Borel-Moore homology of the positive part of a regular semisimple homogeneous affine Springer fiber for $\GL_N$ of slope $c=M/N$.    
    \item The $\C^*\times \C^*$-equivariant $K$-theory class of a coherent sheaf on the Hilbert scheme of points on $\C^2$ arising from a character $D$-module. 
    \item The symmetric function lift of the Khovanov-Rozansky (KhR) link homology of the torus link $T_{M,N}$. 
    \item The bigraded Frobenius character of representations of the type $A_{N-1}$ rational Cherednik algebra $\Hrm_{M/N}$ of minimal support.
\end{enumerate}

It was conjectured in \cite{gors} in the coprime case that the four characters above are equal. This was established for coprime $M,N$ by the work of \cite{hikita,oy,hogancamp-mellit,gmv,CM-shuffle,ma2}. In addition, computations of (1) and (4) were carried out in the case $M = kN$ is an integer multiple of $N$ in \cite{CM-hessenberg} and \cite{LBA}, respectively, and (3) was computed in terms of combinatorics of invariant sets in \cite{wilson} using recursive formulas for KhR knot homology in \cite{hogancamp-mellit,gmv}. Connections were also drawn between (1) and (4) in \cite[Conjecture 2]{TX} and between (1) and (3) in \cite{GGG} in the setting of the $\Delta$-conjecture. Moreover, \cite{gko} attaches to every affine Springer fiber at a semisimple element a quasi-coherent sheaf on $\Hilb^n(\C\times \C^*)$. 


The object of this paper is to compute the invariants (1), (2), and (4) in the non-coprime setting. In a companion paper, Donghyun Kim and Jaeseong Oh prove a version of the Shuffle Theorem for the non-coprime case \cite{KO1,KO2}, generalizing \cite{CM-shuffle,mellitshuffle}, which concludes the proof of the following:

\begin{thm}(Informal version of Theorem \ref{thm:character-comparison})
    \label{thm:main_intro}
    Up to some normalizations, (1)--(4) above are equal for all pairs of positive integers $M,N$.
\end{thm}

For (1), this involves computing the Hessenberg paving of \cite{GKM} in terms of the combinatorics of invariant sets \cite{wilson,gm,gmv}. A byproduct of this calculation is a natural generalization of the chromatic polynomial of a graph to cyclic partite graphs, and an affine paving for a family of generalized Hessenberg varieties. 

For (2) and (4), we generalize the methods of \cite{ma2} to the non-coprime setting. A new feature is a $S_d$-action coming from a Springer action on the associated $D$-module. We compute the Frobenius character of this action in terms of shuffle algebra generators by reducing to the trivial isotypic components.

\subsection{Generalized Hessenberg varieties and partite chromatic functions}

To compute the Borel-Moore homology of affine Springer fibers, we will use the Hessenberg paving of \cite{GKM}. A key difference with the coprime case considered in \cite{hikita} is the presence of generalized Hessenberg varieties.

Let $n = N/d$, $m = M/d$, and $0\leq b<n$ the integer such that $b = m\mod n$. Consider the action of $H=\GL_d^n$ on the set of $n$-tuples of $d\times d$ matrices by
\[
\gfrak_{b}:= \{A=(A_1,\dots, A_n)\mid A_i\text{ is a $d\times d$ matrix}\},\qquad h\cdot A = (h_iA_ih_{b+i}^{-1})_{i=1}^n,
\]
where the index of $h$ is taken modulo $n$. A generalized Hessenberg space $\fm\subset \fg_b$ is a product $(\fm^{(1)},\dots, \fm^{(n)})$, where each $\fm^{(i)}$ is a standard Hessenberg space in $\mathfrak{gl}_d$ determined by a $d\times d$ Dyck path. For $x\in \fg_{b}$, the cyclic Hessenberg variety $\cH(x,\fm)$ is the subvariety of cosets $h$ in the flag variety of $H$ such that $h^{-1}\cdot x\in \fm$. 

When $d=N$, $\cH(x,\fm)$ is a Hessenberg variety, which has been studied for regular $x$ in many works such as \cite{sw2,julianna,precup,BrosnanChow}. In particular, its cohomology or Borel-Moore homology admits an action of the symmetric group, and when $x$ is regular, semisimple, the $S_N$-Frobenius character of the Hessenberg variety computes the chromatic polynomial of a corresponding unit interval graph. By interpreting this action as the monodromy action of a family, Brosnan and Chow show that this interpretation proves the unimodality of the coefficients of the chromatic polynomial of such graphs \cite{BrosnanChow}.

The generalized Hessenberg varieties $\cH(x,\fm)$ fit into natural families which endow the cohomology of $\cH(x,\fm)$ with an $S_d$-action for regular, semisimple $x$. Following techniques of \cite{BrosnanChow}, we show the Frobenius character of this action is determined by the Poincar\'e polynomials of $\cH(x,\fm)$ for regular, but not necessarily semisimple $x$, and following \cite{julianna}, these Poincar\'e polynomials are computed in terms of combinatorial functions by studying an affine paving of $\cH(x,\fm)$ (see Propositions \ref{prop: affine paving of Hessenberg variety} and Proposition \ref{prop: combinatorial cell dim for hessenberg}).

In fact, we show that the Poincar\'e polynomial can be expressed in terms of natural generalizations of the chromatic polynomial of the associated unit interval graph. We generalize the digraph definition of the chromatic polynomial given in \cite[\S3.1]{BrosnanChow} to the setting of cyclic $n$-partite graphs (see \S\ref{sec: chromatic}). In particular, a generalized Hessenberg space produces a natural cyclic $n$-partite digraph and we call the associated chromatic polynomial $\Xi_\fm(Y,t)$. 

\begin{thm}
\label{thm: frob character of hessenberg_intro}
    ( = Theorem \ref{thm: frob character of hessenberg}) For $x\in \fg_{b}$ whose matrix representation \eqref{eqn: normal form for x} is regular, semisimple, all odd degree cohomology of $\cH(x,\fm)$ vanishes, and under the monodromy action, 
    \begin{align*}
        \sum_i t^i\Frob_{S_d}\left(H^{2i}(\cH(x,\fm))\right)=\Xi_{\fm}(Y,t).
    \end{align*}
\end{thm}

We show, using Theorem \ref{thm: frob character of hessenberg_intro}, that $\Xi_\fm(Y,t)$ is a symmetric function in $Y$, and that moreover, it is palindromic in $t$.

\subsection{Homogeneous affine Springer fibers and the Hessenberg paving}

Let $K = \C(\!(t)\!)$ be the field of Laurent series with $\cO\subset K$ the ring of power series in $t$. Let $\Fl_\gamma$ denote the Iwahori affine Springer fiber at an element $\gamma\in \mathfrak{gl}_N(K)$. The symmetric group $S_N$ acts on the Borel-Moore homology of $\Fl_\gamma$ through the Springer action.

An element $\gamma\in \mathfrak{gl}_N(K)$ is called \emph{homogeneous} of slope $c = M/N$ if $s^{-M}\cdot\gamma(s^Nt)$ is conjugate to $\gamma$ for all $s\in \C^\times$. For fixed $M,N$, it is known that all regular, semisimple homogeneous $\gamma\in \mathfrak{gl}_N(K)$ of slope $c = M/N$ have the same Borel-Moore homology  \cite{oy}. In the case when $M,N$ are coprime, the Frobenius character of the Borel-Moore homology of the affine Springer fiber of slope $c=M/N$ was computed in \cite{hikita} (in the setting of $\SL_N$). Computing this Borel-Moore homology explicitly in terms of the combinatorics of \cite{wilson} is one of the main results of this paper.

We fix a regular semisimple homogeneous element $\gamma$ of slope $c = M/N$. It was shown in \cite{GKM} that the affine Springer fiber $\Fl_\gamma$ admits a paving by generalized Hessenberg varieties. We show that the Hessenberg cells in $\Fl_\gamma$ are indexed by dominant labeled $(M,N)$ invariant sets $(\Delta,f)$ with $f$ injective. In character computations, we consider the \emph{positive part} $\Fl_\gamma^+\subset \Fl_\gamma$, which is the closed subset of Hessenberg cells labeled by positive invariant sets $\Delta\subset \Z_{>0}$. The area statistic gives a filtration on $H_{*}^{BM}(\Fl_\gamma^+)$ which we call the \emph{Hikita filtration} whose filtered quotients are finite dimensional. The Hikita filtration is conjecturally related to the perverse filtration up to an involution and reindexing; see \cite[\S8]{KT}.

We note that when $d=\gcd(M,N)$ is greater than 1, the Borel-Moore homology carries an additional structure:  In \cite{oy}, the authors show there is a family of regular semisimple homogeneous affine Springer fibers of fixed slope whose monodromy gives an action of $S_d$ on $H_{*}^{BM}(\Fl_\gamma)$. Hence, one can compute the bi-graded Frobenius character $\Frob_{S_N\times S_d}H_{*}^{BM}(\Fl_\gamma^+)$. Moreover, this is compatible with the $S_d$-actions on the generalized Hessenberg varieties in the paving of $\Fl_\gamma^+$. This was used by Carlsson and Mellit to compute $\Frob_{S_N\times S_N}H_{*}^{BM}(\Fl_\gamma^+)$ in terms of the $\nabla$-operator in the case when $M=kN$ is an integer multiple of $N$ \cite{CM-hessenberg}.

Our main result for affine Springer fibers is the following.

\begin{thm}
    \label{thm: ASF main_intro}
    ( = Theorem \ref{thm: ASF main}) Let $\gamma$ be a regular, semisimple, slope $c = M/N$ homogeneous element in $\mathfrak{gl}_N(K)$ and $d = \gcd(M,N)$. The $(q,t)$-bigraded $S_N\times S_d$ Frobenius character of $H^{BM}_*(\Fl_\gamma^+)$ is
    \begin{equation}
    \label{eqn:main_ASF_intro}
\Frob_{S_N\times S_d}\big(\mathrm{gr}_F H^{BM}_*(\Fl_\gamma^+) \big) = \sum_{\mu\vdash N}\sum_{(\Delta,f)}q^{\area(\Delta)}t^{\codinv(\Delta,f)}\frac{\Xi_{\fm_{(\Delta,f)}}(Y,t)}{\Aut_t(\alpha_{(\Delta,f)})}m_\mu(X).
    \end{equation}
    In particular, if we pair \eqref{eqn:main_ASF_intro} with $e_1^d[Y]$, then we recover $\omega_X L_{M,N}^{\mathrm{Wil}}(X,q,t^{-1})$, where $L_{M,N}^{\mathrm{Wil}}(X,q,t)$ is Wilson's symmetric function lift of torus link homology (see \S\ref{sec: combo}).
\end{thm}

\subsection{Minimal support representations of the rational Cherednik algebra}
For $n>0$, we let $\hfrak_n$ denote the $(n-1)$-dimensional standard representation of $S_n$ and $\Dcal(\hfrak_n)$ the ring of differential operators on $\hfrak_n$. Introduced in \cite{etingof-ginzburg}, the rational Cherednik algebra $\Hrm_c(S_n,\hfrak_n)$ is a deformation of $\Dcal(\hfrak_n)\rtimes S_n$ at a complex parameter $c$.
\cite{ma2} considers the unique finite-dimensional simple module $\Lrm_{m/n}$ of $\Hrm_{m/n}(S_n,\hfrak_n)$ when $m>0$ is coprime to $n$. It is shown that the bigraded Frobenius character of the associated graded of $\Lrm_{m/n}$ is given by the distinguished element $Q_{m,n}$ of the elliptic Hall algebra acting on $1$ in the polynomial representation. The computation is via the cuspidal Hodge module $\Nbb_{m,n}$ on $\mathfrak{sl}_n$. The associated graded of $\Nbb_{m,n}$ naturally gives rise to a coherent sheaf on the commuting variety and hence a coherent sheaf on the Hilbert scheme of points on the plane.

In the second part of the paper, we run the same machine as \cite{ma2} to obtain parallel results in the non-coprime case. For $\Hrm_{M/N}(S_N,\hfrak_N)$ for $\gcd(M,N)=d$, the role of the finite-dimensional simple modules is played by simple modules of minimal support, which are in bijection with partitions of $d$. For $\tau\vdash d$, let $\Lrm_{M/N}^{(N)}(n\tau)$ be the corresponding minimally supported module. As shown in \cite{egl}, the graded character of $\Lrm_{M/N}^{(N)}(n\tau)$ encodes the $\tau$-colored HOMFLY-PT polynomials of the $(m,n)$-torus knot. 

Let $V_\tau$ be the simple $S_d$-module associated to $\tau$, and let $\tau^t$ denote the transpose partition. Define
\begin{align*}
    \Lrm_{M,N}:=\bigoplus_{\tau\vdash d}V_{\tau^t}\otimes\Lrm_{M/N}^{(N)}(n\tau).
\end{align*}
By \cite{cee}, $\Lrm_{M,N}$ can be realized in terms of the multiplicity space of $(\C^M)^{\otimes N}$ in the $d$-th convolution power $\Nbb_{N,M}$ of the $\GL_m$-enhancement of $\Nbb_{n,m}$, and hence inherits the normalized Hodge filtration $F$. On $\Nbb_{N,M}$ one has a Springer action from \cite{BLM}, which descends to the $S_d$-action on $\Lrm_{M,N}$. Expressing the associated graded of $\Nbb_{M,N}$ in terms of the pushforward of a DG module as in \cite{isospec}, we compute the left hand side of the following theorem, which is then matched with the right hand side via shuffle algebra techniques of \cite{negutshuffle}, \cite{negutwalgebra}, \cite{negut-integral}, \cite{BHMPS},  \cite{gn_commuting_idempotent} and others.
\begin{thm}( = Theorem \ref{thm:Sd-joint-character})
		With respect to the Euler and filtration gradings, 
		\begin{align}
        \label{eqn:intro ch}
			\ch_{S_N\times S_d\times A}(\gr^F\Lrm_{M,N})(X,Y;q,h)
			=q^{(d-1)/2}(1-q)
			\left.\omega_X\bigl(h_d[(t-1)X^{m,n}Y]\cdot1\bigr)\right|_{t=h/q}.
		\end{align}
	\end{thm}
We also obtain a coherent sheaf $\Fcal_{M,N}$ on the Hilbert scheme of points on $\C^2$ whose equivariant K-theory class is represented by the right hand side of \eqref{eqn:intro ch}. We have the cohomology vanishing property of the tensor product of $\Fcal_{M,N}$ with hook part of the dual Procesi bundle, and we propose a vanishing conjecture for the full tensor product. 
\begin{rem}
In \cite{gko}, a quasi-coherent sheaf on $\Hilb^n(\C\times \C^*)$ is constructed for every affine Grassmannian Springer sheaf at a semisimple element. Such a sheaf is conjectured to be coherent. Their construction uses the $\Z$-algebra construction of \cite{gs}. Corollary 7.4 of \cite{ma2} generalizes to the non-coprime case and shows that our sheaf $\Fcal_{M,N}$ coincides with the image of $\Lrm_{M,N}^{S_N}$ under the Gordon-Stafford functor. As a result, our coherent sheaf $\Fcal_{M,N}$ should coincide with the sheaf in \cite{gko}. 
\end{rem}
\subsection{A potential future generalization} 
A natural extension is to consider mixed slope affine Springer fibers. Although there is no corresponding Cherednik algebra modules or character $D$-modules, one can still take parabolic induction of coherent sheaves on commuting varieties of mixed slopes. The relevant elliptic Hall algebra element is product of  $Q_{m,n}$'s for different $m,n$. The connection between the coherent sheaf and the elliptic Hall algebra element is automatic. However, GKM theory does not provide a cell decomposition of such affine Springer fibers, and there is currently no such shuffle theorem. 

\subsection{Organization of the paper}

This paper is split into two primary parts:  In part 1, we compute the Borel-Moore homology of regular semisimple homogeneous affine Springer fibers, and so prove Theorem \ref{thm: ASF main_intro}. To do so, we first study the geometry of generalized Hessenberg spaces by first producing an affine paving in \S\ref{sec: affine paving of hessenberg} and then studying families of generalized Hessenberg varieties in \S\ref{sec: families of hessenbergs}. In \S\ref{sec: combo}, we recall the combinatorics of labeled invariant sets and prove that Wilson's definitions \cite{wilson} differ from those appearing in this paper by the canonical involution on symmetric functions. Finally in \S\ref{sec: ASF}, we consider the Borel-Moore homology of affine Springer fibers by studying Goresky-Kottwitz-MacPherson's Hessenberg paving in terms of the combinatorics of $M,N$-invariant sets and the geometry of generalized Hessenberg spaces.

In part 2, we adapt the methods of \cite{ma2} to construct lifts of the minimal support representations of the rational Cherednik algebra to character $D$-modules, whose associated graded defines our sheaves $\Fcal_{M,N}$ on the Hilbert scheme. In \S\ref{sec:coherent-characters} we compute the bigraded character of the associated sheaves on the commuting variety and on the Hilbert scheme using methods similar to \emph{loc. cit.}; in particular, we also collect partial results on vanishing of higher cohomology of $\Fcal_{M,N}$. In \S\ref{sec:Sd-refinement}, we study the additional $S_d$-action and compute the joint character. The comparison with the constructible side is stated in Theorem \ref{thm:character-comparison}.

\subsection{Acknowledgements}We thank Donghyun Kim and Jaeseong Oh for helpful discussions and their proof of the shuffle theorem, Minh-T\^{a}m Trinh for inspiring discussions, Eugene Gorsky and Nicolle Gonzalez for bringing \cite{KO1} to our attention and Eugene Gorsky and Sean Griffin for helpful communications. 
The third-named author was partially supported by the NSF grant DMS-2503504. 

ChatGPT was used to assist the proofs of Proposition \ref{prop: char triv part} and Theorem \ref{thm:Sd char} and audit numerous drafts. Besides cosmetic editing, it corrected errors in the proof of Proposition \ref{prop: combinatorial poincare duality} and statement of Lemma \ref{lem: hikita vs wilson labeling convention}, proposed an alternate proof strategy for Proposition \ref{prop: affine paving of Hessenberg variety}, and fixed normalizations throughout part 2, such as the factor $q^{r/2}$ in Theorem \ref{thm:Sd char}, the factor $q^{d-N}$ in Theorem \ref{thm:character-comparison}, etc. 
\part{Constructible side}

\section{Cyclic Hessenberg varieties}
We will introduce some generalized Hessenberg varieties which will appear in the description of our affine Springer fibers using \cite{GKM}. It is an interesting feature that the Hessenberg varieties discussed here are natural generalizations of those arising from symmetric varieties.
\subsection{Hessenbergs of finite order automorphisms}
Let $G$ be a reductive group with an algebraic automorphism $\theta$ of finite order $n$. We also denote by $\theta$ the induced action on $\fg$. We have a decomposition
\[
\fg = \bigoplus_{\xi\in \mu_n} \fg_\xi
\]
where $\mu_n$ denotes the set of $n$-th roots of unity and $\fg_\xi$ is the eigenspace of the root of unity $\xi$ in $\fg$. Let $H = G^\theta$ be the fixed point subgroup of $G$. The adjoint action of $H$ acts naturally on the eigenspaces $\fg_\xi$. Fix a primitive $n$-th root of unity $\zeta$. We will be interested in the following varieties associated to the representation of $H$ on $\fg_\xi$.

\begin{defn}
    Fix a maximal abelian Lie algebra $\fa\subset \fg_\xi$ and a Borel $B\subset H$. A vector subspace $\fm\subset \fg_\xi$ is called a \emph{(generalized) Hessenberg space} if $\fm$ is $B$ stable and contains $\fa$.

    Let $x\in \fg_\xi$, and let $\fm\subset \fg_\xi$ be a Hessenberg space in $\fg_\xi$. The Hessenberg variety $\cH(x,\fm)$ is the subvariety of the flag variety $H/B$ given by
    \[
\cH(x,\fm) = \{ g\in H/B\mid g^{-1}\cdot x\in \fm \}.
    \]

    If $P\supset B$ is a standard parabolic of $H$, then we denote $\cH_P(x,\fm)$ for the image of $\cH(x,\fm)$ in $H/P$.
\end{defn}

The varieties $\cH(x,\fm)$ fall into the framework of generalized Hessenberg varieties arising in \cite{GKM}, and when $\theta$ is an algebraic involution, \cite{tvx} studies this particular type of Hessenberg variety for the example of the symmetric space $\SL_n/\SO_n$. While the authors believe the class of Hessenbergs defined above are a natural family for further study, in this paper we restrict to the special case below. We note that this case is simpler than that of \cite{tvx} as the automorphism $\theta$ below is inner.

\begin{defn}
\label{def: theta and g zeta}
    Fix a primitive $n$-th root of unity $\zeta\in \mu_n$. From here on, throughout Part I, we denote by $G = \GL_{N}$ for $N = nd$ and take
    \[
    \theta(g) = J^{-1}gJ\quad \text{ for }J = \begin{pmatrix}
        I_d & & & \\
         & \zeta I_d & & \\
         & & \ddots & \\
         & & & \zeta^{n-1}I_d
    \end{pmatrix}
    \]
    In this case, $H = \GL_d^n$ is the block diagonal Levi subgroup of $\GL_N$, which acts on the subspace $\fg_{\xi}\subset \fg$ for $\xi = \zeta^{b}$.
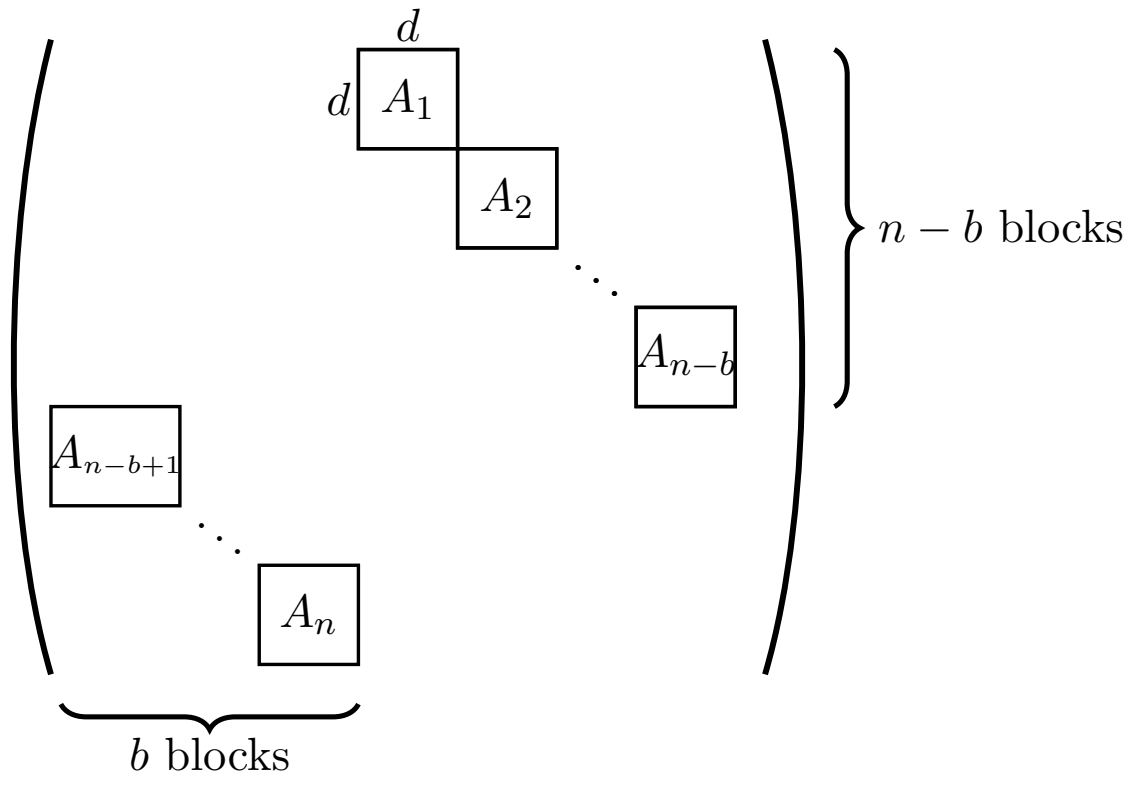
    \begin{equation}
    \label{eqn: normal form for x}
    \scalebox{.8}{\begin{tikzpicture}[scale=.82]
\draw[thick] 
  (3.8,2.4) rectangle (2.8,3.4);

\draw[thick] 
  (2,4) rectangle (1,5);

\node at (2.4,3.8) {$\ddots$};
  
\draw[thick] 
  (1,5) rectangle (0,6);

\draw[thick] 
  (-1.8,2.4) rectangle (-3.1,1.4);

\node at (-1.4,1.2) {$\ddots$};

\draw[thick] 
  (-1,0.8) rectangle (-0,-0.2);

\draw[line width=1.4pt]
  (-3.1,6.1)
  .. controls (-3.6,4.2) and (-3.6,1.4) ..
  (-3.1,-0.3);

\draw[line width=1.4pt]
  (4.1,6.1)
  .. controls (4.6,4.2) and (4.6,1.4) ..
  (4.1,-0.3);

\node at (-.2,5.5) {$d$};

\node at (.5,6.25) {$d$};

\draw[thick] 
  (3.8,2.4) rectangle (2.8,3.4);
\node at (3.3,2.9) {$A_{n-b}$};

\draw[thick] 
  (2,4) rectangle (1,5);
\node at (1.5,4.5) {$A_{2}$};

\node at (2.4,3.8) {$\ddots$};
  
\draw[thick] 
  (1,5) rectangle (0,6);
\node at (0.5,5.5) {$A_1$};

\node at (-0.5,0.3) {$A_n$};
\node at (-2.45,1.9) {$A${\footnotesize{${}_{n-b+1}$}}};

\draw[
  decorate,
  decoration={brace,amplitude=6pt},
  line width=1.2pt
]
  (4.8,6.0) -- (4.8,2.4);

\node at (6.48,4.2) {$n-b$ blocks};
\draw[
  decorate,
  decoration={brace,mirror,amplitude=6pt},
  line width=1.2pt
]
  (-3.0,-0.6) -- (0,-0.6);
  
\node at (-1.5,-1.1) {$b$ blocks};
\end{tikzpicture}}
    \end{equation}
\end{defn}

We say that an element $x\in \fg_\xi$ is \emph{regular} if its centralizer $C_H(x)$ in $H$ is of minimal dimension. We denote the set of all regular elements in $\fg_\xi$ by $\fg_\xi^\reg$. We will assume throughout that $b$ is coprime to $n$.

\begin{lem}
\label{lem: regularity agrees}
    An element $x\in \fg_\xi$ is regular with respect to the $H$ action if and only if it is regular as an element of $\fg$ (with respect to the action of $G$). That is, we have
    \[
    \fg_\xi^\reg = \fg_\xi\cap \fg^\reg
    \]
\end{lem}
\begin{proof}
    This follows from \cite[Proposition 11]{vinberg}.
\end{proof}

We will be specifically interested in elements $x\in \fg_\xi$ which are regular and which have no 0 eigenvalue. The following says that, to treat such cases, it suffices to consider cases where $x$ has identical blocks.

\begin{lem}
\label{lem: normal form of x}
    If $x\in \fg_\xi^\reg$ has all eigenvalues nonzero, then there exists an element of $H$ which conjugates $x$ to an element of the form \eqref{eqn: normal form for x} with all $A_i=A$ for $A$ being a $d\times d$ regular matrix in standard Jordan normal form with all eigenvalues nonzero and $A^n$ regular. Moreover, $A$ is unique up to reordering of Jordan blocks and multiplication of eigenvalues by $n$-th roots of unity.

    In this case, we say that $x$ is of type $\lambda\vdash d$ if the matrix $A$ has Jordan type $\lambda\vdash d$.
\end{lem}
\begin{proof}
    Consider $x$ in the form \eqref{eqn: normal form for x}. In particular, we label the blocks of $x$ by $A_i$ according to the picture \eqref{eqn: normal form for x}. We write subscripts cyclically, so that $A_{n+i}=A_i$. The characteristic polynomial of $x$, $\chi_x(t)$ can be computed in terms of the characteristic polynomial of the product $A_1 A_{b+1} A_{2b+1}\cdots A_{(n-1)b+1}$ by
    \[
    \chi_x(t) = \chi_{A_1 A_{b+1} A_{2b+1}\cdots A_{(n-1)b+1}}(t^n).
    \]
    Since $x$ is regular with no zero eigenvalues, we in particular have that $A_i$ is invertible for all $i$, and in particular, there exists an invertible matrix $A$ such that $A^n = A_1 A_{b+1} A_{2b+1}\cdots A_{(n-1)b+1}$. By conjugating by a block diagonal matrix $g = (g_1,\dots,g_{n-1},g_n)\in H$, with $g_{1+ib} = A^{-i}A_1A_{1+b}\cdots A_{1+(i-1)b}$ for all $0\leq i<n$, with indices understood modulo $n$, we may assume that $x$ has $A_i = A$ for all $i$. Then, diagonal conjugation (by elements $g = (g_1,g_1,\cdots,g_1)$) allows us to put $A$ into Jordan form, unique up to reordering of factors and multiplying eigenvalues by $n$-th roots of unity. Note that the characteristic polynomial of this matrix is $\chi_{A^n}(t^n)$, so in particular, regularity of $x$ implies regularity of $A^n$.
\end{proof}

We note that there is a perhaps more intuitive way of representing elements $x\in \fg_\xi^\reg$ where $A_i = I_d$ for all $i>1$ and $A_1$ is in Jordan form. The advantage of our choice here is that this choice will behave nicely in relation to the Bruhat decomposition of $H$, as we will see in \S\ref{sec: affine paving of hessenberg}.

Let $\ft_{\Delta}$ denote the set of diagonal matrices 
\[
\{\mathrm{diag}(\delta,\delta,\cdots,\delta)\mid \delta \text{ diagonal in }\fgl_d\}\subset \Lie(H).
\]
We fix $\fa = \fc_{\fg_\xi}(\ft_\Delta)$ the centralizer of this diagonal torus in $\fg_\xi$. With respect to $\fa$ and the standard upper triangular Borel of $H$, a Hessenberg space $\fm\subset \fg_\xi$ is determined by a tuple of Hessenberg functions $\mathrm{m}_i$; explicitly, $\fm$ corresponds to an $n$-tuple of Hessenberg functions $(\mathrm{m}^{(1)},\dots, \mathrm{m}^{(n)})$ if $\fm$ consists of all matrices of the form \eqref{eqn: normal form for x} with $A_\ell$ having all entries $(i,j)$ with $i>m^{(\ell)}_{j}$ equal to zero.

\subsection{Affine paving of $\cH(x,\fm)$}
\label{sec: affine paving of hessenberg}

In this section, we construct an affine paving for elements $x\in \fg_\xi$ in the standard form of Lemma \ref{lem: normal form of x} by intersecting with Bruhat cells for $H$. We find a classification of which cells are nonempty (Proposition \ref{prop: nonempty cells for hessenberg}) and a formula for the dimension of nonempty cell (Propositions \ref{prop: dimension of cells for hessenberg} and \ref{prop: combinatorial cell dim for hessenberg}).

Let $B: = B_H$ denote the standard upper triangular Borel in $H = \GL_d^n$. Let $U\subset B$ denote the unipotent radical of $B$, $T$ the diagonal maximal torus in $H$ (and so also in $G$), and $W_H$ the Weyl group of $H$ with respect to $T$. We fix representatives $\dot w\in N_G(T)$ for $w\in W$. We denote by $\Phi$ the set of roots of $G$ with respect to $T$, $\Phi^+$ the positive roots for a Borel of $G$ extending $B$ (i.e. $B_G\cap H = B$), and $\Sigma$ the set of simple roots in $\Phi^+$. We adopt similar notation for $\GL_d$ with respect to the standard Borel, with associated sets of roots $\Sigma_d\subset \Phi^+_d\subset \Phi_d$. For $1\leq i,j\leq N$, we denote by $E_{i,j}\in \fg$ the matrix with a 1 in the $(i,j)$ entry and zeros elsewhere, and we write $E_{\alpha_{i,j}} = E_{i,j}$ whenever $\alpha_{i,j} = x_i - x_j$ is the corresponding root of $\GL_N$ with respect to the standard diagonal maximal torus.

Let $W_H = S_d^n$ denote the Weyl group of $H$. Recall the standard Bruhat decomposition
\[
H/B = \bigsqcup_{w\in W_H}B\dot w B/B
\]
where each cell $B\dot w B/B$ is an affine space whose dimension is equal to the length $\ell(w)$.
\begin{defn}
    Fix $x$ to be of the form in Lemma \ref{lem: normal form of x}. For each $w\in W_H$, we let $\cH^w(x,\fm)$ denote the intersection
    \[
    \cH^w(x,\fm) = \cH(x,\fm)\cap B\dot w B/B
    \]
\end{defn}

We note that, in what follows, it is important that $x$ and $B$ be chosen to be compatible, as we have done here.

\begin{lem}
\label{lem: unipotent action on p}
    Let $x\in \fg_\xi$ be of the form described in Lemma \ref{lem: normal form of x}, where we express $A$ in Jordan normal form, with nonzero entries only on the diagonal and upper super-diagonal. Let $\Sigma_x$ denote the set of simple roots $\alpha_{i,i+1}$  in $\Sigma_d$ such that $A$ contains a nonzero entry in position $(i,i+1)$, and for any root $\alpha_{i,j}\in \Phi_d$, let $\Phi_{i,j}$ denote the set of all roots $\alpha_{i',j'}$ of $G = \GL_N$ such that $(i',j')\equiv (i,j)\mod d$ and $\fg_{\alpha_{i',j'}}\subset \fg_{\zeta^b}$. Then, we have 
    \[
    U\cdot x = x + \mathcal{U} +\mathcal{V},
    \]
    where
    \[
    \mathcal{U} = \bigoplus_{i\colon \;\alpha_{i,i+1}\in \Sigma_x}\left\{ \left.\sum_{\alpha\in \Phi_{i,i+1}} c_\alpha E_\alpha \; \right| \; \sum_{\alpha\in \Phi_{i,i+1}}c_\alpha = 0 \right\},\quad 
    \mathcal{V} = \bigoplus_{\alpha_{i,j}\in \Phi^+_d \setminus \Sigma_x}\bigoplus_{\alpha\in \Phi_{i,j}} \fg_\alpha
    \]
\end{lem}
\begin{proof}
For $u=(u_p)\in U$, put $B_p=u_p^{-1}Au_{p+b}$. Each $B_p$ is upper triangular with the same diagonal as $A$, and
\[
(B_p)_{i,i+1} = A_{i,i+1}+A_{i,i}(u_{p+b})_{i,i+1}-A_{i+1,i+1}(u_p)_{i,i+1}.
\]
If $\alpha_{i,i+1}\in\Sigma_x$, then $A_{i,i}=A_{i+1,i+1}$, so summing over $p$ gives $\sum_p(B_p)_{i,i+1}=nA_{i,i+1}$. Thus $U\cdot x\subseteq x+\mathcal U+\mathcal V$.

Both sides have dimension $\dim(U)-d+r$, where $r$ is the
number of distinct Jordan factors in $A$. Indeed, the stabilizer of $x$ in $U$ is isomorphic to $Z_{U_d}(A^n)$, of dimension $d-r$, while $\mathcal U+\mathcal V$ has codimension $|\Sigma_x|=d-r$ in the space of strictly upper triangular block tuples.

By the Kostant--Rosenlicht theorem, $U\cdot x$ is closed.
The inclusion above is therefore an equality, since
$x+\mathcal U+\mathcal V$ is irreducible and both sides
have the same dimension.
\end{proof}

\begin{prop}
    \label{prop: nonempty cells for hessenberg}
    A cell $\cH^w(x,\fm)$ is nonempty if and only if $w\cdot \fm$ contains the diagonal and for all $\alpha_{i,i+1}\in \Sigma_x$, there exists some $j$ such that $w\cdot \fm$ contains $E_{i+jd+1,i+(b+j)d}$.
\end{prop}
\begin{proof}
    By definition, $\cH^w(x,\fm)$ is nonempty if there exists $u\in U$ such that $u^{-1}\cdot x\in w\cdot \fm$. By Lemma \ref{lem: unipotent action on p}, this is equivalent to the condition
    \[
    (x+\Ucal+\Vcal)\cap w\cdot \fm\neq \emptyset.
    \]
    In particular, as $w\cdot \fm$ and $\Vcal$ are sums of root spaces and the roots in $\Vcal$ are disjoint from the nonzero entries of $x$, this is equivalent to $x+\Ucal\cap w\cdot \fm\neq\emptyset$. The set $x+\Ucal$ has nonzero entries on the diagonal entries of each block, and for any $\alpha_{i,i+1}\in \Sigma_x$, the sum of terms in the roots spaces in $\Phi_{i,i+1}$ is nonzero. It is easy to see that this is nonempty if and only if $w\cdot \fm$ includes all diagonal terms and for every $\alpha_{i,i+1}\in \Sigma_x$, there exists some $\alpha\in \Phi_{i,i+1}$ such that $w\cdot \fm$ is not zero in the root space $\fg_\alpha$.
\end{proof}

\begin{prop}
    \label{prop: dimension of cells for hessenberg}
    Suppose that $w$ satisfies the conditions of Proposition \ref{prop: nonempty cells for hessenberg}, which is to say $\cH^w(x,\fm)$ is nonempty. Then $\cH^w(x,\fm)$ is smooth and 
    \begin{equation}
    \label{eqn: first dim formula for Hessenberg cells}
    \dim \cH^w(x,\fm) = \dim\left( U/U\cap w\cdot U \right) - \dim\left(U\cdot x/U\cdot x\cap w\cdot \fm \right)
    \end{equation}
\end{prop}
\begin{proof} 
    We follow a very similar argument to \cite[Proposition 3.7]{precup}. Consider the set
    \[
    U(x,w\cdot \fm) = \{u\in U\mid u^{-1}\cdot x\in w\cdot \fm\}.
    \]
    There is a surjective map
    \[
    U(x,w\cdot \fm)\to \cH^w(x,\fm),\quad u\mapsto u\dot{w}
    \]
    This morphism is surjective with fibers given by cosets of $U\cap w\cdot U$ in $U$. In particular, it suffices to show that $U(x,w\cdot \fm)$ is smooth and 
    \[
    \dim U(x,w\cdot \fm) = \dim U - \dim\left(U\cdot x/U\cdot x\cap w\cdot \fm \right)
    \]
    We consider the action map
    \[
    U\to U\cdot x, \quad u\mapsto u^{-1}\cdot x
    \]
    This map is smooth of relative dimension $\dim Z_U(x)$. We form the Cartesian diagram
    \[
    \xymatrix{
    U\times_{U\cdot x}(U\cdot x\cap w\cdot \fm)\ar[r]\ar[d] & U\cdot x\cap w\cdot \fm \ar[d] \\
    U\ar[r] & U\cdot x
    }
    \]
    We note that by Lemma \ref{lem: unipotent action on p}, the intersection $U\cdot x\cap w\cdot \fm$ is an affine translate of a linear subspace of the vector space $w\cdot \fm$, and so is smooth. Therefore, the fiber product is smooth of dimension  $\dim Z_U(x)+\dim(U\cdot x\cap w\cdot \fm)$. Furthermore, we can identify the fiber product with
    \[
    U\times_{U\cdot x}(U\cdot x\cap w\cdot \fm) = U(x,w\cdot \fm).
    \]
    Hence, we conclude that 
    \begin{align*}
    \dim U(x,w\cdot \fm) &= \dim Z_U(x)+\dim(U\cdot x\cap w\cdot \fm) \\
    &= \dim U - \dim U\cdot x + \dim(U\cdot x\cap w\cdot \fm) \\
    & = \dim U - \dim\left( U\cdot x/U\cdot x\cap w\cdot \fm \right) \qedhere
    \end{align*}
\end{proof}

The dimension formula of Proposition \ref{prop: dimension of cells for hessenberg} is strongly reminiscent of \cite[Corollary 5.5]{precup}. We now give a more combinatorial description of this cell dimension similar in spirit to that of Tymoczko's cell dimension in type $A$ \cite[Theorem 4.3]{julianna} (see also \cite[Theorems 33 and 35]{BrosnanChow}).

\begin{prop}
    \label{prop: combinatorial cell dim for hessenberg}
    If $x\in \fg_\xi^\reg$ has no zero eigenvalues and is of type $\lambda\vdash d$ and if $\cH^w(x,\fm)$ is nonempty, then $\dim \cH^w(x,\fm)$ is computed in the following combinatorial way:  First, consider $w^{-1} = (w_1^{-1},\cdots, w_{n}^{-1})$ where $w_i^{-1}$ is an element of $S_d$ considered as a function on $\{1,\dots,d\}$. Then, the dimension of $\cH^w(x,\fm)$ is the number of triples $(i,j,\ell)$ with $i<j$ such that
    \begin{equation}
    \label{eqn: combinatorial dim for hessenberg cells}
    w_\ell^{-1}(j)< w_\ell^{-1}(i)\leq m_{w^{-1}_{\ell+b}(j)}^{(\ell)}
    \end{equation}
    where we consider the subscript $\ell+b$ mod $n$.
\end{prop}
\begin{proof}
    We compute the dimension of cells using the description of Proposition \ref{prop: dimension of cells for hessenberg}. Let $\fu_\xi$ denote 
    \[
    \fu_\xi := \bigoplus_{\substack{\ell = 1,\dots, n \\ i<j}} \fg_{\alpha_{i+(\ell-1)d,j+(b+\ell-1)d}}\subset \fg_\xi.
    \]
    where we view subscripts mod $N$. By the description of $U\cdot x$ in Lemma \ref{lem: unipotent action on p}, when $\cH^w(x,\fm)$ is nonempty we have
    \[
\dim\left( U\cdot x/U\cdot x\cap w\cdot \fm \right) = \dim(\mathfrak{u}_\xi/\mathfrak{u}_\xi\cap w\cdot \fm)
    \]
    Likewise, $\dim(U/U\cap w\cdot U) = \dim(\mathfrak{u}/\mathfrak{u}\cap w\cdot \mathfrak{u})$, so we can rewrite \eqref{eqn: first dim formula for Hessenberg cells} as
    \[
    \dim(\fu_\xi\cap w\cdot \fm) - \dim(\mathfrak{u}\cap w\cdot \mathfrak{u})
    \]
    Since $\mathfrak{u} = \bigoplus_{\ell=1}^n\bigoplus_{i<j}\fg_{(\ell-1)d+i,(\ell-1)d+j}$, we can rewrite the above in terms of triples:
    \[
\#\left\{(i,j,\ell)\mid i<j,\; w_\ell^{-1}(i)\leq m_{w_{\ell+b}^{-1}(j)}\right\} - \#\left\{ (i,j,\ell)\mid i<j,\; w_\ell^{-1}(i)<w_\ell^{-1}(j) \right\}
    \]
    which returns \eqref{eqn: combinatorial dim for hessenberg cells}.
\end{proof}

\begin{prop}
\label{prop: affine paving of Hessenberg variety}
    For any $w$, the intersection $\cH^w(x,\fm)$ is an affine space whenever it is nonempty. Therefore, $\cH(x,\fm)=\sqcup_w\cH^w(x,\fm)$ forms an affine paving.
\end{prop}

\begin{proof}
The proof of Proposition \ref{prop: affine paving of Hessenberg variety} will be analogous to the argument in \cite{julianna}. The strategy will be to construct a family of affine space bundles interpolating between a point and the full Hessenberg cell $\cH_w(x,\fm)$ by adding conditions one row at a time. First, we need some auxiliary notation (which will not be used elsewhere in the paper). We may assume throughout that $\cH_w(x,\fm)$ is
nonempty.

Recall $U\subset H$ denotes the standard unipotent subgroup of $H$ and $U^-$ its opposite. Put $U_w=U\cap w\cdot U^-$. Recall that we have an identification $BwB/B=U_wwB/B$, and that furthermore, any element $x\in U_wwB/B$ determines a \emph{unique} element $u\in U_w$ so that $x=uwB/B$.

For any $u\in U$, we can write $u=(u_1,\dots,u_n)$, where $u_j\in\GL_d$ is a $d\times d$ upper triangular unipotent matrix for all $j$. Let $U^{\geq i}\subset U$ denote the subset of elements $u$ such that, for all $j$, $I_d-u_j$ only has nonzero entries in rows $\geq i$. Furthermore, for each $i$, we denote by $\mathrm{pr}_i\colon\fgl_d\to\C^{d(d-i+1)}$ the projection of a $d\times d$ matrix to its rows $i,\dots,d$, with $\mathrm{pr}_{d+1}=0$. We use the same notation for the corresponding blockwise projection on $\fg_\xi$. If $x$ is in the normal form of Lemma \ref{lem: normal form of x}, we put
\[
Z_i:=
\left\{
u\in U_w\cap U^{\geq i}
\ \middle|\
\big(\mathrm{pr}_i(u_p^{-1}Au_{p+b})\big)_{p=1}^n
\in\mathrm{pr}_i(w\cdot\fm)
\right\}.
\]
That is, $Z_i$ imposes the condition $u^{-1}xu\in w\cdot\fm$ for each $d\times d$ block in rows $\geq i$. In particular, $Z_{d+1}={\color{blue}\{1\}}$ while $Z_1\simeq\cH_w(x,\fm)$.

We can decompose $u=(u_1,\dots,u_n)\in U^{\geq i}$ by $u_p=v_p(1+R_p)$, where $v=(v_1,\dots,v_n)\in U^{\geq(i+1)}$ and $R_p$ is strictly upper triangular, with nonzero entries only in row $i$. It is easy to check that $R_p^2=0$ and so $u_p^{-1}=(1-R_p)v_p^{-1}$. Put $W_p=v_p^{-1}Av_{p+b}$. Since $W_p$ is upper triangular and $R_p,R_{p+b}$ are supported in row $i$ and columns $>i$, we have $R_pW_pR_{p+b}=0$. We have
\[
u_p^{-1}Au_{p+b} = W_p+W_pR_{p+b}-R_pW_p.
\]
Therefore, in rows $\geq(i+1)$, $u_p^{-1}Au_{p+b}$ and $W_p=v_p^{-1}Av_{p+b}$ agree. Hence, we have shown that the morphism $u\mapsto v$ defines a map
\[
\pi_i\colon Z_i\to Z_{i+1}.
\]

We wish to show that $\pi_i$ is an affine space bundle. Fix an integer $1\leq i\leq d$. Let $\cR$ denote the space of tuples $(R_1,\dots,R_n)$ such that $R_p$ is strictly upper triangular, with all entries zero except possibly those $(i,j)$ with $w_p^{-1}(j)\leq w_p^{-1}(i)$. We call triples $(p,i,j)$ where $R_p$ is allowed a nonzero $(i,j)$ entry \emph{free}. The projection
\[
U_w\cap U^{\geq i} \to U_w\cap U^{\geq(i+1)}, \qquad u\mapsto v
\]
has fiber $\cR$ at every point $v\in U_w\cap U^{\geq(i+1)}$.

On the other hand, we call a triple $(p,i,j)$ \emph{constrained} if the $(i,j)$ entry of the $p$-th block is required to vanish in $w\cdot\fm$. The fiber of $\pi_i$ can be understood in terms of constrained triples: for every constrained triple $(p,i,j)$, we have a relation
\[
(W_p)_{i,j}+A_{i,i}(R_{p+b})_{i,j}-(R_pW_p)_{i,j}=0.
\]
Let $\cT$ be a vector space with basis indexed by all constrained triples $(p,i,j)$ with $j>i$, and define
\[
\Psi\colon\cR\to\cT,\qquad \Psi(R_1,\dots,R_n) =\big(A_{i,i}(R_{p+b})_{i,j}-(R_pW_p)_{i,j}\big)_{(p,i,j)}\in\cT.
\]
$\Psi$ is a linear map, and the fiber of $\pi_i$ at $v$ consists of solutions $R\in\cR$ to the system of equations
\begin{equation}
\label{eqn: fiber of pi equations}
\Psi(R)=\big(-(W_p)_{i,j}\big)_{(p,i,j)}
\end{equation}
for $W_p=v_p^{-1}Av_{p+b}$ as before. Here and below, the coordinates of $\cT$ are indexed only by constrained triples with $j>i$. The conditions with $j\leq i$ are automatic under our nonemptiness assumption, as verified in part~(1) of the following lemma.

\begin{lem}
    \label{lem: constrained and free}
    Let $1\leq p\leq n$ and $1\leq i,j,j'\leq d$ be integers. We consider $p+b$ mod $n$.
    \begin{enumerate}
        \item If $j<i$ and $(p,i,j)$ is constrained, then the condition $(u_p^{-1}Au_{p+b})_{i,j}=0$ is automatic. If $(p,i,i)$ is constrained, then $Z_i$, and hence $\cH_w(x,\fm)$, are empty.
        \item Assume $\cH_w(x,\fm)$ is nonempty. If $(p,i,j)$ is constrained and $j>i$, then $(p,i,j)$ and $(p+b,i,j)$ are free.
        \item Assume $\cH_w(x,\fm)$ is nonempty. Suppose $v\in Z_{i+1}$ and use it to define $W_p=v_p^{-1}Av_{p+b}$. If $j>i$, $(p,i,j)$ is \emph{not} free and $(p,i,j')$ is constrained, then $(W_p)_{j,j'}=0$.
    \end{enumerate}
\end{lem}

\begin{proof}
    The Hessenberg space $\fm$ is determined by a tuple of Hessenberg functions $(m^{(1)},\dots,m^{(n)})$. The triple $(p,i,j)$ is constrained when $m_{w_{p+b}^{-1}(j)}^{(p)}<w_p^{-1}(i)$ and free when $j>i$ and $w_p^{-1}(j)\leq w_p^{-1}(i)$.

    Suppose $j<i$. Since $W_p$ is a product of upper triangular matrices, $(W_p)_{i,j}=0$ and $(R_pW_p)_{i,j}=0$. Therefore, the corresponding entry condition reads $0=0$ in this case. If we take $j=i$, the corresponding entry condition reads $0=-A_{i,i}$, which contradicts our assumption that $x$, and hence $A$, has no zero eigenvalue. Hence, we have proved (1).

    For (2), since $\cH_w(x,\fm)$ is nonempty, $(p,\ell,\ell)$ is not constrained for all $\ell$, and so
    \[
    w_p^{-1}(j)\leq m^{(p)}_{w_{p+b}^{-1}(j)}<w_p^{-1}(i)
    \]
    and $(p,i,j)$ is free. Similarly,
    \[
    m^{(p)}_{w_{p+b}^{-1}(j)}<w_p^{-1}(i)\leq m^{(p)}_{w_{p+b}^{-1}(i)},
    \]
    and so $w_{p+b}^{-1}(j)<w_{p+b}^{-1}(i)$ since $m^{(p)}$ is a nondecreasing function. Hence, $(p+b,i,j)$ is also free.

    Finally, if $(p,i,j)$ is not free, $i<j$, and $(p,i,j')$ is constrained, then we have
    \[
    w_p^{-1}(j)>w_p^{-1}(i)>m^{(p)}_{w_{p+b}^{-1}(j')}
    \]
    so that $(p,j,j')$ is constrained. Therefore, since $j>i$, the condition from the definition of $v\in Z_{i+1}$ implies that $(W_p)_{j,j'}=0$, proving (3).
\end{proof}

We claim now that the map $\Psi$ is surjective. Let $\xi\in\cT^*$ be a linear functional for which $\xi\circ\Psi=0$. We will show that $\xi=0$. We view $\xi$ as a tuple $\xi=(\xi_1,\dots,\xi_n)$, where $\xi_p$ is a $d$ dimensional column vector supported in rows $j$ with $j>i$ such that $(p,i,j)$ is constrained. We have
\[
\xi\circ\Psi(R)=\sum_p\big(A_{i,i}R_{p+b}\xi_p-R_pW_p\xi_p\big) = \sum_p R_p(A_{i,i}\xi_{p-b}-W_p\xi_p).
\]
So, $\xi\circ\Psi(R)=0$ if and only if $A_{i,i}\xi_{p-b}[j]=(W_p\xi_p)[j]$ for all $p$ and $j$ such that $(p,i,j)$ is free.

If $(p,i,j)$ is not free and $j>i$, then by Lemma \ref{lem: constrained and free}(2), $(p-b,i,j)$ and $(p,i,j)$ are not constrained. Therefore, $\xi_{p-b}[j]=0$. On the other hand,
\[
(W_p\xi_p)[j] = \sum_{j'}(W_p)_{j,j'}\xi_p[j']
\]
where when $\xi_p[j']\neq0$, $(p,i,j')$ is constrained so since $(p,i,j)$ is not free, $(W_p)_{j,j'}=0$ by Lemma \ref{lem: constrained and free}(3). Therefore, $(W_p\xi_p)[j]=0$ as well. We conclude that $A_{i,i}\xi_{p-b}[j]=(W_p\xi_p)[j]$ for all $p$ and all $j>i$, regardless of whether $(p,i,j)$ is free.

We will denote by an overline the projection to coordinates
$j>i$. For upper triangular matrices, an overline denotes the
lower-right submatrix indexed by $j>i$; this operation
commutes with multiplication and inversion. We have shown that
$A_{i,i}\ol{\xi}_{p-b}=\ol{W}_p\ol{\xi}_p$.
Composing then gives
\[
A_{i,i}^n\ol\xi_{p-b} = \ol{W_pW_{p+b}\cdots W_{p+(n-1)b}}
\cdot\ol\xi_{p-b} = \ol{v_p^{-1}A^n v_p}\cdot\ol\xi_{p-b}.
\]
Thus
\[
\ol\xi_{p-b}\in\ol v_p^{-1}\ker\big(\ol A^{\,n}-A_{i,i}^{\,n}I\big).
\]
Since $A^n$ is regular, it has a unique Jordan block with eigenvalue $A_{i,i}^n$. Note that if $J$ is a Jordan block with eigenvalue $a$, then
\[
J^n = a^nI+na^{n-1}N+\binom{n}{2}a^{n-2}N^2+\cdots,
\]
where $N$ is the standard principal nilpotent given by the sum of simple root vectors. Since $a\neq0$, the kernel of $J^n-a^nI$ is the kernel of $N$. Consequently, if the Jordan block of $A$ containing $i$ ends at $i$, then $\ker(\ol A^{\,n}-A_{i,i}^{\,n}I)=0$. Otherwise, $i$ and $i+1$ belong to the same Jordan block and
\[
\ker\big(\ol A^{\,n}-A_{i,i}^{\,n}I\big) = \C\cdot e_{i+1}.
\]
Moreover, $\ol v_p$ is upper triangular unipotent, so it fixes the first basis vector $e_{i+1}$ of the truncated space. It follows that either $\xi=0$, or $i$ and $i+1$ belong to the same Jordan block and $\xi_p=\lambda_pe_{i+1}$ for every $p$, with $\lambda_p\in\C$. If $\lambda_p\neq0$, then the equation $A_{i,i}\ol\xi_{p-b}=\ol W_p\ol\xi_p$ forces $\lambda_{p-b}=\lambda_p$ since $(W_p)_{i+1,i+1}=A_{i+1,i+1}=A_{i,i}$. Therefore, since $b$ is coprime to $n$, $\lambda_p\neq0$ for all $p$, and $(p,i,i+1)$ is constrained for all $p$. Since $i$ and $i+1$ are part of the same Jordan block, this contradicts the non-emptiness condition of Proposition \ref{prop: nonempty cells for hessenberg}, and so $\xi=0$.

As $v$ varies over $Z_{i+1}$, the coefficients of $\Psi$ vary regularly, giving a morphism of vector bundles
\[
\Psi\colon\mathcal O_{Z_{i+1}}\otimes_\C\cR\longrightarrow\mathcal O_{Z_{i+1}}\otimes_\C\cT.
\]
The preceding argument proves that this morphism is surjective on every fiber. Hence $K=\ker(\Psi)$ is a vector bundle of rank $\dim\cR-\dim\cT$. The right-hand side of \eqref{eqn: fiber of pi equations} is a regular section of the target bundle, and its inverse image is $Z_i$. Thus $\pi_i\colon Z_i\to Z_{i+1}$ is a torsor under $K$.

If $Z_{i+1}\simeq\Abb^r$ is an affine space, then $H^1(Z_{i+1},K)=0$, so this torsor is trivial. By the Quillen--Suslin theorem, $K$ is a trivial vector bundle on $\Abb^r$. Therefore,
\[
Z_i\simeq \Abb^{\,r+\dim\cR-\dim\cT}.
\]
By downwards induction on $i$, we deduce that $Z_1\simeq\cH_w(x,\fm)$ is an affine space, as claimed.
\end{proof}

The affine paving above will be used to give a combinatorial formula for the Poncar\'e polynomial and Frobenius characteristic of $\cH(x,\fm)$. For the parabolic case $\cH_P(x,\fm)$, we use the following.

\begin{lem}
    \label{lem: parabolic case of hessenberg is fibration}
    Suppose that $P = P_\alpha$ is a (standard) parabolic of $H$ with $\alpha = (\alpha^{(1)},\dots, \alpha^{(n)})$ for $\alpha^{(\ell)}$ a composition of $d$, and let $\fm\subset \fg_\xi$ be a Hessenberg space stable under the action of $P$. Let $x$ have all eigenvalues nonzero. The projection map
    \[
    \cH(x,\fm)\to \cH_P(x,\fm)
    \]
    is a trivial fibration with fibers
    \[
    \prod_{\ell = 1}^n\prod_{i}\Fl_{\alpha_i^{(\ell)}}
    \]
    where $\Fl_{\alpha_i^{(\ell)}} = \GL_{\alpha_i^{(\ell)}}/B_{\alpha_i^{(\ell)}}$ is the flag variety for $G = \GL_{\alpha_i^{(\ell)}}$.
\end{lem}
\begin{proof}
    The projection $G/B\to G/P$ is a Zariski locally trivial fibration with fibers as above, so it suffices to observe that, for any $g\cdot P\in G/P$, if $g\cdot P\in \cH_P(x,\fm)$, then $g'\cdot B\in \cH(x,\fm)$ for all $g'\in g\cdot P$.
\end{proof}

\section{Cyclic chromatic functions}
\label{sec: chromatic}

\begin{defn}
    For each $1\leq i\leq n$, let $V_i = \{\mathbf{1}_i,\dots, \mathbf{d}_i\}$ be a set of $d$ vertices. A cyclic $n$-partite digraph on $V = \bigsqcup_i V_i$ with step size $b$ is a directed graph $D$ where, if $v\to w$ is an arrow of $D$ and $v\in V_i$, then $w\in V_{i+b}$. Here, all indices on the sets $V_i$ are considered mod $n$.
\end{defn}

\begin{ex}
\label{ex: digraph from hessenberg}
    Suppose $m=an+b$ for $b$ coprime to $n$. Let $\mathbf{m}$ be a tuple of $n$ Hessenberg functions of size $d$ each, written as
\[\mathbf{m}=\left( (m_1^{(1)},\dots, m_d^{(1)}), \dots, (m_1^{(n)},\dots, m_d^{(n)})\right).\]
    Define a cyclic $n$-partite digraph $D$ with only arrows from $\mathbf{i}_\ell$ to $\mathbf{j}_{\ell+b}$ for $i> m^{(\ell)}_j$. 

    For example, when $(dm,dn)=(2,4)$ and the Hessenberg function tuple is $((1,2),(1,2))$, we have arrows $(\mathbf{2}_1\to \mathbf{1}_2),(\mathbf{2}_2\to \mathbf{1}_1)$. 
\end{ex}

The complement digraph $\ol{D}$ is the complement of $D$ in the sense above:  That is, $\ol D$ is the cyclic $n$ partite digraph which contains the edge $\mathbf{i}_\ell$ and $\mathbf{j}_{\ell+b}$ if and only if $D$ does not contain that edge. For example \ref{ex: digraph from hessenberg}, $\ol{D}$ simply has arrows between $\mathbf{i}_\ell$ and $\mathbf{j}_{\ell+b}$ when $i\leq  m^{(\ell)}_j$. For consistency with \cite{BrosnanChow}, we will work mainly with $\ol D$ rather than $D$ in our notation.

Take a composition $(\beta_1,\dots, \beta_r)=\beta \models d$, and put $s_j = \sum_{i\leq j}\beta_i$.

\begin{defn}
\label{def: suitable}
	A sequencing $\eta=(\eta_\ell(1),\dots, \eta_\ell(d))_{\ell=1}^n$ of the vertices is $(\ol{D},\beta)$-suitable if the following two conditions hold:
	\begin{enumerate}
		\item For every $1\leq t\leq d$, there is a loop $\eta_1(t)\to \eta_{1+b}(t)\to\dots \to \eta_{1-b}(t) \to \eta_{1}(t)$ in $\ol D$. Denote such a loop by $L_t$.  (This condition is independent of $\beta$.)
		\item For every $0\leq j\leq r-1$, there is at least one directed path connecting $L_{s_{j}+1}\to  \dots \to L_{s_{j+1}}$ in $\ol D$. (We write $s_0=\beta_0=0$.) 
	\end{enumerate}
\end{defn}
We can rephrase this condition as follows:  Condition (1) above states that $\eta$ determines disjoint loops $L = (L_1,\dots,L_d)$ in $\ol D$ (so that every loop $L_i$ is of length $n$ and every vertex of $\ol D$ appears in a unique loop $L_i$). We may form the digraph $\ol D/L$ on $d$ vertices by identifying the vertices in each loop $L$. We have an edge $L_i\to L_j$ in $\ol D/L$ whenever there is an edge from some vertex in $L_i$ to some vertex in $L_j$. By construction, $\eta$ descends to an ordering of the vertices of $\ol D/L$. Then, $\eta$ is $(\ol D,\beta)$ suitable if and only if $\eta$ is an ordered path cover for $\ol D/L$ in the sense of \cite[Definition 23]{BrosnanChow}.
\begin{defn}
    Define the \emph{descent}, $\des(\eta)$, of a $(\ol D,\beta)$-suitable sequencing $\eta$ to be the number of arrows $\eta_\ell(\mathbf{i})_\ell\to \eta_{\ell+b}(\mathbf{j})_{\ell+b}$ in $\ol D$ such that $\mathbf{i}<\mathbf{j}$ and $\eta_\ell(i)>\eta_\ell(j)$. 
\end{defn}
For cyclic $n$-partite graphs as in example \ref{ex: digraph from hessenberg}, we have
\begin{equation}
\label{eqn: descent for hessenberg case}
\des(\eta)=\#\left\{(\ell,i,j), 1\leq \ell\leq n, 1\leq i<j \leq d, m^{(\ell)}_{\eta_{\ell+b}(j)}\geq \eta_\ell (i)>\eta_\ell(j)\right\}.
\end{equation}
\begin{defn}
    The \emph{cyclic chromatic function} of a cyclic $n$-partite digraph $\ol D$ is the function
    \[
    \Xi_{\ol{D}}(Y,t):=\sum_{(\eta,\beta)} t^{\des(\eta)}M_\beta
    \]
    where the sum is over compositions $\beta \models d$ and $\eta$ which are $(\ol D,\beta)$-suitable. Here, $M_\beta$ is the monomial quasisymmetric function for the composition $\beta$ in the variables $Y=(y_1,y_2,\ldots)$.
\end{defn}

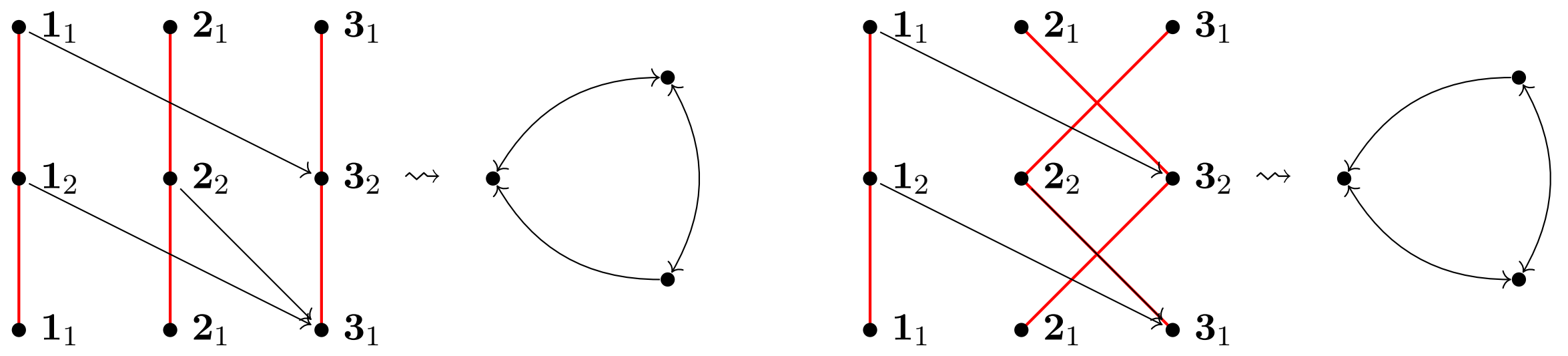
\begin{figure}
    \centering    
    \begin{tikzpicture}
        \draw[red, thick] (-1.5,1.5) -- (-1.5,-1.5);
        \draw[red,thick] (0,1.5) -- (0,-1.5);
        \draw[red,thick] (1.5,1.5) -- (1.5,-1.5);
        \node[circle,fill=black,inner sep=0pt,minimum size=4pt,label=right:{$\mathbf{1}_1$}] (a) at (-1.5,1.5) {}; 
        \node[circle,fill=black,inner sep=0pt,minimum size=4pt,label=right:{$\mathbf{2}_1$}] (a) at (0,1.5) {}; 
        \node[circle,fill=black,inner sep=0pt,minimum size=4pt,label=right:{$\mathbf{3}_1$}] (a) at (1.5,1.5) {}; 
        \node[circle,fill=black,inner sep=0pt,minimum size=4pt,label=right:{$\mathbf{1}_2$}] (a) at (-1.5,0) {}; 
        \node[circle,fill=black,inner sep=0pt,minimum size=4pt,label=right:{$\mathbf{2}_2$}] (a) at (0,0) {}; 
        \node[circle,fill=black,inner sep=0pt,minimum size=4pt,label=right:{$\mathbf{3}_2$}] (a) at (1.5,0) {}; 
        \node[circle,fill=black,inner sep=0pt,minimum size=4pt,label=right:{$\mathbf{1}_1$}] (a) at (-1.5,-1.5) {}; 
        \node[circle,fill=black,inner sep=0pt,minimum size=4pt,label=right:{$\mathbf{2}_1$}] (a) at (0,-1.5) {}; 
        \node[circle,fill=black,inner sep=0pt,minimum size=4pt,label=right:{$\mathbf{3}_1$}] (a) at (1.5,-1.5) {}; 
        \draw[->] (-1.4,1.45) -- (1.4,.05);
        \draw[->] (-1.4,-.05) -- (1.4,-1.45);
        \draw[->] (0.1,-.1) -- (1.4,-1.4);
        \draw (2.5,0) node {$\rightsquigarrow$};
        \node[circle, fill=black, inner sep=0pt, minimum size=4pt] (v1) at (3.2,0) {}; 
        \node[circle, fill=black, inner sep=0pt, minimum size=4pt] (v2) at (4.932,1) {}; 
        \node[circle, fill=black, inner sep=0pt, minimum size=4pt] (v3) at (4.932,-1) {}; 
        \draw[<->] (v1) to[bend left=30] (v2);
        \draw[<->] (v2) to[bend left=30] (v3);
        \draw[->] (v3) to[bend left=30] (v1);
    \end{tikzpicture}
    \hspace{.5in}
    \begin{tikzpicture}
        \draw[red, thick] (-1.5,1.5) -- (-1.5,0) -- (-1.5,-1.5);
        \draw[red, thick] (0,1.5) -- (1.5,0) -- (0,-1.5);
        \draw[red,thick] (1.5,1.5) -- (0,0) -- (1.5,-1.5);
        \node[circle,fill=black,inner sep=0pt,minimum size=4pt,label=right:{$\mathbf{1}_1$}] (a) at (-1.5,1.5) {}; 
        \node[circle,fill=black,inner sep=0pt,minimum size=4pt,label=right:{$\mathbf{2}_1$}] (a) at (0,1.5) {}; 
        \node[circle,fill=black,inner sep=0pt,minimum size=4pt,label=right:{$\mathbf{3}_1$}] (a) at (1.5,1.5) {}; 
        \node[circle,fill=black,inner sep=0pt,minimum size=4pt,label=right:{$\mathbf{1}_2$}] (a) at (-1.5,0) {}; 
        \node[circle,fill=black,inner sep=0pt,minimum size=4pt,label=right:{$\mathbf{2}_2$}] (a) at (0,0) {}; 
        \node[circle,fill=black,inner sep=0pt,minimum size=4pt,label=right:{$\mathbf{3}_2$}] (a) at (1.5,0) {}; 
        \node[circle,fill=black,inner sep=0pt,minimum size=4pt,label=right:{$\mathbf{1}_1$}] (a) at (-1.5,-1.5) {}; 
        \node[circle,fill=black,inner sep=0pt,minimum size=4pt,label=right:{$\mathbf{2}_1$}] (a) at (0,-1.5) {}; 
        \node[circle,fill=black,inner sep=0pt,minimum size=4pt,label=right:{$\mathbf{3}_1$}] (a) at (1.5,-1.5) {}; 
        \draw[->] (-1.4,1.45) -- (1.4,.05);
        \draw[->] (-1.4,-.05) -- (1.4,-1.45);
        \draw[->] (0.1,-.1) -- (1.4,-1.4);
        \draw (2.5,0) node {$\rightsquigarrow$};
        \node[circle, fill=black, inner sep=0pt, minimum size=4pt] (v1) at (3.2,0) {}; 
        \node[circle, fill=black, inner sep=0pt, minimum size=4pt] (v2) at (4.932,1) {}; 
        \node[circle, fill=black, inner sep=0pt, minimum size=4pt] (v3) at (4.932,-1) {}; 
        \draw[<-] (v1) to[bend left=30] (v2);
        \draw[<->] (v2) to[bend left=30] (v3);
        \draw[<->] (v3) to[bend left=30] (v1);
    \end{tikzpicture}
    \caption{The cyclic $2$-partite digraphs for the cases in example \ref{ex: chromatic}. The black arrows denote the arrows in $D$. The red lines indicate possible choices of disjoint cycles $L = (L_1,L_2,L_3)$ in $\ol D$. The graphs drawn are the quotient digraphs $\ol D/L$.}
    \label{fig: examples of D}
\end{figure}
\begin{defn}
    When $D$ is constructed from a tuple of Hessenberg functions $\fm$ as in example \ref{ex: digraph from hessenberg}, we denote $\Xi_\fm(Y,t):=\Xi_{\ol D}(Y,t)$.
\end{defn}
\begin{ex}
\label{ex: chromatic}
We compute $\Xi_{\fm}(Y,t)$ for $n=2$, $m=1$ and $d=3$. Here, $m_\lambda$ denotes the monomial symmetric function for the partition $\lambda$. The cyclic $n$-partite digraph $D$ for the first example below is pictured in Figure \ref{fig: examples of D}, with its possible loop decompositions.

	For $\mathbf{m}=((2,3,3),(1,3,3))$, \[\Xi_{\fm}(Y,t)=(1+5t+5t^2+t^3)m_{(1^3)}+(1+4t+4t^2+t^3)m_{(2,1)}+(1+3t+3t^2+t^3)m_{(3)}.\]
	
        For $\mathbf{m}=((2,2,3),(1,3,3))$,
	\[\Xi_{\fm}(Y,t)=(1+4t+t^2)m_{(1^3)}+(1+3t+t^2)m_{(2,1)}+(1+2t+t^2)m_{(3)}.\]
	
	For $\mathbf{m}=((2,3,3),(3,3,3))$,
	\[\Xi_{\fm}(Y,t)=(1+4t+7t^2+7t^3+4t^4+t^5)(m_{(1^3)}+m_{(2,1)}+m_{(3)}).\]
\end{ex}

\begin{lem}
    \label{lem: hessenberg poincare polynomial}    
    Suppose that $x\in \fg_\xi^\reg$ is of type $\lambda$ expressed in the form of Lemma \ref{lem: normal form of x}. Let the matrix $A$ appearing in the normal form for $x$ have Jordan blocks appearing in an order specified by the composition $\beta$ (which is a rearrangement of the parts of $\lambda$). Let $\fm$ have associated digraph $\ol D$ as in Example \ref{ex: digraph from hessenberg}. Then, the cohomology of $\cH(x,\fm)$ is concentrated in even degrees, and the Poincare polynomial $P_\fm(t)= \sum_i \dim\left (H^{2i}\big(\cH(x,\fm)\big)\right)t^{i}$ of the Hessenberg variety $\cH(x,\fm)$ is
    \[
P_\fm(t) = \sum_{(\eta,\beta)} t^{\des(\eta)}
    \]
    where the sum is over all $\eta$ which are $(\ol D,\beta)$ suitable.
\end{lem}
\begin{proof}
    Since $\cH(x,\fm)$ is projective with an affine paving, it is well known that its cohomology is concentrated in even degrees and that $H^{2i}\big(\cH(x,\fm)\big)$ is of dimension equal to the number of $i$ dimensional affine cells. For $w\in S_d^n$, we view $\eta = w^{-1}$ as an $n$ tuple $\eta = (\eta^{(1)},\dots, \eta^{(n)})$. It is immediate by Proposition \ref{prop: combinatorial cell dim for hessenberg} and \eqref{eqn: descent for hessenberg case} that, if $\cH^w(x,\fm)$ is nonempty, it is affine of dimension $\des(\eta)$. Therefore, it suffices to check that the condition in Proposition \ref{prop: nonempty cells for hessenberg} is the same as the condition that $\eta$ be $(\ol D,\beta)$ suitable. In terms of the digraph $\ol D$, the condition that $w\cdot \fm$ contains all diagonal entries is equivalent to condition (1) in Definition \ref{def: suitable} while the condition that for all $\alpha_{i,i+1}\in \Sigma_x$, there exists some $j$ such that $w\cdot \fm$ contains $E_{i+(b+j)d,i+jd+1}$ is equivalent to condition (2) in Definition \ref{def: suitable}. Note that, for the latter set of equivalant conditions, the composition $\beta$ appearing is determined by the ordering of the Jordan factors as this determines which super-diagonal entries of the $d\times d$ blocks must be in $w\cdot \fm$.
\end{proof}

\begin{cor}
    \label{cor: combinatorial poincare duality}
     The quasisymmetric function $\Xi_{\fm}(Y,t)$ is a symmetric function in the variables $Y$.
\end{cor}
\begin{proof}
    The coefficient of the monomial quasisymmetric function $M_\beta$ in $\Xi_{\ol D}(Y,t)$ is
    \[
    \sum_{(\eta,\beta)} t^{\des(\eta)}
    \]
    where the sum is over all $\eta$ which are $(\ol D,\beta)$ suitable. By Lemma \ref{lem: hessenberg poincare polynomial}, this sum is the Poincar\'e polynomial of the regular Hessenberg variety $\cH(x,\fm)$, where $x$ is of type $\lambda$. In particular, since $\cH(x,\fm)$ is depnds only on the conjugacy class of $x$, and therefore does not depend on the ordering of the Jordan blocks in the normal form of $x$, we see that the coefficient of $M_\beta$ is equal to the coefficient of $M_{\beta'}$ for $\beta'$ any rearrangement of $\beta$.
\end{proof}

In fact, $\Xi_{\fm}(Y,t)$ is a palindromic polynomial in $t$, as we will now show.

\begin{prop}
    \label{prop: combinatorial poincare duality}
     The function $\Xi_{\fm}(Y,t)$ is palindromic in $t$, in the sense that it satisfies the identity
     \[
     \Xi_{\fm}(Y,t) =  t^{D_\fm}\Xi_{\fm}(Y,t^{-1})
     \]
     where $D_\fm = |\mathbf{m}|-\frac{N(d+1)}{2} = \sum_{\ell,i}m_i^{(\ell)}-\frac{nd(d+1)}{2}$.
\end{prop}
    \begin{proof}
    Our proof is a direct analogue to \cite[Proposition 2.6]{sw1}. We describe an involution on pairs $(\eta,\beta)$ such that $\eta$ is $(\beta,D)$-suitable, which will reverse the descent statistic. First, we define a straightening procedure for words in the loops determined by $\eta$.

    Let $\cL$ be the set of these loops, and write $v_\ell(L)$ for the index of the vertex of $L$ in the vertex set $V_\ell$. Thus, for every $1\leq \ell\leq n$ and $L\in \cL$, we have $v_\ell(L)\leq m_{v_{\ell+b}(L)}^{(\ell)}$. Define the order $A\succ B$ if the arrow $A\to B$ is absent from the graph $D/\mathcal L$, or equivalently if $v_\ell(A)>m_{v_{\ell+b}(B)}^{(\ell)}$ for all $\ell$. This is a strict partial order. Indeed, $A\succ B$ implies
    $v_\ell(A)>v_\ell(B)$ for every $\ell$, and if
    $A\succ B\succ C$, then
    \[
        v_\ell(A)>v_\ell(B)>
        m_{v_{\ell+b}(C)}^{(\ell)}
    \]
    for every $\ell$, so $A\succ C$.

    To straighten a word, repeatedly replace adjacent loops $AB$ by $BA$ whenever $A\succ B$. Each exchange decreases the inversion number of the word of $V_1$-indices, so this procedure terminates. Moreover, the resulting word is independent of the choices of exchanges. Exchanges in disjoint positions commute, while overlapping exchanges occur only in a subword $ABC$ with $A\succ B\succ C$. By transitivity, either initial exchange can be followed by further exchanges to obtain $CBA$. Since the procedure terminates, these observations imply uniqueness of the resulting word, which we denote by $\operatorname{st}(W)$. A word is straightened precisely when its consecutive loops form a directed path in $D/\mathcal L$.

    The straightening preserves the descent
    statistic. Namely, for distinct loops $A,B$, put
    \[
        c(A,B)=
        \#\left\{\ell:
        v_\ell(B)<v_\ell(A)
        \leq m_{v_{\ell+b}(B)}^{(\ell)}
        \right\}.
    \]
    If $\eta$ orders the loops as $L_1,\dots,L_d$, then
    \[
        \des(\eta)=\sum_{i<j}c(L_i,L_j).
    \]
    Whenever $A\succ B$, we have $c(A,B)=c(B,A)=0$.
    Exchanging adjacent $A,B$ changes no other relative
    orders, and hence preserves $\des$.

    Now write $\beta=(\beta_1,\dots,\beta_r)$, and recall
    $s_0=0$ and $s_j=\beta_1+\cdots+\beta_j$. Thus
    \begin{equation}
        \label{eqn: beta partition}
        \{1,\dots,d\}
        =\bigsqcup_{j=1}^r\{s_{j-1}+1,\dots,s_j\}.
    \end{equation}
    Let
    \[
        W_j=(L_{s_{j-1}+1},\dots,L_{s_j})
    \]
    be the corresponding words in the loops. Define
    $\iota(\eta,\beta)=(\wt\eta,\beta^{\mathrm{rev}})$ by
    ordering the loops of $\wt\eta$ as
    \[
        \operatorname{st}(W_r^{\mathrm{rev}})
        \operatorname{st}(W_{r-1}^{\mathrm{rev}})
        \cdots
        \operatorname{st}(W_1^{\mathrm{rev}}).
    \]
    Here, reversing or exchanging loops means performing
    the same rearrangement in every $\eta_\ell$.
    The loops themselves are unchanged, and each resulting
    word is a directed path. Thus $\wt\eta$ is
    $(\beta^{\mathrm{rev}},D)$-suitable.

    This construction is an involution. Indeed, consider words
    equivalent if they differ by adjacent exchanges of
    $\succ$-comparable loops. Each equivalence class has a
    unique straightened representative, and reversal
    preserves this equivalence relation. Since each original
    word $W_j$ is already straightened, reversing and
    straightening twice returns $W_j$. The order of the
    blocks is also reversed twice.

    Finally, let $\eta^{\mathrm{rev}}$ denote the sequencing
    obtained by reversing the entire order of the loops.
    Complete reversal exchanges the order of every pair, so
    \[
        \des(\eta)+\des(\eta^{\mathrm{rev}})
        =\sum_{A\neq B}c(A,B).
    \]
    For fixed $\ell$ and $B$, the indices $v_\ell(A)$ run
    through $\{1,\dots,d\}$. Consequently,
    \begin{align*}
        \sum_{A\neq B}c(A,B)
        &=
        \sum_{\ell,B}
        \left(m_{v_{\ell+b}(B)}^{(\ell)}-v_\ell(B)\right)\\
        &=
        \sum_{\ell,i}m_i^{(\ell)}
        -\frac{nd(d+1)}{2}
        =D_\fm.
    \end{align*}
    Since straightening preserves the descent statistic,
    we obtain
    \[
        \des(\wt\eta)=D_\fm-\des(\eta).
    \]

    Let $P_\beta(t)$ be the coefficient of $M_\beta$ in
    $\Xi_\fm(Y,t)$. The involution gives
    \[
        P_\beta(t)
        =t^{D_\fm}P_{\beta^{\mathrm{rev}}}(t^{-1}).
    \]
    Since $\Xi_\fm(Y,t)$ is symmetric in $Y$, we have
    $P_{\beta^{\mathrm{rev}}}(t)=P_\beta(t)$.
    The desired identity follows.
\end{proof}

The following is immediate from Proposition \ref{prop: combinatorial poincare duality} and Lemma \ref{lem: hessenberg poincare polynomial}.

\begin{cor}
    \label{cor: poincare duality for hessenberg}
    Suppose that $P_{\fm,x}(t)= \sum_i \dim\left (H^{2i}\big(\cH(x,\fm)\big)\right)t^{i}$ is the Poincar\'e polynomial of the Hessenberg variety $\cH(x,\fm)$, where $x\in \fg_\xi^\reg$ has all eigenvalues nonzero. Then, $P_{\fm,x}(t) = t^{D_\fm}P_{\fm,x}(t^{-1})$ is a palindromic polynomial in $t$.
\end{cor}

Finally, we also collect the Poincar\'e polynomial for the parabolic case.
\begin{cor}
\label{cor: poincare polynomial for parabolic hessenberg}
    If $P = P_\alpha$ is a parabolic of $H$ corresponding to an $n$ tuple of compositions of $\alpha^{(\ell)}\models d$, if $x\in \fg_\xi^\reg$ has all eigenvalues nonzero, and if $\fm$ is $P$ stable, then the odd cohomology of $\cH_P(x,\fm)$ vanishes and the Poincar\'e polynomial for the even cohomology of $\cH_P(x,\fm)$ is
    \[
    \frac{P_\fm(t)}{\Aut_t(\alpha)},
    \]
    where $P_\fm(t)$ is the Poincar\'e polynomial of $\cH(x,\fm)$ (see Lemma \ref{lem: hessenberg poincare polynomial}) and
    \[
    \Aut_t(\alpha) = \prod_{\ell=1}^n\prod_i [\alpha_i^{(\ell)}]_t!
    \]
    is a product of $t$-factorials $[r]_t! := \prod_{k=1}^r\frac{1-t^k}{1-t}$ depending on $\alpha$.
\end{cor}
\begin{proof}
    Follows from Lemma \ref{lem: parabolic case of hessenberg is fibration} since the Poincar\'e polynomial of $\Fl_{r}$ is $[r]_t!$.
\end{proof}

\subsection{Grothendieck-Springer type maps and centralizers of regular, semisimple elements}
\label{sec: GS and centralizers}

Let $\fa\subset \fg_\zeta$ denote the subset of matrices where, in the representation of \eqref{eqn: normal form for x}, we have $A = A_1=A_2=\cdots=A_n$ is a diagonal $d\times d$ matrix. Recall that $\mu_n$ denotes the group of $n$-th roots of unity. The subset $\fa$ carries an action of the group $W_\fa:=N_H(\fa)/C_H(\fa) \simeq \mu_n^d\rtimes S_d$, where we may view $W_\fa\subset W$ as a subgroup of the usual Weyl group (where the embedding is induced by the inclusion of $\fa\subset \ft = \fc_\fg(\fa)$, with $\ft$ being a Cartan of $\fg$). In particular, $S_d$ acts by permuting diagonal entries of $A$ and $\mu_n^d$ acts by multiplying entries by roots of unity.

\begin{lem}
\label{lem: chevalley}
    Consider the GIT quotient 
    \[
    \fc:=\fg_\zeta\git H := \mathrm{Spec} (\C[\fg_\zeta]^H).
    \]
    Then, the restriction map of functions to $\fa\subset \fg_\zeta$ induces an isomorphism $\fg_\zeta\git H\simeq \fa\git W_\fa$.
\end{lem}
\begin{proof}
    This can be found in \cite[Theorem 7]{vinberg}.
\end{proof}

\begin{defn}
    The Grothendieck-Springer resolution for the representation of $H$ on $\fg_\zeta^\reg$ is the closed subvariety $\wt\fg_\zeta^\reg\subset \fg_\zeta^\reg\times \Fl_H$ defined by
    \[
    \wt\fg_\zeta^\reg = \{(x,B)\in \fg_\zeta^\reg\times \Fl_H\mid C_H(x)\subset B\}.
    \]
\end{defn}

There are obvious projections
\[
\fg_\zeta^\reg\xleftarrow{\pi_\zeta}\wt\fg_\zeta^\reg\xrightarrow{\pi_{\Fl}}\Fl_H.
\]

\begin{defn}
    The universal centralizer $I\subset H\times \fg_\zeta^\reg$ is the group scheme over $\fg_\zeta^\reg$ whose fiber over $x\in \fg_\zeta^\reg$ is the centralizer
\[
I_x = \{h\in H\mid h\cdot x = x\}.
\]
with the action $\cdot$ by the restriction of the adjoint action on $\fg$. We will also denote by $I = I^{\mathrm{rs}}$ be the restriction of $I$ to $\fg_\zeta^{\mathrm{rs}}$, which is a family of tori in $H$, when there is no confusion. For any $x\in \fg_\zeta^{\mathrm{rs}}$, the torus $I_x$ acts on the Flag variety $\Fl_H$ through the conjugation action of $H$.
\end{defn}

\begin{prop}
    \label{prop: GS res over rss locus}
    There is a Cartesian diagram
    \[
    \xymatrix{
    \wt\fg_\zeta^\mathrm{rs}\ar[r]\ar[d] & \coprod_{W_H/S_d} \fa\git\mu_n^d \ar[d] \\
    \fg_\zeta^{\mathrm{rs}}\ar[r] & \fc
    }
    \]
    In particular, the Weyl group of $H$, $W_H = S_d^n$, acts on $\wt\fg_\zeta^\mathrm{rs}$ over $\fg_\zeta^\mathrm{rs}$, and this action is free and transitive on fibers.
\end{prop}
\begin{proof}
    Fix $B_H\subset H$ the standard upper triangular Borel of $H$. For $w\in W_H/S_d$, we consider the space $H\times^{B_H}B_H\cdot w\cdot \fa^\rs$, which is smooth and comes equipped with a morphism
    \begin{equation}
    \label{eqn: induced space component of GS}
H\times^{B_H}B_H\cdot w\cdot \fa^\rs\to \fg_\zeta^\rs, \quad (h,x)\mapsto h\cdot x,
    \end{equation}
    where $\cdot$ above refers to the adjoint action. The morphism \eqref{eqn: induced space component of GS} is finite, \'etale, with fiber over $y\in \fg_\zeta^\rs$ given by
    \begin{align*}
\{(h,b\cdot w\cdot x)\mid (hbw)\cdot x = y\}/B_H&\simeq \{(h,w\cdot x)\mid x\in \fa,h\cdot w\cdot x = y\}/(I_x\rtimes \mu_n^d)\\
&\simeq \{x\in \fa\mid \text{there exists $h$ such that }h\cdot x = y\}/\mu_n^d \\
&\simeq S_d
    \end{align*}
    Furthermore, there is a commutative diagram
    \[
    \xymatrix{
H\times^{B_H}B_H\cdot w\cdot \fa^\rs \ar[r]\ar[d] & \fa\git \mu_n^d\ar[d] \\
\fg_\zeta^\rs\ar[r] & \fc
    }
    \]
    with the top horizontal arrow given by $(h,b\cdot w\cdot x)\mapsto \big(b\cdot w\cdot x\mod B_H\big)\in B_H\cdot \fa/B_H \simeq \fa\git \mu_n^d$. Since both vertical arrows are finite, \'etale of the same degree, if follows that the diagram above is Cartesian.

    Take the map
    \begin{equation}
    \label{eqn: w component of GS resolution}
H\times^{B_H}B_H\cdot w\cdot \fa^\rs\to \wt\fg_\zeta^\rs,\quad (h,x)\mapsto (h\cdot x,h\cdot B_H).
    \end{equation}
    It is clear this map commutes with the projection maps to $\fg_\zeta^\rs$. Moreover, \eqref{eqn: w component of GS resolution} is injective:  For any $(y,B)\in \wt\fg_\zeta^\rs$, if $(y,B)$ is in the image of \eqref{eqn: w component of GS resolution}, then it must be the image of $(h,h^{-1}\cdot y)$ where $B = h\cdot B_H$. Hence, the image of \eqref{eqn: w component of GS resolution} defines a component of $\wt\fg_\zeta^\rs$ which we denote by $\wt\fg_{\zeta,w}^\rs$. It is clear that $\wt\fg_{\zeta,w}^\rs = \wt\fg_{\zeta,w'}^\rs$ whenever $w$ and $w'$ belong to the same coset in $W_H/S_d$ (with $S_d$ embedded diagonally), so we may further label these components by $[w]\in W_H/S_d$. Moreover, since $w\cdot \fa\cap \fa$ contains no regular elements whenever $w\in W_H\setminus S_d$, these components $\wt\fg_{\zeta,[w]}^\rs$ are distinct for distinct cosets of $W_H/S_d$. We obtain an inclusion 
    \begin{equation}
    \label{eqn: inclusion}
    \coprod_{[w]\in W_H/S_d}\wt\fg_{\zeta,[w]}^\rs\subset \wt\fg_{\zeta}^\rs
    \end{equation}
    defined over $\fg_\zeta^\rs$. Now, for any $x\in \fg_\zeta^\rs$, $x$ is $H$ conjugate to an element $x_0$ of the form \eqref{eqn: normal form for x} with the matrices $A$ diagonal, which evidently satisfies that $C_H(x_0)$ is a subtorus of the diagonal torus $H$ which contains regular semisimple elements in $\Lie(H)$. Hence, the fiber over any $x\in \fg_\zeta^\rs$ is isomorphic to the fiber over $x_0$, which consists of all Borels containing the diagonal maximal torus in $H$. Therefore, both sides of the inclusion \eqref{eqn: inclusion} above have $|W_H| = (d!)^n$ many points in each fiber over $\fg_\zeta^\rs$, we conclude that \eqref{eqn: inclusion} is an equality, and collecting the Cartesian diagrams above, we obtain
     \[
    \xymatrix{
    \wt\fg_\zeta^\mathrm{rs}\ar[r]\ar[d] & \coprod_{[w]\in W_H/S_d} \fa\git\mu_n^d \ar[d] \\
    \fg_\zeta^{\mathrm{rs}}\ar[r] & \fc
    }
    \]
    In particular, $W_H$ acts naturally on fibers of $\pi_\zeta$.
\end{proof}

\begin{lem}
    \label{lem: GS cover gives localization family}
    The Grothendieck-Springer space $\wt\fg_\zeta^{\mathrm{rs}}$ is isomorphic to the space 
    \[
    (\Fl_H\times \fg_\zeta^{\mathrm{rs}})^{I}:= \{(B,x)\mid B\in (\Fl_H)^{I_x}\}
    \]
\end{lem}
\begin{proof}
    It is immediate that $B\in (\Fl_H)^{I_x}$ if and only if $I_x\subset B$, which is the definition of $\wt\fg_\zeta^\reg$.
\end{proof}

Let $C = C_H(\fa)$ be the centralizer of $\fa$. More explicitly, $C$ is the rank $d$ torus given by the diagonal embedding of the set of $d\times d$ diagonal matrices in $\GL_d$ into $H = \GL_d^n$. There is a natural action of the diagonal $S_d\subset W_H$ on $C$.

We will describe $I^{\mathrm{rs}}$ using a distinguished component of the Grothendieck-Springer cover. To this end, consider the diagonal embedding $\Fl_{\GL_d}\to \Fl_H$ sending a Borel $B_d\subset \GL_d$ to $B_d\times \cdots \times B_d\subset H$. We call the image of this map $\Fl_\Delta\subset \Fl_H$. Let 
\[
(\wt\fg_\zeta^\rs)_\Delta = \wt\fg_\zeta^\rs\cap\pi_{\Fl}^{-1}(\Fl_\Delta).
\]
That is, $(\wt\fg_\zeta^\rs)_\Delta\subset \wt\fg_\zeta^\rs$ consists of pairs $(x,B)$ such that $B = B_d\times \cdots B_d\subset \GL_d^n = H$ is a diagonal Borel. Then, it is easy to check that $(\wt\fg_\zeta^\rs)_\Delta$ is a component of $\wt\fg_\zeta^\rs$ and that the map $\pi_\Delta\colon (\wt\fg_\zeta^\rs)_\Delta\to \fg_\zeta^\rs$ is an \'etale $S_d$ cover fitting into a Cartesian diagram
\[
\xymatrix{
(\wt\fg_\zeta^\rs)_\Delta\ar[r]\ar[d] & \fa\git \mu_n^d \ar[d] \\
\fg_\zeta^\rs\ar[r] & \fc
}
\]

\begin{prop}
    \label{prop: I is trivial over GS cover}
    We have an isomorphism $\pi_\Delta^*I^{\mathrm{rs}}\simeq C\times \wt\fg_\zeta^\mathrm{rs}$, where $S_d$ acts diagonally on the product. 
    
    In particular, $\pi_\zeta^*I^\rs$ is also trivial with fibers $C$.
\end{prop}
\begin{proof}
    Let $\bT$ denote the canonical torus of $H$. We claim first that $\bT$ contains a canonical subtorus, which we denote by $\bC\subset \bT$, such that $C$ is identified with $\bC$ under any identification $C\subset T\simeq \bT$. Indeed, let $T$ be the diagonal maximal torus, which clearly contains $C$, and let $B\supset T$ be any diagonal Borel that contains $T$. Then, $B$ determines an isomorphism $T\to \bT$, and the image of $C$ under this map is independent of choice of diagonal Borel since $C$ is stable under the $S_d$ action.

    Now, consider the map 
    \[
    \pi^*_\Delta I\to \bT\times (\wt\fg_\zeta^\reg)_\Delta
    \]
    sending $(x,B,h)\mapsto h\mod [B,B]\in \bT$. It is easy to check that this map factors through $\bC\times (\wt\fg_\zeta^\reg)_\Delta$ when $B$ is diagonal and that it is an isomorphism over $(\wt\fg_\zeta^\rs)_\Delta$.

    By considering the $W_H$ action from Proposition \ref{prop: GS res over rss locus}, 
    \[
    \wt\fg_\zeta^\rs = \bigcup_{w\in W_H/S_d}w\cdot (\wt\fg_\zeta^\rs)_\Delta
    \]
    and so the description of $\pi_\Delta^*I^\rs$ also carries over to a description on $\wt\fg_\zeta^\rs$.
\end{proof}

We end this section by giving a construction of a map playing the role of Kostant section. Fix a section $\epsilon_d\colon \fgl_d\git \GL_d\to \fgl_d^\reg$. For example, we may take $\epsilon_d$ to be the Kostant slice associated to an $\mathfrak{sl}_2$ triple in $\fgl_d$ or the companion section
\[
\epsilon_d(a_1,\dots, a_d) =  \begin{pmatrix}
    0 &  &  &\cdots & (-1)^{d+1}a_d \\
    1 & 0 &  & \cdots& (-1)^da_{d-1} \\
    & \ddots & \ddots & & \vdots \\
    & & 1 & 0 & -a_2 \\
     & & & 1 & a_1 \\
\end{pmatrix}
\]
Define $\epsilon\colon \fc\to \fg_\zeta^\reg$ by identifying $\fc\simeq \mathfrak{gl}_d\git \GL_d$ (taking the characteristic polynomial $P_a(t)$ to $Q_a(t)$ where $P_a(t) = Q_a(t^n)$) and using the formula
\[
\begin{tikzpicture}[scale=.8]
\node at (-5,3) {$\epsilon(a) = $};

\node at (2.4,3.8) {$\ddots$};
  
\draw[thick] 
  (-1.8,2.4) rectangle (-2.8,1.4);

\node at (-1.4,1.2) {$\ddots$};

\draw[thick] 
  (-1,0.8) rectangle (-0,-0.2);

\draw[line width=1.4pt]
  (-3.1,6.1)
  .. controls (-3.6,4.2) and (-3.6,1.4) ..
  (-3.1,-0.3);

\draw[line width=1.4pt]
  (4.1,6.1)
  .. controls (4.6,4.2) and (4.6,1.4) ..
  (4.1,-0.3);



\draw[thick] 
  (3.8,2.4) rectangle (2.8,3.4);
\node at (3.3,2.9) {$I_d$};

\draw[thick] 
  (2,4) rectangle (1,5);
\node at (1.5,4.5) {$I_d$};

\node at (2.4,3.8) {$\ddots$};
  
\draw[thick] 
  (1,5) rectangle (-.1,6);
\node at (0.45,5.5) {$\epsilon_d(a)$};

\node at (-0.5,0.3) {$I_d$};
\node at (-2.3,1.9) {$I_d$};

\draw[
  decorate,
  decoration={brace,amplitude=6pt},
  line width=1.2pt
]
  (4.8,6.0) -- (4.8,2.4);

\node at (6.5,4.2) {$n-b$ blocks};
\draw[
  decorate,
  decoration={brace,mirror,amplitude=6pt},
  line width=1.2pt
]
  (-3.0,-0.6) -- (0,-0.6);
  
\node at (-1.5,-1.1) {$b$ blocks};
\end{tikzpicture}
\]

\begin{prop}
    \label{prop: kostant section}
    The map $\epsilon$ gives a section of the map $\fg_\zeta^\reg\to \fc$.
\end{prop}
\begin{proof}
    Follows immediately as the characteristic polynomial of $\epsilon(a)$ is equal to the characteristic polynomial of $\epsilon_d(a)$ evaluated at $t^n$.
\end{proof}

\begin{rem}
    Note that, as opposed to the case of the usual Chevalley map $\fg^\reg\to \fg\git G$, the Kostant section above is \emph{not} unique up to $H$ conjugation.
\end{rem}

\subsection{The monodromy Action and Families of Hessenberg Varieties}\label{sec: families of hessenbergs}

\begin{defn}
    Let $\fq\subset \fg_\xi$ denote the subvariety of $\fg_\xi$ consisting of elements of nonzero eigenvalues. We denote by $\cH(\fm)\subset \fq\times \Fl_H$ the family consisting of
    \[
    \cH(\fm) = \{(x,g)\mid x\in \fq,\; g\in \cH(x,\fm)\}
    \]
    We denote the family morphism by
    \[
    \pi_\fm\colon \cH(\fm)\to \fq.
    \]
\end{defn}
We will also denote by $\wt\fq^\reg\subset \wt\fg_\zeta^\reg$ for the preimage of $\fq^\reg$ in $\wt\fg_\zeta^\reg$, and similarly for $\wt\fq^\rs$.

\begin{lem}
\label{lem: family is smooth and family map is flat}
    The family $\cH(\fm)$ is smooth of dimension
    \[
    \frac{nd(d-1)}{2}+|\mathbf{m}| = \frac{nd(d-1)}{2}+\sum_\ell\sum_i m_i^{(\ell)},
    \]
    and the family morphism $\pi_\fm$ is flat over $\fq^\reg$.
\end{lem}
\begin{proof}
    The projection $\cH(\fm)\to \Fl_H$ taking $(x,g)\mapsto g\in H/B$ is an open subset of a vector bundle whose rank is
    \[
    |\mathbf{m}| = \sum_\ell\sum_i m_i^{(\ell)}.
    \]
    Since $\Fl_H$ is smooth of dimension $\frac{nd(d-1)}{2}$, the first statement follows.

    To show that $\pi_\fm$ is flat, by miracle flatness it suffices to show that the fibers have the same dimension. We have $\pi_\fm^{-1}(x) = \cH(x,\fm)$, which has dimension $D_\fm$ for every $x\in \fq^\reg$ by Lemma \ref{lem: hessenberg poincare polynomial}.
\end{proof}

\begin{cor}
    \label{cor: smoothness of rss hessenberg}
    If $x\in \fq^{\mathrm{rs}}$ is regular semisimple, then $\cH(x,\fm)$ is smooth.
\end{cor}
\begin{proof}
    By Lemma \ref{lem: family is smooth and family map is flat}, the map $\pi_\fm$ is flat over $\fq^\reg$ and $\cH(\fm)$ and $\fq^\reg$ are smooth.

    Consider the morphism
    \[
    f\colon H\times_{B_H}\fm\to \fg_\xi,\quad \quad f(h,y) = h\cdot y
    \]
    with the dot denoting adjoint action. We claim this map is smooth $\fg_\xi^\rs$. To this end, for $y\in \fg_\xi^\rs$, we compute that the differential of $f$ at $(1,y)$ has image $[\fh,y]+\fm$, where $\fh = \Lie(H)$. Suppose that $z\in \fg_\xi^*\simeq \fg_{\xi^{-1}}$ annihilates $[\fh,y]+\fm$. Then, for all $a\in \fh$,
    \[
    0=\mathrm{tr}(z[a,y]) = \mathrm{tr}([y,z]a),
    \]
    so that by nondegeneracy of the pairing on $\fh$, $[y,z]=0$. Since $z$ commutes with $y\in \fg_\xi^\rs$ and $z\in \fm^\perp\subset \fg_{\xi^{-1}}$, we must have $z = 0$. Therefore, the differential of $f$ at $(1,y)$ is surjective for any $y\in \fg_\xi^\rs$, and by acting by $H$ on the left, we see that $f$ is smooth over $\fg_\xi^\rs$, of relative dimension
    \[
    |\fm|+\dim(H/B_H)-\dim\fg_\xi = |\fm|+\frac{nd(d-1)}{2}-d^2n = |\fm|-\frac{nd(d+1)}{2} = D_\fm.
    \]
    The corollary now follows since for any $y\in \fg_\xi^\rs$, the projection $\pi\colon H\times_{B_H}\fm\to \Fl_H$ to the first factor produces an isomorphism from the fiber $f^{-1}(y)$ to $\cH(x,\fm)$.    
\end{proof}

\begin{defn}
    If $X$ is a complex manifold and $Y$ is a closed, nowhere dense analytic subspace of $X$ which contains $x$, then the \emph{local fundamental group} of $X\setminus Y$ at $x$ is the isomorphism class of $\pi_1(U\setminus Y,p)$ for $U$ a good neighborhood of $x$ (in the sense of \cite[Definition 46]{BrosnanChow}) and $p\in U\setminus Y$.
\end{defn}

The local fundamental group is an important object as it controls the monodromy action on nearby fibers to $x$. When $X'\to X$ is a Galois map of smooth complex algebraic varieties with Galois group $\Gamma$, let $Y\subset X$ be the ramification locus in $X$. For $x\in X$, we select a good neighborhood $U$ of $x$ and $p\in U\setminus Y$, so that $\pi_1(U\setminus Y,p)$ computes the local fundamental group of $X\setminus Y$ at $p$. A choice of preimage $p'\in X'$ determines an action of $\pi_1(U\setminus Y,p)$ a monodromy homomorphism
\[
\pi_1(U\setminus Y,p)\to \Gamma.
\]
The \emph{local monodromy group} of $X'\to X$ at $x$ is the image of this monodromy map, which is a subgroup of $\Gamma$ well defined up to conjugacy in $\Gamma$.

\begin{defn}
    We say that a point $a\in \fc = \fg_\zeta\git H$ is of type $\lambda$, where $\lambda\vdash d$ is a partition of $d$, if the associated characteristic polynomial $P_a(x)$ is of the form
    \[
    P_a(x) = \prod_i \prod_{\mu\in \mu_n}(x-\mu x_i)^{\lambda_i},
    \]
    where the $x_i^n\neq x_j^n$ for all $i\neq j$ and the $x_i$ are nonzero.
\end{defn}

It is clear from the definitions that $x\in \fg_\zeta^\reg$ is of type $\lambda$ (see Lemma \ref{lem: normal form of x}) if and only if its image in $\fc$ is of type $\lambda$. 

\begin{lem}
\label{lem: local monodromy}
    The local monodromy group of the map $(\fa\git \mu_n^d)^\rs\to \fc$ at a point $a\in \fc$ of type $\lambda$ is the Young subgroup $S_\lambda = \prod_i S_{\lambda_i}\subset S_d$. Similarly, the local monodromy of the map $\wt\fq^\rs\to \fq$ at a point $x\in \fq^\reg$ of type $\lambda$ is conjugated to $S_\lambda$.
\end{lem}
\begin{proof}
    By \cite[Proposition 106]{BrosnanChow}, for the first claim, it is enough to compute the stabilizer in $S_d$ of a preimage of $a$ in $\fa\git \mu_n^d$. If $a$ is of type $\lambda$, then any preimage of $a$ in $\fa$ is of the form \eqref{eqn: normal form for x} with 
    \[
    A_1 = \cdots  = A_n = \mathrm{diag}(\underbrace{x_1,\dots, x_1}_{\lambda_1\text{ times}}, \cdots ,\underbrace{x_k,\dots, x_k}_{\lambda_{k}\text{ times}}).
    \]
    The group $S_d$ acts on this element by permuting the coordinates of the diagonal matrix, and so the lemma follows immediately.

    To transfer this computation of local monodromy to $x\in \fg_\zeta$, we make use of the Cartesian diagram of Proposition \ref{prop: GS res over rss locus}. Namely, we may use the Kostant section $\epsilon\colon \fc\to \fg_\zeta$ (see Proposition \ref{prop: kostant section}) and an identical argument to \cite[Proposition 113]{BrosnanChow} to deduce that the local monodromy group at $x$ is the same as the local monodromy group of its image in $\fc$ with respect to the map
    \[
    \coprod_{W_H/S_d} (\fa\git \mu_n^d)^\rs\to \fc.
    \]
    The latter is the same as the local monodromy of a single of these $S_d$ covers, which was computed in the first half of this lemma.
\end{proof}

Consider the centralizer group scheme as in \S\ref{sec: GS and centralizers} but restricted over $\fq$, which by abuse of notation we also denote by
\[
I = \{(h,x)\in H\times \fq\mid h^{-1}\cdot x = x \}.
\]
The group $I$ acts on $\cH(\fm)$ via left multiplication. The restriction of $I$ to $\fq^\mathrm{rs}$ forms a torus bundle. We may therefore form the local system $\scrL(\fm)$ on $\fq^\mathrm{rs}$, whose fiber over $x\in \fq^\mathrm{rs}$ is the $I_x$ equivariant cohomology
\[
\scrL(\fm)_x = H^*_{I_x}\big(\cH(x,\fm)\big).
\]
More explicitly, $\scrL(\fm)$ is the pushforward of the trivial local system $\ul{\C}$ on the stack $[\cH^\rs(\fm)/I]$ to $\fq^\rs$.

The local system $\scrL(\fm)$ is related to the ordinary cohomology as follows. Let $\alpha\colon \mathbf{B}(I^{\mathrm{rs}})\to \fq^{\mathrm{rs}}$ be the fiberwise classifying stack of the torus bundle $I^{\mathrm{rs}}$ over $\fq^{\mathrm{rs}}$. The cohomology $\scrA = \bR \alpha_*\ul{\C}$ is a sheaf of polynomial algebras over $\fq^{\mathrm{rs}}$. Let $\scrA_+\subset \scrA$ denote the positive degree part
\[
\scrA_+ = \bigoplus_{k>0}\bR^{2k}\alpha_*\ul{\C}.
\]
Then, since generalized Hessenberg varieties have an affine paving (Proposition \ref{prop: affine paving of Hessenberg variety}), they are equivariantly formal and hence
\begin{equation}
\label{eqn: GKM for scrL}
    \scrL(\fm)/\scrL(\fm)\scrA_+ \simeq \bR\pi_{\fm,*}\ul{\C}.
\end{equation}

Now, fix $x_0\in \fq^\mathrm{rs}$. We will study the monodromy action of $\pi_1(\fq^{\mathrm{rs}},x_0)$ on the equivariant cohomology group $\scrL(\fm)_{x_0}$. The following reduces the monodromy action to an $S_d$ action.

\begin{prop}
\label{prop: local monodromy is thru Sd}
    The fundamental group $\pi_1(\fq^\mathrm{rs},x_0)$ acts on $\scrL(\fm)_{x_0}$ through the canonical quotient map $\pi_1(\fq^\mathrm{rs},x_0)\to S_d$.
\end{prop}
\begin{proof}
    Let $I$ act on $\wt\fq^\rs$ trivially, and let 
    \[
    \pi_\fq\colon [\wt\fq^\rs/I]\simeq \wt\fq^\rs\times \mathbf{B}(I)\to \fq^\rs
    \]
    be the structure map. By localization and Lemma \ref{lem: GS cover gives localization family}, we have an embedding of local systems
    \[
\scrL(\fm)\hookrightarrow \pi_{\fq,*}\ul\C.
    \]
    We claim that the action of $\pi_1(\fq^\rs,x_0)$ on the latter is trivial on the image of $\pi_1(\wt\fq^\rs,\wt x_0)$, where $\wt x_0$ is any preimage of $x_0$. Indeed, by Proposition \ref{prop: I is trivial over GS cover} $\mathbf{B}(I)$ is trivialized after pulling back along $\wt\fq^\rs\to \fq^\rs$. Hence, the action of $\pi_1(\fq^\rs,x_0)$ factors through the quotient map
    \[
    \pi_1(\fq^\rs, x_0)\to \pi_1(\fq^\rs, x_0)/\pi_1(\wt\fq^\rs,\wt x_0)\simeq S_d.\qedhere
    \]
\end{proof}

By \eqref{eqn: GKM for scrL}, we realize an action of $S_d$ on the ordinary cohomology group $H^*(\cH(x_0,\fm))$ as a corollary.

\begin{rem}
    One can also define a dot action on $H^*(\cH(x_0,\fm))$ mimicking \cite{julianna} and show that it coincides with the monodromy action defined above. 
\end{rem}
Fix elements $x_0\in \fq^\rs$ and $x_\lambda\in \fq^\reg$ of type $\lambda$. Take $B(x_\lambda)$ to be a sufficiently small ball centered at $x_\lambda$. 
The local invariant cycle theorem of Beilinson-Bernstein-Deligne says that for any integer $i$ the natural map 
\begin{align}
    \label{eqn: local invariant cycle}
H^i(\cH(x_\lambda,\fm))\to H^0\left(B(x_\lambda)\cap \fq^\rs,\bR^i\pi_{\fm,*}\C\right).
\end{align}
is surjective, and $H^0\left(B(x_\lambda)\cap \fq^\rs,\bR^i\pi_{\fm,*}\C\right)\cong H^i(\cH(x_0,\fm))^{\pi_1(B(x_\lambda)\cap \fq^\rs)}$. Since the Betti numbers of $\cH(x_\lambda,\fm)$ are palindromic by Proposition \ref{prop: combinatorial poincare duality}, the map \eqref{eqn: local invariant cycle} is an isomorphism by \cite[Theorem 102]{BrosnanChow}.
\begin{thm}
\label{thm: frob character of hessenberg}
    For $x_0\in \fq^\rs$, under the monodromy action, 
    \begin{align}
        \label{eqn: frob} \sum_i t^i\Frob\left(H^{2i}(\cH(x_0,\fm))\right)=\Xi_{\fm}(Y,t).
    \end{align}
    where $\Frob$ denotes the Frobenius character with respect to the $S_d$ action. Moreover, if $P = P_\alpha$ is the parabolic of $H$ corresponding to a tuple of compositions of $d$, $\alpha = (\alpha^{(\ell)})$, and if $\fm$ is stable under the action of $P$, then
    \begin{align}
        \label{eqn: frob parabolic} \sum_i t^i\Frob\left(H^{2i}(\cH_P(x_0,\fm))\right)=\frac{\Xi_{\fm}(Y,t)}{\Aut_t(\alpha)}.
    \end{align}
    where $\Aut_t(\alpha) = \prod_\ell\prod_i [\alpha^{(\ell)}_i]_t!$. In both cases, all odd degree cohomology of $\cH(x_0,\fm)$ or $\cH_P(x_0,\fm)$ vanishes.
\end{thm}
\begin{proof}    
    Odd degree cohomology vanishes by the existence of affine pavings (Proposition \ref{prop: affine paving of Hessenberg variety}). First assume that $P = B$ is a Borel. We prove that both sides of \eqref{eqn: frob} have the same pairing with the homogeneous symmetric function $h_\lambda$ for any $\lambda \vdash d$. By Frobenius reciprocity, 
\[
\left(\Frob\!\left(H^{2i}(\cH(x_0,\fm))\right),h_\lambda\right)=\dim\left(H^{2i}(\cH(x_0,\fm))^{S_\lambda}\right),
\]
    which equals $\dim(
    H^{2i}(\cH(x_\lambda,\fm)))$ by the discussion above, as $\pi_1(B(x_\lambda)\cap \fq^\rs)$ acts through its image in $S_d$ by Proposition \ref{prop: local monodromy is thru Sd}. It remains to recall that $(m_\mu,h_\lambda)=\delta_{\lambda,\mu}$ and use Lemma \ref{lem: hessenberg poincare polynomial}. 

    For the parabolic case, we argue identically using Corollary \ref{cor: poincare polynomial for parabolic hessenberg}.
\end{proof}
Since $\cH(x,\fm)$ is projective with an affine paving, Theorem \ref{thm: frob character of hessenberg} remains true when one replaces $H^i$ with $H^{BM}_i$, which is the statement we will use later.  Moreover, since the odd cohomology of $\cH_P(x,\fm)$ vanishes, it is equivariantly formal. This gives the following
formulas for their equivariant Frobenius series.

\begin{cor}
\label{cor: equivariant frob character of hessenberg}
    Keep notations from Theorem \ref{thm: frob character of hessenberg}. With the action of the torus $T=C_H(x_0)$, the graded Frobenius characters of the $T$-equivariant cohomology are
    \begin{align*}
        \sum_i t^i\Frob\left(H_T^{2i}(\cH(x_0,\fm))\right)=\Xi_{\fm}\left[\frac{Y}{1-t},t\right],\qquad
        \sum_i t^i\Frob\left(H_T^{2i}(\cH_P(x_0,\fm))\right)=\frac{\Xi_{\fm}\left[\frac{Y}{1-t},t\right]}{\Aut_t(\alpha)},
    \end{align*}
    where the square brackets denote plethystic substitution.
\end{cor}
\begin{proof}
    By equivariant formality,
    \[
H_T^{*}(\cH_P(x_0,\fm))\simeq H^{*}(\cH_P(x_0,\fm)) \otimes H_T^{*}(\mathrm{pt}),
    \]
    as graded $S_d$-modules. The equivariant parameters raise the cohomological degree by 2, and so the formula follows.
\end{proof}

On Borel-Moore homology, the equivariant parameters lower the homological degree by 2, so we get the analogous formulas
\begin{align}
\label{eqn:equivariant BM of hessenberg}
        \sum_i t^i\Frob\left(H^{T,BM}_{2i}(\cH(x_0,\fm))\right)=\Xi_{\fm}\left[\frac{Y}{1-t^{-1}},t\right],\qquad
        \sum_i t^i\Frob\left(H_{2i}^{T,BM}(\cH_P(x_0,\fm))\right)=\frac{\Xi_{\fm}\left[\frac{Y}{1-t^{-1}},t\right]}{\Aut_t(\alpha)}.
\end{align}

\section{Invariant sets and affine permutations}
\subsection{Invariant sets}
\label{sec: combo}
\begin{defn}
		An $M,N$ invariant set is a nonempty subset $\Delta\subset \Z$ satisfying the following conditions
        \begin{itemize}
            \item If $r\in \Delta$ then $r+M,r+N\in \Delta$;
            \item $\Delta$ is bounded from below;
            \item Let $a$ be the minimal element in $\Delta$. Then $ \#(\Z_{\geq a}\setminus \Delta)<\infty$.
        \end{itemize}
	\end{defn}
    Note that when $M$ and $N$ are coprime, the third condition follows from the first two. We may often drop the $M,N$ when it is clear from context and just refer to $\Delta$ as an invariant set. 
    
    Let $I_{M,N}$ be the set of all $M,N$-invariant sets. Let $I_{M,N}^+$ be the set of all $M,N$-invariant sets contained in $\Z_{>0}$. We call invariant sets $\Delta\in I_{M,N}^+$ \emph{positive}. For every $\Delta\in I_{M,N}^+$, we associate to it statistics of area and dinv. We begin with the area of $\Delta$.
	
        \begin{defn}
		For every positive $\Delta\in I_{M,N}^+$, the \emph{area} of $\Delta$ is the cardinality of its complement in $\Z_{>0}$, denoted 
		\[
		\area(\Delta) = \#(\Z_{>0}\setminus \Delta).
		\]
	\end{defn}
	
	To define $\codinv(\Delta)$, we define first the set of $N,M$ generators of $\Delta$.	
	
	\begin{defn}
		The set $\Ngen(\Delta)$ is the set of $r\in \Delta$ such that $r-N\not\in \Delta$.
	\end{defn}

	We will often picture an invariant set as a tuple of paths as follows. We imagine $d = \gcd(M,N)$ sheets, labeled $0,\dots, d-1$. For each sheet, we draw an array extending infinitely in the horizontal directions and having $n = N/d$ rows. On sheet $i$, we label the box at the origin by $i$, and let rows increase to the left by $N$ and columns increase upwards by $M$. (See Figure \ref{fig: labeling}.) The set $\Delta$ can be described pictorially as the unique tuple of lattice paths, one in each sheet, which moves upwards or leftwards only, begins on each sheet with a vertical segment, and for which $\Ngen(\Delta)$ is the set of elements containing a vertical step on their rightmost boundary.

The following definition is given in \cite{gmv} and can be seen to agree with the difference $\pdinv - \mathrm{ldinv}$ for the canonical labeling of \cite{wilson}. 
	\begin{defn}
		For any $\Delta\in I_{M,N}$, $\codinv(\Delta)$ is the number
		\[
\codinv(\Delta) = \sum_{g\in \Ngen(\Delta)}\#[g,g+M-1]\cap \Delta^c.
		\]
	\end{defn}
    
	We can further define labelings on an invariant set $\Delta$. These functions will play an important role in the combinatorics of parahoric affine Springer fibers in \S\ref{sec: parahoric case}.
    
	\begin{defn}
    \label{def: labeled invariant set}
		A \emph{labeled invariant set} is the data of an invariant set $\Delta$ and labeling function
		\[
	f\colon \Ngen(\Delta)\to \Z_{> 0}
		\]
		such that if $g,g+M\in \Ngen(\Delta)$, then $f(g)\leq f(g+M)$.
	\end{defn}

    We may reinterpret the above closer to the spirit of \cite{hikita} as follows. To an $(M,N)$ invariant set $\Delta$, we have an associated tuple of ribbons given by remembering only the squares labeled by elements of $\Ngen(\Delta)$. (See Figure \ref{fig: labeling}.) A function $f$ on $\Ngen(\Delta)$ is a labeling of $\Delta$ if and only if labels are non-decreasing moving upwards whenever the boxes they label are directly above each other. Note that boxes shown on the top row of each sheet below still have a box on top of them; each sheet should be viewed as an infinite grid which repeats itself every vertical $n$ rows, with a horizontal shift by $m$ squares.

    \begin{figure}
        \centering
        \begin{tikzpicture}[scale=1]

\node at (-4,1.5) {Sheet 2};

\def\leftcols{-2}   
\def\rightcols{3}   
\def\rows{2}        
\def\M{4}           
\def\N{6}          

\foreach \i in {\leftcols,...,\rightcols} {
  \foreach \j in {0,...,\rows} {

    \draw (\i,\j) rectangle ++(1,1);

    \pgfmathtruncatemacro{\label}{(\rightcols-\i-2)*\N + \j*\M+2}

    \draw[very thick,red] (0,0) -- (0,2);
    \draw[very thick,red] (2,2) -- (0,2);
    \draw[very thick,red] (2,2) -- (2,3);

    \node at (\i+0.5,\j+0.5) {\label};
  }
}

\end{tikzpicture}\quad
\begin{tikzpicture}[scale=1]

\draw[very thick,->] (-3,1.5) -- (-2,1.5);
\node at (-2.5,1.8) {$f$};

\draw (-1,0) rectangle ++(1,1);
\draw (-1,1) rectangle ++(1,1);
\draw (1,2) rectangle ++(1,1);

    \draw[very thick,red] (0,0) -- (0,2);
    \draw[very thick,red] (2,2) -- (0,2);
    \draw[very thick,red] (2,2) -- (2,3);

\node at (-0.5,0.5) {1};
\node at (-0.5,1.5) {2};
\node at (1.5,2.5) {1};

\end{tikzpicture}

\vspace{.2cm}

\begin{tikzpicture}[scale=1]

\node at (-4,1.5) {Sheet 1};

\def\leftcols{-2}   
\def\rightcols{3}   
\def\rows{2}        
\def\M{4}           
\def\N{6}          

\foreach \i in {\leftcols,...,\rightcols} {
  \foreach \j in {0,...,\rows} {

    \draw (\i,\j) rectangle ++(1,1);

    \pgfmathtruncatemacro{\label}{(\rightcols-\i-2)*\N + \j*\M+1}

    \draw[very thick,red] (2,0) -- (2,3);
    
    \node at (\i+0.5,\j+0.5) {\label};
  }
}
\end{tikzpicture}\quad
\begin{tikzpicture}[scale=1]

\draw[very thick,->] (-3,1.5) -- (-2,1.5);
\node at (-2.5,1.8) {$f$};

\draw (-1,0) rectangle ++(1,1);
\draw (-1,1) rectangle ++(1,1);
\draw (-1,2) rectangle ++(1,1);

    \draw[very thick,red] (0,0) -- (0,3);

\node at (-0.5,0.5) {2};
\node at (-0.5,1.5) {3};
\node at (-0.5,2.5) {4};

\node at (1.8,1.5) {};
\end{tikzpicture}
        \label{fig: labeling}
    \end{figure}

	We can attach the statistic $\lcodinv$ to a labeled invariant set in an analogous way to $\codinv(\Delta)$.

\begin{defn}
    For any pair of $\Delta\in I_{M,N}$ and a labeling $f$, define its labeled codinv $\lcodinv(\Delta,f)$ by the number of pairs $a,b\in \Ngen(\Delta)$ such that $a<b<a+M$ and $f(a)>f(b)$.
\end{defn}

Note that in the above definition, the pair $a,b\in \Ngen(\Delta)$ may be taken to lie in different sheets.

\begin{ex}
    Two sheets are pictured for $N=6$ and $M=4$ in Figure \eqref{fig: labeling}. The red line indicates the invariant set $\Delta$ containing all numbers to the left of the paths in each sheet, and $\Ngen(\Delta) = \{1,5,9,10,14,18\}$. One can verify $\area(\Delta) = 6$ and $\codinv(\Delta) = 6$. A labeling $f$ of the invariant set $\Delta$ from Figure \ref{fig: labeling}. For this labeling, we have $\lcodinv(\Delta,f) = 1$ counting the pair of $N$-gens $9$ and $10$.
\end{ex}

For a labeled invariant set $(\Delta,f)$, we denote by
\[
\codinv(\Delta,f) := \codinv(\Delta)+\lcodinv(\Delta,f)
\]

    There are differing conventions for the labelings in Definition \ref{def: labeled invariant set} and the statistic $\lcodinv$. The definition given above is the same as that given by Hikita in \cite{hikita} for ``negative labels.'' There is an analogous definition given in \cite{wilson} which replaces $f(g)\leq f(g+M)$ with $f(g) < f(g+M)$. Let us call $(\Delta,f)$ satisfying this latter condition a \emph{$W$-labeled invariant set}. Similarly, we define Wilson's statistic $\lcodinv_W(\Delta,f)$ to count the number of pairs $a,b\in \Ngen(\Delta)$ such that $a<b<a+M$ and $f(a)<f(b)$, and $\codinv_W(\Delta,f) = \codinv(\Delta)+\lcodinv_W(\Delta,f)$. We consider the functions
    \[
D_\Delta(X,t) = \sum_{f\colon \text{labeling of $\Delta$}} t^{\lcodinv(\Delta,f)}x^f,\qquad D_\Delta^{\mathrm{Wil}}(X,t) = \sum_{f\colon \text{ $W$-labeling of $\Delta$}} t^{\lcodinv_W(\Delta,f)}x^f
    \]
    where $x^f = \prod_{i\in\Ngen(\Delta)} x_{f(i)}$ is the monomial corresponding to the labeling $f$. We also put
    \[
    e_\Delta:=\#\{(a,b)\in \Ngen(\Delta)\mid a<b<a+M\}.
    \]
    We define the functions
    \[
L_{M,N}(X,q,t) = \sum_{\Delta\in I_{M,N}^+} q^{\area(\Delta)}t^{\codinv(\Delta)}D_\Delta(X,t),
\]
\[
L_{M,N}^{\mathrm{Wil}}(X,q,t) = \sum_{\Delta\in I_{M,N}^+} q^{\area(\Delta)}t^{-\codinv(\Delta)-e_\Delta}D_\Delta^{\mathrm{Wil}}(X,t).
    \]
    It is an easy exercise to see that $L_{M,N}^{\mathrm{Wil}}(X,q,t)$ is equal to the function that Wilson calls $L_{M,N} = L(0^M,0^N)$. In particular, $L_{M,N}^{\mathrm{Wil}}(X,q,t)$ recovers the triply graded superpolynomial of the $(M,N)$ torus link via a Schröder inner product (see Corollary 1 of \cite{wilson}).

    For our geometric setting, $L_{M,N}(X,q,t)$ will be the functions which appear naturally. They are related to Wilson's functions by the following lemma.
    \begin{lem}
        \label{lem: hikita vs wilson labeling convention}
        Let $\Delta$ be an invariant set. Then,
        \[
        D_\Delta(X,t) = t^{e_\Delta}\cdot \omega_X D_\Delta^{\mathrm{Wil}}(X,t^{-1})
        \]
        where $\omega_X$ is the canonical involution on (quasi)symmetric functions in the variables $X$. Moreover,
        \[
        L_{M,N}(X,q,t) = \omega_X L_{M,N}^{\mathrm{Wil}}(X,q,t^{-1}).
        \]
    \end{lem}
    \begin{proof}
The identity relating $D_\Delta$ and $D_\Delta^{\mathrm{Wil}}$ follows from an identical argument to \cite[Lemma 3.11]{hikita}. After this, the reciprocity for $L_{M,N}$ and $L_{M,N}^{\mathrm{Wil}}$ is immediate.
    \end{proof}

\subsection{Associated elements of the affine Weyl groups}

The affine Weyl group for $\GL_N$ is the group $W^\aff$ of bijective maps $w\colon \Z\to \Z$ such that $w(i+N) = w(i)+N$ for all $i$. For any element $w\in W^\aff$, we have a window notation for $w$ given by listing the $N$ tuple $[w(1),w(2),\dots,w(N)]$.

We will now attach two elements of the affine Weyl group:  $w_{\cell}$ and $w_{\hess}$ to a given labeled invariant set $(\Delta,f)$. We picture an invariant set $\Delta$ as a path tuple (with a sequence of $d$ sheets, each containing a path). We have a sequence $\vla\in \Z^N$ so that $\vla_i$ is the number of boxes in between the left side of box $i$ and the North step in the same row. Equivalently, $\vla_i = \frac{x_i - i}{N}$ where $x_i\in \Ngen(\Delta)$ is the unique $N$ gen with residue $i$ mod $N$. For $\vla\in \mathbb{X}_*(T)$, we let $w_{\vla}\in W^\aff$ be the corresponding element of the affine Weyl group determined by $w_\vla(i) = i+N\vla_i$ for $1\leq i\leq N$.

We say that a labeled invariant set $(\Delta,f)$ has content $\mu = (\mu_1,\dots, \mu_l)\vdash N$ if $\mu_i = \#f^{-1}(i)$. The affine Weyl group $W^\aff$ contains the finite Weyl group $W = S_N$, and so for any composition $\mu\models N$, we have an embedding $S_\mu = \prod_{i}S_{\mu_i}\subset W\subset  W_{\aff}$. For a labeling $f\colon \Ngen(\Delta)\to \Z_{>0}$ with content $\mu$, we determine a coset of the finite Weyl group $w_\fin\in S_\mu\backslash S_N$. When $\mu = (1^N)$, we let $w_\fin(i) = f(x_i)$ where $x_i\in \Ngen(\Delta)$ is the unique element such that $x_i = i\mod N$. In general, we use the following.

\begin{lem}
\label{lem: labelings to finite weyl group}
    There is a natural bijection between cosets in $S_\mu\backslash S_N$ and functions $f\colon \{1,\dots, N\}\to \Z$ such that $\#f^{-1}(i) = \mu_i$ for all $i$.
\end{lem}
\begin{proof}
    For a coset $S_\mu\cdot w$, we consider the function $f$ which sends $w^{-1}(1),\dots, w^{-1}(\mu_1)$ to $1$; sends $w^{-1}(\mu_1+1),\dots, w^{-1}(\mu_1+\mu_{2})$ to 2; etc. In general, we put 
    \begin{equation}
    \label{eqn: w fin for content mu}
        f^{-1}(j) = \left\{w^{-1}\!\left(1+\sum_{i\leq j-1}\mu_i\right),\dots, w^{-1}\left(\sum_{i\leq j}\mu_i\right)\right\}
    \end{equation}
    In the other direction, given a function $f$, we consider any element $w\in S_N$ such that \eqref{eqn: w fin for content mu} holds. This choice of $w$ is well defined up to reordering of the indices $f^{-1}(j)$, which is to say that $w$ is well defined up to the left action of $S_\mu$. It is easy to see these two constructions are mutually inverse.
\end{proof}

For $(\Delta,f)$ a labeled invariant set of content $\mu$, we therefore have a coset $[w_\fin] = S_\mu\cdot w_\fin\in S_\mu\backslash S_N$.

\begin{defn}
    For a labeled invariant set $(\Delta,f)$ of content $\mu$, we put 
    \[
    w_{\cell} =w_\cell(\Delta,f)= w_{-\vla}\cdot w_{\fin}^{-1}\in W^\aff/S_\mu \quad\text{and} \quad w_{\hess} =w_\hess(\Delta,f)= w_\fin \cdot w_{-\vla} \in S_\mu\backslash W^\aff.
    \]
\end{defn}

In particular, when $\mu = (1^N)$, we have $w_\cell(\Delta,f),w_\hess(\Delta,f)\in W^\aff$.

We will later use these constructions in conjunction with choices of representatives of cosets in $S_d^n\backslash W^\aff/S_\mu$ which satisfy the following dominance property.

\begin{defn}
\label{def: dominance}
    A labeled invariant set $(\Delta,f)$ is \emph{dominant} if:
\begin{enumerate}
    \item $\check\lambda$ is non-increasing in blocks of size $d$. That is, we have 
    \[
    \check\lambda_{1+(\ell-1)d}\geq \check\lambda_{2+(\ell-1)d}\geq \cdots \geq \check\lambda_{d+(\ell-1)d}
    \]
    for $1\leq \ell\leq n$.
    \item Whenever $\check\lambda_i = \check\lambda_{i+1}$ for $(i,i+1)$ in the same $d$ block, we have $f(x_i) \leq f(x_{i+1})$, where $x_i,x_{i+1}$ are the $N$ gens of $\Delta$ whose residue mod $N$ is $i$ and $i+1$, respectively.
\end{enumerate}
    We call $w_\cell = w_\cell(\Delta,f)$ \emph{dominant} whenever $(\Delta,f)$ is dominant.
\end{defn}

\section{Homology of affine Springer Fibers}\label{sec: ASF}

Throughout this section, we work with $G=\GL_N$. Let $K=\C(\!(t)\!)$ and $\Ocal=\C[\![t]\!]$. Inside the loop group $G_K$, we have the arc group $G_\Ocal$ and the Iwahori subgroup $I\subset G_\cO$. Let $W$ denote the finite Weyl group for $G$, $\mathbb{X}_*(T)$ be the coweight lattice and $W^\aff=W\ltimes \mathbb{X}_*(T)$ be the extended affine Weyl group, which we view as the space of bijective maps $w\colon \Z\to \Z$ such that $w(i+N) = w(i)+N$.

Let $u\in \gfrak_\Ocal$ denote an elliptic homogeneous element of slope $1/N$. Without loss of generality, we may choose $u=\alpha_0+\alpha_1+\cdots +\alpha_N = \sum_{i=1}^{N-1}E_{i,i+1}+tE_{N,1}$ where $\alpha_i$'s are the usual  affine simple roots ($\alpha_0=tw_0\alpha_h$ for the highest root $\alpha_h$ and the longest element $w_0\in W$), and let $v = \mathrm{diag}(v',\dots,v')$ where $v'$ is a diagonal $d\times d$ matrix whose eigenvalues are nonzero and have pairwise distinct $n$-th powers. Take 
\begin{equation}
\label{eqn: gamma}
\begin{tikzpicture}[scale=.6]
\node at (-5.7,3) {$\gamma=vu^M = $};

\node at (-3.9,3) {$t^k$};

\draw[thick] 
  (3.8,2.4) rectangle (2.8,3.4);

\draw[thick] 
  (2,4) rectangle (1,5);

\node at (2.4,3.8) {$\ddots$};
  
\draw[thick] 
  (1,5) rectangle (0,6);

\draw[thick] 
  (-1.8,2.4) rectangle (-2.8,1.4);

\node at (-1.4,1.2) {$\ddots$};

\draw[thick] 
  (-1,0.8) rectangle (-0,-0.2);

\draw[line width=1.4pt]
  (-3.1,6.1)
  .. controls (-3.6,4.2) and (-3.6,1.4) ..
  (-3.1,-0.3);

\draw[line width=1.4pt]
  (4.1,6.1)
  .. controls (4.6,4.2) and (4.6,1.4) ..
  (4.1,-0.3);

\node at (-.3,5.5) {$d$};

\node at (.5,6.35) {$d$};

\draw[thick] 
  (3.8,2.4) rectangle (2.8,3.4);
\node at (3.3,2.9) {$v'$};

\draw[thick] 
  (2,4) rectangle (1,5);
\node at (1.5,4.5) {$v'$};

\node at (2.4,3.8) {$\ddots$};
  
\draw[thick] 
  (1,5) rectangle (0,6);
\node at (0.5,5.5) {$v'$};

\node at (-0.5,0.3) {$tv'$};
\node at (-2.3,1.9) {$tv'$};

\draw[
  decorate,
  decoration={brace,amplitude=6pt},
  line width=1.2pt
]
  (4.8,6.0) -- (4.8,2.4);

\node at (7,4.2) {$n-b$ blocks};
\draw[
  decorate,
  decoration={brace,mirror,amplitude=6pt},
  line width=1.2pt
]
  (-3.0,-0.6) -- (0,-0.6);
  
\node at (-1.5,-1.3) {$b$ blocks};
\end{tikzpicture}
\end{equation}

We write ``$\cdot$" to denote the adjoint action. Consider the map 
\[a:G_K\times^I \Lie(I)\to \gfrak_K, \quad (g,x)\mapsto g\cdot x.\]
Suppose $d=gcd(M,N)$, $m=M/d$, $n=N/d$. Following \cite{hikita} we are interested in studying the BM homology of the affine Springer fiber $\Fl_\gamma:=a^{-1}(\gamma)$ (endowed with the reduced scheme structure) inside the affine flag variety $\Fl:=G_K/I$ when $d>1$. 

For any composition $\alpha$ of $N$, we may consider the parahoric subgroup $\Pfrak_\alpha\subset G_\cO$ which consists of all elements $g\in G_\cO$ which reduce mod $t$ to the standard parabolic of type $\alpha$ in $G(\C)$. We will further consider the parahoric affine Springer fiber $\Fl_\gamma^\alpha = a_\alpha^{-1}(\gamma)$ for 
\[a_\alpha:G_K\times^{P_\alpha} \Lie(\Pfrak_\alpha)\to \gfrak_K, \quad (g,x)\mapsto g\cdot x.\]
In particular, when $\alpha = (N)$, we recover the affine Grassmannian Springer fiber, which denote by $\Sp_\gamma = \Fl_{\gamma}^{(N)}$.

A new feature of this affine Springer is that it carries a natural $S_d$ action independent of its usual Springer action. Our main result will compute the graded Frobenius character for the $S_d\times S_N$ action on the BM homology of $\Fl_\gamma$.

\subsection{GKM stratification of homogeneous affine Springer fibers}
For any $w\in W^\aff$, we may express $w = w_{\check\lambda}\cdot w_\fin$, where $w_{\check\lambda}\in \mathbb{X}_*(T)$ and $w_\fin\in W$, and we may embed $w\in G_K$ by $t^{\check\lambda}\cdot w_\fin$, where $t^{\check\lambda} = \check\lambda(t)$ is the image of the uniformizer under the cocharacter $\check\lambda$ and $w_\fin$ is viewed as the permutation matrix with nonzero entries at $(w_\fin(i),i)$ for $1\leq i\leq N$. For compositions $\nu,\mu$ of $N$, we have a description of the two sided cosets of the loop group of $G$ as follows
\[
G_K/\Pfrak_\mu = \bigsqcup_{w\in S_\nu \backslash W^\aff/S_\mu} \Pfrak_\nu w \Pfrak_\mu/\Pfrak_\mu,
\]
with closures given by the order on $W^\aff$. Moreover, when $\Pfrak_\nu=I$ is an Iwahori, $I w \Pfrak_\mu/\Pfrak_\mu$ is isomorphic to an affine space, which is of dimension $\ell(w)$ when $\Pfrak_\mu = G_\cO$. 

Recall that $\check\rho$ is the half sum of all positive roots, given explicitly by $\check\rho = \left(\frac{N-1}{2},\frac{N-3}{2},\dots, \frac{1-N}{2}\right)$. For a composition $\mu$, we also define $\check\rho_\mu$ to be the cocharacter such that $\langle\check\rho_\mu,\alpha_i\rangle$ is 1 if the simple root $\alpha_i$ not contained in the Levi factor for $\Pfrak_\mu$ (i.e. $i = \sum_{j\leq r}\mu_j$ for some $r$) and 0 otherwise. For example, $\check\rho_{(1^N)} = \check\rho$ and $\check\rho_{(N)} = 0$.

Under suitable restrictions on $\nu$, \cite{GKM} show that such a decomposition will determine cells in parahoric affine Springer fibers by intersecting $\Fl_{\gamma}^\mu$ with the two sided cosets $\Pfrak_\nu w \Pfrak_\mu/\Pfrak_\mu$. We summarize their results for $\gamma$ the homogeneous slope $\frac{M}{N}$ element defined in \eqref{eqn: gamma}.

\begin{thm}
\label{thm: GKM}
    (\cite[Thm 1.2 and \S4.8]{GKM}) For $\gamma$ as in \eqref{eqn: gamma}, the intersection
    \[
\Fl_\gamma^\mu\cap (\Pfrak_\nu w \Pfrak_\mu/\Pfrak_\mu)
    \]
    is an iterated affine space bundle over a generalized Hessenberg variety which is Zariski locally trivial with affine fibers whenever $\Pfrak_\nu\subset \Pfrak := \Pfrak_{(d^n)}$, i.e. whenever $\nu$ is a refinement of the partition $(d^n)$.

    When $\nu = (1^N)$ (so that $\Pfrak_\nu = I$ is an Iwahori), for any $w\in W^\aff/S_\mu$ the cell $\Fl_\gamma^\mu\cap (Iw\Pfrak_\mu/\Pfrak_\mu)$ is an affine spaces of dimension
    \[
\#\{\wt{\alpha}\text{ affine root }\mid 0\leq \wt{\alpha}(\check\rho/N)<M/N\text{ and } \wt{\alpha}(w\cdot (\varepsilon\check\rho_\mu))<0\}
    \]

    When $\nu = (d^n)$, the cell $\Fl_\gamma^\mu\cap (\Pfrak_{(d^n)}w\Pfrak_\mu/\Pfrak_\mu)$ is an affine space bundle over the Hessenberg variety $\cH_P\!\left(\ol{\gamma},(\fg_{\zeta^b})_{w\cdot (\varepsilon\check\rho_\mu) - \frac{1}{n}\check\rho_{(d^n)},-\frac{M}{N}}\right)$ where:
    \begin{itemize}
        \item $\ol\gamma = \mathrm{ev}_{t=1}(\gamma)\in \fg_{\zeta^b}^\reg$ is the evaluation of $\gamma$ at $t = 1$;
        \item $0<\varepsilon<\!\!<1$ is a small parameter;
        \item For $\varphi\in \mathbb{X}_*(T)$ viewed as a map $\varphi\colon \bG_m\to T$ and $a\in \Q$, the space $(\fg_{\zeta^b})_{\varphi,a}$ denotes
        \[
(\fg_{\zeta^b})_{\varphi,a} = \bigoplus_{i\geq a}\{ X\in\fg_{\zeta^b}\mid \varphi(z)\cdot X = z^iX \},
        \]
        that is, the graded part in $\fg_{\zeta^b}$, under the grading induced by $\varphi$, bounded below by $a$;
        \item The parabolic $P\subset H = \GL_d^n$ is given by the cocharacter $w\cdot (\varepsilon\check\rho_\mu)$ in the sense that
        \[
        \Lie(P) = \bigoplus_{i\geq 0}\{X\in \fgl_d^n\mid (w\cdot\check\varepsilon\rho_\mu)(z)\cdot X = z^iX\text{ for all }z\in \bG_m\}
        \]
        is the non-negatively graded part of $\Lie(H)$ under the grading induced by $w\cdot (\varepsilon\check\rho_\mu)$.
    \end{itemize}
\end{thm}

\subsection{Families and the $S_d$-action}

We construct the $S_d$-action by monodromy, following the family construction of \cite[\S5.3.2]{oy}. The compatibility with the Hessenberg varieties of Theorem~\ref{thm: GKM} is the relative version of \cite[\S5.4.1 and Theorem~5.4.2]{oy}.

Recall $M=aN+bd$, where $0\leq b<n$, and recall that $N=nd$ and $M=md$. For $x\in\fg_\xi$, write $A_p$ for its block in position $(p,p+b)$, with block indices read modulo $n$, as in \eqref{eqn: normal form for x}. Define $\gamma_x\in\fg_K$ by replacing $A_p$ with
\[
(\gamma_x)_{p,p+b}=t^{a+\delta_p}A_p, \qquad \delta_p=\begin{cases}
0,&p\leq n-b,\\
1,&p>n-b.
\end{cases}
\]
Thus $x\mapsto\gamma_x$ is $H$-equivariant, $\gamma_x|_{t=1}=x$, and $\gamma_{\bar\gamma}=\gamma$ for $\gamma$ of the form \eqref{eqn: gamma} and $\bar\gamma = \gamma|_{t=1}$.

For a composition $\mu$ of $N$, consider the family
\[
q_\mu\colon\Fl^\mu\longrightarrow\fq^\rs,\qquad\Fl^\mu=\left\{(x,g\Pfrak_\mu)\mid \Ad(g^{-1})\gamma_x\in\Lie(\Pfrak_\mu)\right\}.
\]
Its fiber over $x$ is $\Fl_{\gamma_x}^\mu$. This family of homogeneous affine Springer fibers is a special case of \cite[\S5.3.2]{oy}, expressed in our block coordinates.

Homogeneity of the elements $\gamma_x$ gives a fiberwise $\bG_m$-action. Explicitly, put
\[
h(z)=\mathrm{diag}(z^{-1}I_d,z^{-2}I_d,\ldots,z^{-n}I_d),
\]
and let
\[
z\cdot(g(t)\Pfrak_\mu)=h(z)g(z^nt)\Pfrak_\mu.
\]
Indeed,
\[
h(z)\gamma_x(z^nt)h(z)^{-1}=z^m\gamma_x(t),
\]
so this action preserves every fiber of $q_\mu$.

Put $\Pfrak=\Pfrak_{(d^n)}$, whose Levi factor is $H$; it comes equipped with a projection $pr_H\colon \Pfrak\to H$. For $w\in W_H\backslash W^\aff/S_\mu$, define
\[
P_w=pr_H(\Pfrak\cap  w\Pfrak_\mu w^{-1}), \qquad
\fm_w=\left\{x\in\fg_\xi\mid \gamma_x\in\Ad(w)\Lie(\Pfrak_\mu)\right\}.
\]
Then $P_w$ is a parabolic of $H$ and the space $\fm_w$ is a $P_w$-stable Hessenberg space. We therefore have the Hessenberg family
\[
\mathcal H_w^\mu=\left\{ (x,hP_w)\in\fq^\rs\times H/P_w\mid h^{-1}\cdot x\in\fm_w \right\}, \qquad \pi^\mu_w\colon \cH^\mu_w\to \fq^\rs,
\]
whose fiber over $x\in \fq^\rs$ is the Hessenberg space $\cH_{P_w}(x,\fm_w)$. The map $hP_w\mapsto hw\Pfrak_\mu$ identifies this family with the corresponding part of the fixed locus:
\[
(\Fl^\mu)^{\bG_m} = \bigsqcup_{w\in W_H\backslash W^\aff/S_\mu} \cH_w^\mu.
\]
This is the identification of \cite[Lemma~5.3.5 and Theorem~5.4.2(1)]{oy}.

Let
\[
\mathcal E_w^\mu = \mathcal F^\mu\cap\bigl(\fq^\rs\times\Pfrak w\Pfrak_\mu/\Pfrak_\mu\bigr).
\]
The construction of Theorem~\ref{thm: GKM} works relatively over $\fq^\rs$. Namely, it gives an iterated affine space bundle
\[
p_w\colon\mathcal E_w^\mu\longrightarrow\mathcal H_w^\mu,
\]
whose projection is the contraction under the above $\bG_m$-action. In particular, this identifies the Hessenberg families appearing here with those studied in \S\ref{sec: families of hessenbergs}, rather than merely identifying their individual fibers.

Fix $x_0=\bar\gamma\in\fq^\rs$. We use the quotient
\[
\theta\colon\pi_1(\fq^\rs,x_0)\longrightarrow S_d
\]
from Proposition~\ref{prop: local monodromy is thru Sd}. Equivalently, if $\det(T-x)=Q_x(T^n)$, this quotient records the permutation of the $d$ distinct roots of $Q_x$.

\begin{prop}
\label{prop: affine Springer monodromy}
The monodromy of the family $\Fl^\mu$ induces $S_d$-actions on the cohomology and Borel--Moore homology of $\Fl_\gamma^\mu$. These actions are compatible with the filtration by $\Pfrak$-orbits. If $r_w$ is the relative dimension of $p_w$, the natural identifications
\[
H^*(E_w^\mu)\simeq H^*(\cH_{P_w}(x_0,\fm_w)),\qquad 
H_k^{\mathrm{BM}}(E_w^\mu)\simeq H_{k-2r_w}^{\mathrm{BM}}(\cH_{P_w}(x_0,\fm_w))
\]
are $S_d$-equivariant, where the action on the Hessenberg variety is the monodromy action of \S\ref{sec: families of hessenbergs}.
\end{prop}

\begin{proof}
We adapt the construction of \cite[Lemma~5.4.7]{oy} to our nonelliptic setting by working on finite projective truncations. Factorization through $S_d$ will follow from Proposition~\ref{prop: local monodromy is thru Sd} and equivariant localization.

Let $Z\subset G_K/\Pfrak_\mu$ be a finite union of closures of $\Pfrak$-orbits. Such varieties are projective, preserved by the $\mathbb G_m$-action, and exhaust $G_K/\Pfrak_\mu$. Intersecting with $\Fl^\mu$ gives a proper family
\[
q_{\mu,Z}\colon\Fl_Z^\mu\longrightarrow\fq^\rs.
\]
Its relative Hessenberg paving has successive strata which are iterated affine space bundles over the smooth proper families $\mathcal H_w^\mu$. Since the Hessenberg varieties have no odd cohomology, induction using the open--closed exact sequences shows that $R^kq_{\mu,Z,*}\mathbb Q$ are local systems, concentrated in even degrees.

Put $X=\Fl_\gamma^\mu\cap Z$. Applying the localization argument of \cite[Corollaries~5.4.4--5.4.5]{oy} to this proper family gives a monodromy-equivariant injection
\[
H^*_{\mathbb G_m}(X)\hookrightarrow H^*(X^{\mathbb G_m})\otimes\mathbb Q[u,u^{-1}],\qquad \deg(u)=2.
\]
Here injectivity follows from equivariant formality, since $H^{\mathrm{odd}}(X)=0$. The fixed locus is a finite disjoint union of the Hessenberg varieties $\cH_{P_w}(x_0,\fm_w)$. By Proposition~\ref{prop: local monodromy is thru Sd}, $\ker(\theta)$ acts trivially on their cohomology; it also fixes $u$. Consequently, it acts trivially on $H^*_{\mathbb G_m}(X)$ and on
\[
H^*(X)\simeq H^*_{\mathbb G_m}(X)/(u).
\]
Thus monodromy factors through $S_d$. Since $X$ is projective, duality gives the same factorization on its Borel--Moore homology.

These actions are compatible with inclusions of truncations. Passing to inverse limits for cohomology and direct limits for Borel--Moore homology therefore gives the actions on $\Fl_\gamma^\mu$.

Finally, the maps $p_w\colon\mathcal E_w^\mu\to\mathcal H_w^\mu$ are defined over $\fq^\rs$, so their pullback and Thom isomorphisms respect monodromy. This identifies the resulting actions on the strata with the Hessenberg actions of \S\ref{sec: families of hessenbergs}. The open--closed exact sequences likewise respect monodromy, giving compatibility with the filtration by $\Pfrak$-orbits.
\end{proof}

The same construction applies to finite closed unions of the strata $E_w^\mu$, and in particular to the positive area truncations considered in \S\ref{sec: positivity in ASF}. The forgetful maps between the partial affine Springer fibers are defined over $\fq^\rs$, so their induced maps on cohomology and Borel--Moore homology also respect these actions.

\subsection{Iwahori Orbits and affine pavings}
When $\Pfrak_\nu = I$ is an Iwahori, we denote the cells considered in the previous section by
\[
C_w^\mu:= \Fl_\gamma^\mu\cap I w \Pfrak_\mu/\Pfrak_\mu
\]

In the case $\Pfrak_\mu = G_\cO$, we have $W^\aff/W \simeq \mathbb{X}_*(T)$, and for $\vla\in \mathbb{X}_*(T)$ we may define the cell
\[
C_{\vla}:= \Sp_\gamma\cap It^\vla G(\cO)/G(\cO)
\]

Recall that, for our specific choice of $\gamma$, Theorem \ref{thm: GKM} says that the cells $C_{\vla}$ and $C_w^\mu$ above are affine spaces when they are nonempty.

\begin{prop}
\label{prop: nonempty cells}
Fix a coweight $\vla$, and let $w_\vla$ be the affine permutation such that $w_\vla(i)=i+N\vla_i$. The following are equivalent.
\begin{enumerate}[(a)]
    \item The cell $C_{-\vla}$ is nonempty.
    \item For all $1\leq i\leq N$, $w_\vla(i)-w_\vla(\ol{i-M})\leq M$. Here $\ol{i-M}\in [1,N]$ is the residue modulo $N$.
    \item $w_\vla(\Z_+)$ is an $M,N$-invariant set.
\end{enumerate}
Furthermore, there is a bijection between nonempty cells and $M,N$-invariant sets. 
\end{prop}
\begin{rem}
    The second condition is known to define the Sommers' region.
\end{rem}
\begin{proof}
    Take $t^{-\vla}\in G_K$. The cell $C_{-\vla}$ is nonempty if and only if $t^{-\vla}\in Sp_\gamma$. Recall that $M=kN+db$ and $N=dn$. By the definition of $\gamma$, $t^{\vla}\cdot \gamma\in \gfrak_\Ocal$ is equivalent to $\vla_1-\vla_{db+1}\geq -k$, $\dots$, $\vla_{N-db}-\vla_N\geq -k$, $\vla_{1+d(n-b)}-\vla_1\geq -k-1$, $\dots$, $\vla_N-\vla_{N-d(n-b)}\geq -k-1$. This family of inequalities is exactly equivalent to (b). 

    We now prove $(b)\Leftrightarrow (c)$. Write $x_i = w_\vla(i)$. Then
    \[
    w_\vla(\Z_+) = \bigcup_{i=1}^N \{x_i+Nj\mid j\geq 0\}.
    \]
    This set is $M$ stable if and only if for all $i$, 
    \[
    x_i+M \in \{x_{\ol{i+M}}+Nj\mid j\geq 0\}
    \]
    which holds if and only if $x_{\ol{i+M}}\leq x_{i}+M$ for all $i$, which is exactly condition (b).

    Thus we have a map from nonempty cells to $I_{M,N}$. We build an inverse map: given an $M,N$-invariant set $\Delta$, let $\{x_i\}_{i=1}^N$ be the set of Ngens such that $x_i\equiv i$ mod $N$. Define $\vla(\Delta)$ by $\vla_i(\Delta):=(x_i-i)/N$. Then $C_{-\vla(\Delta)}$ is nonempty.
\end{proof}

We note that while the affine Springer fiber $\Sp_\gamma$ depends only on the conjugacy class of $\gamma$, the cells $C_{\vla}$ depend on the specific choice of $\gamma$, and Proposition \ref{prop: nonempty cells} in particular depends on the form of $\gamma$ we chose.

In the special case when $\Pfrak_\mu = G_\cO$, it is immediate that Theorem \ref{thm: GKM} simplifies to the following. 
\begin{cor}
\label{cor_GKM}
     For $\vla\in \mathbb{X}_*(T)$, when $C_{\vla}$ is nonempty, it is an affine space of dimension
    \[
    \dim(C_\vla) = \#\{ (\alpha,u)\in \Phi\times \Z\mid 0\leq \langle \check\rho/N,\alpha\rangle + u <M/N\text{ and }\langle -\vla,\alpha\rangle+u<0 \}
    \]
\end{cor}

Due to our correspondence negating $\vla$, we will occasionally write $C_\Delta=C_{-\vla}$ when $\Delta = \Delta_{\vla}$.

\begin{prop} 
    $\dim(C_{-\vla})=\codinv(\Delta_{\vla})$.
\end{prop}
\begin{proof}
    Below, for an integer $z$, we write $\ol{z}\in[1,N]$ for the residue modulo $N$. Then, for $x\in \Ngen(\Delta),$ $w(\ol{x})=x$.

    Recall that $\codinv(\Delta)=\sum |[x,x+M-1]\cap \Delta^c|$ where $x$ runs over all $N$-gens. For a fixed pair of $N$-gens $x,y$ with $x<y$, we count all elements in the form of $y-uN$ for $u>0$ that lie in $[x,x+M-1]\cap \Delta^c$. In other words, the pair $(x,y)$ contributes to $\codinv(\Delta)$ by $\codinv_{x,y}$ which counts integers $u$ satisfying 
    \begin{align}
        \label{ineq: codinv}
    \frac{y-x-M}{N}<u<\frac{y-x}{N},\quad u>0.
    \end{align}
    On the dimension side, let $\alpha(x,y)=\epsilon_x-\epsilon_y$. We have 
    \[\langle \check{\rho},\alpha(x,y)\rangle=\ol{y}-\ol{x} \]
    and \[\langle \vla,\alpha(x,y)\rangle =-\frac{y-\ol{y}-(x-\ol{x})}{N}.\]
    Therefore, the pair $(x,y)$ contributes to $\dim(C_{-\vla})$ by $\dim_{x,y}$ which counts integers $u$'s satisfying 
    \begin{align}
        \label{ineq: dim}
    -\frac{\ol{y}-\ol{x}}{N}<u<\frac{M-(\ol{y}-\ol{x})}{N},\quad u<\frac{y-\ol{y}-(x-\ol{x})}{N}.\end{align}
    Clearly, $\dim_{y,x}=0$ for $x<y$. 
    Thus proposition would follow from the claim that $\codinv_{x,y}=\dim_{x,y}$ for $x<y$.

    By doing integer shifts, we see that the number of $u$'s satisfying \eqref{ineq: dim} equals the number of $u$'s satisfying 
    \[-\frac{y-{x}}{N}<u<\frac{M-({y}-{x})}{N},\quad u<0,\]
    which is symmetric to \eqref{ineq: codinv} with respect to $\{u=0\}$. Therefore the claim follows. 
\end{proof}
\subsection{The Parahoric Case}
\label{sec: parahoric case}
	Let $\mu$ be a partition of $N$, and let $\Pfrak_\mu$ denote the corresponding standard parahoric subgroup in $G_\cO$. Namely, $\Pfrak_\mu$ consists of all matrices $g\in G_\cO$ whose reduction mod $t$ lies in the standard parabolic subgroup of type $\mu$ in $G_\cO$. Recall that we denote the parahoric affine Springer fiber in $G_K/\Pfrak_\mu$ by $\Fl^\mu_{\gamma}$. Let $\pi_\mu: \Fl^\mu_\gamma\to \Sp_\gamma$ be the quotient map. From Theorem \ref{thm: GKM}, we have an affine paving of $\Fl^\mu_{\gamma}$ with cells $\Fl^\mu_{\gamma}\cap (It^{\check{\lambda}}w\Pfrak_\mu/\Pfrak_\mu) = (It^{-\vla}w^{-1}\Pfrak_\mu/\Pfrak_\mu)\cap \pi_\mu^{-1}(C_{\Delta})$ for $\vla\in \mathbb{X}_*(T)$ corresponding to an invariant set $\Delta$ and $w\in W/W_\mu$.
	
	\begin{prop}\label{prop: parahoric cells}
		For each $M,N$ invariant set $\Delta$ and for $\vla\in \mathbb{X}_*(T)$ corresponding to $\Delta$ under the bijection of Proposition \ref{prop: nonempty cells}, there is a bijection 
		\[\{w\in W_\mu\backslash S_N \mid  (It^{-\check{\lambda}}w^{-1}\Pfrak_\mu/\Pfrak_\mu) \cap \pi_\mu^{-1}(C_{-\vla})\neq \emptyset \} \Leftrightarrow \mathrm{Labelings}(\Delta,\mu)\]
		where $\mathrm{Labelings}(\Delta,\mu)$ denotes the set of labelings $f$ of $\Delta$ such that $\#f^{-1}(i) = \mu_i$ (see \S\ref{sec: combo}).
	\end{prop}

\begin{proof}
    As before, for an integer $z$, we write $\ol{z}\in\{1,\dots,N\}$ for the residue modulo $N$.

    For any invariant set $\Delta$, we may identify $\Ngen(\Delta)$ with $\{1,\dots,N\}$ by matching $1\leq i\leq N$ with the unique $N$ gen $x_i$ such that $x_i=i\mod N$. Then, by Lemma \ref{lem: labelings to finite weyl group}, we have a bijection between functions $f\colon\Ngen(\Delta)\to\Z$ of content $\mu$ and cosets $w\in S_\mu\backslash S_N$ defined so that \eqref{eqn: w fin for content mu} holds. We claim that this bijection restricts to the desired bijection.

    The element $t^{\vla}\cdot\gamma$ is a sum of elements in the root spaces $e_{i,i+M}$ where the subscripts are taken to be the corresponding representatives $1,\dots,N$ mod $N$. For $\vla\in\mathbb{X}_*(T)$, we let $w_{\vla}$ be the affine Weyl group element satisfying $w_{\vla}(i)=i+N\vla_i$. By the proof of Proposition \ref{prop: nonempty cells}, the entry of $t^{\vla}\cdot\gamma$ in
    $e_{i,i+M}$ has $t$ valuation $0$ if and only if $w_{\vla}(\ol{i+M})-w_{\vla}(i)=M,$ which is true if and only if the $N$ gens $x_i=w_{\vla}(i)$ and
    $x_{\ol{i+M}}=w_{\vla}(\ol{i+M})$ are directly above each other in the path notation for $\Delta$ (i.e. $x_{\ol{i+M}}=x_i+M$). In this case, the condition that $w\cdot e_{i,i+M}\in\mathfrak b$ be upper triangular is equivalent to the condition $w(i)\leq w(\ol{i+M})$. If $\mu=(1^N)$, then the labeling $f$ is given by $x_i\mapsto w(i)$, and the condition that $w(i)\leq w(\ol{i+M})$ whenever $x_{\ol{i+M}}=x_i+M$ is the condition that $f$ is a labeling of $\Delta$.

    For general $\mu$, refine $f$ to a function $f'$ of content $(1^N)$ by ordering the generators first by their labels under $f$, and then, among those with equal labels, by their integer values. The permutation $w'(i)=f'(x_i)$ lies in $S_\mu\cdot w$. For $x,x+M\in\Ngen(\Delta)$,
\[
f(x)\leq f(x+M)\quad\Longleftrightarrow\quad f'(x)<f'(x+M).
\]
    Thus every labeling $f$ admits a labeling refinement of content $(1^N)$; conversely, merging the labels of any such refinement gives a labeling of content $\mu$.

    Under the natural surjection $\Fl_\gamma\to\Fl_\gamma^\mu$, the inverse image of the intersection in the statement is
\[
\bigsqcup_{v\in S_\mu\cdot w}(It^{-\vla}v^{-1}I/I)\cap\pi^{-1}(C_{-\vla}).
\]
    By the case $\mu=(1^N)$, this is nonempty if and only if $f$ admits a refinement which is a labeling of content $(1^N)$. By the preceding paragraph, this is equivalent to $f$ being a labeling of content $\mu$.
\end{proof}    

	\begin{prop}\label{prop: parahoric dimension}
	When $C_{w_\cell(\Delta,f_w)}^\mu$ is nonempty, its dimension is
    \[
   \dim(C_{w_\cell(\Delta,f_w)}^\mu) = \dim((It^{-\check{\lambda}}w^{-1}\Pfrak_\mu/\Pfrak_\mu) \cap \pi_\mu^{-1}(C_{-\vla}))=\codinv(\Delta_\vla,f_w).
    \]
	\end{prop}
    \begin{proof} 
    Take $\ybf_\mu =c\check\rho_\mu$ for $0<c<\!\!<1$.  We have that $(t^{-\check{\lambda}}w^{-1})\cdot \ybf_\mu=w^{-1}\ybf_\mu+\vla$.

    Take two $N$-gens $x<y$.
    First, suppose $\langle {w}^{-1}\ybf_\mu ,\alpha(x,y)\rangle \geq 0$, or equivalently $\langle \ybf_\mu ,w\alpha(x,y)\rangle \geq 0$. By the proof of Proposition \ref{prop: parahoric cells}, the labeling $f_w$ satisfies $f_w(x)\leq f_w(y)$. Therefore, the contribution of $(x,y)$ to $\lcodinv(\Delta_\vla,f_w)$ is $0$. Using Theorem \ref{thm: GKM} in this case, we also have that the pair $(x,y)$ contributes to $\dim(C_{w_\cell(\Delta,f_w)}^\mu)$ by the number of $u$ satisfying the same inequalities in \eqref{ineq: dim}. 
    
    Now suppose $\langle {w}^{-1}\ybf_\mu ,\alpha(x,y)\rangle < 0$, or equivalently the labeling $f_w$ satisfies $f_w(x)> f_w(y)$. The second inequality in \eqref{ineq: dim} becomes
    \[u\leq \frac{y-\ol{y}-(x-\ol{x})}{N}.\]
    Therefore the pair $(x,y)$ has $1$ more contribution to $\dim(C_{w_\cell(\Delta,f_w)}^\mu)$ on top of $\codinv_{x,y}(\Delta_\vla)$ if and only if 
    \[\frac{y-\ol{y}-(x-\ol{x})}{N}<\frac{M-(\ol{y}-\ol{x})}{N}\]
    which is equivalent to $y<x+M$ or that $(x,y)$ contributes $1$ to $\lcodinv(\Delta_\vla,f_w)$.    
\end{proof}
\subsection{Positivity and the Hikita filtration}
\label{sec: positivity in ASF}
By Proposition \ref{prop: parahoric cells}, we have
\[
H^{BM}_*(\Fl_\gamma^\mu) = \mathrm{Span}_\C\!\left\{ C_{w_\cell(\Delta,f)}^\mu\mid (\Delta,f)\text{ is a labeled invariant set of content $\mu$} \right\}
\]
\begin{defn}
    The positive (flag) affine Springer fiber $\Fl_\gamma^{\mu,+}$ is the union
    \[
    \Fl_\gamma^{\mu,+} = \bigcup_{\substack{
        (\Delta,f)\colon \text{ content }\mu, \\
        \Delta\in I_{M,N}^+}} C_{w_\cell(\Delta,f)}^\mu
    \]
    where we recall $I_{M,N}^+$ was defined in \S\ref{sec: combo}. We endow the Borel-Moore homology $H^{BM}_*(\Fl_\gamma^{\mu,+})$ with an increasing filtration $F_\bullet$ given by
    \[
    F_aH^{BM}_*(\Fl_\gamma^{\mu,+}) = \mathrm{Span}_\C\!\left\{ C_{w_\cell(\Delta,f)}^\mu\;\;\left|\;\; \Delta\in I_{M,N}^+,\; \text{content}(\Delta,f) = \mu,\; \area(\Delta)\leq a \right. \right\}
    \]
    We call the filtration $F_\bullet$ the \emph{Hikita filtration}.
\end{defn}

This filtration is induced by taking BM homology of the subvariety of $\Fl_\gamma^{\mu,+}$ defined by
\[
\Fl_{\gamma,\leq a}^{\mu,+} = \bigcup_{\substack{
        (\Delta,f)\colon \text{ content }\mu, \\
        \Delta\in I_{M,N}^+,\; \area(\Delta)\leq a}} C_{w_\cell(\Delta,f)}^\mu.
\]
We prove that $\Fl_{\gamma,\leq a}^{\mu,+}$ is closed below, so that the area statistic induces a genuine stratification on $\Fl_{\gamma}^{\mu,+}$.

\begin{lem}
\label{lem: Bruhat order respects area}
    Suppose that $w = w_{-\vla}w_\fin^{-1}$ and $w' = w_{-\vla'}w_\fin'^{-1}$ are two elements of the affine Weyl group, with $\vla,\vla'\in \mathbb{X}_*(T)$ and $w_\fin,w_\fin'\in S_N$. Furthermore, let $\vla^{\mathrm{sort}}$ denote 
    \[
    \vla^{\mathrm{sort}} = \sigma(\vla) = (\vla_{\sigma(1)},\dots,\vla_{\sigma(N)})
    \]
    for $\sigma\in S_N$ so that $\vla_{\sigma(1)}\geq \vla_{\sigma(2)}\geq \cdots\geq \vla_{\sigma(N)}$. If $w\leq w'$ in the Bruhat order, then $\sum_{i\leq j}\vla_{i}^{\mathrm{sort}}\leq \sum_{i\leq j}\vla'^{\mathrm{sort}}_{i}$ for all $1\leq j\leq N$.
\end{lem}
\begin{proof}
Inversion preserves Bruhat order, so $w^{-1}\leq(w')^{-1}$. By \cite[Propositions 1.1(a) and 1.8(b)]{stembridge}, projection to the two-sided quotient by the finite Weyl group is order-preserving. Since $w^{-1}=w_\fin w_{\vla}$ and $(w')^{-1}=w_\fin' w_{\vla'}$, \cite[Theorem 4.10]{stembridge} then gives $\vla^{\mathrm{sort}}\preceq(\vla')^{\mathrm{sort}}$ in dominance order, which is precisely the asserted inequalities.
\end{proof}

\begin{cor}
    The subvariety $\Fl_{\gamma,\leq a}^{\mu,+}$ is closed.
\end{cor}
\begin{proof}
    Note that if $\Delta = \Delta_\vla$, then $\Delta$ is positive if and only if $\vla_i\geq0$ for every $i$. When $w=w_\cell(\Delta,f)=w_{-\vla}w_\fin^{-1}$ for positive $\Delta$, we have
    \[
    \area(\Delta)=\sum_{i=1}^N\vla_i.
    \]
    Suppose $w'=w_\cell(\Delta',f')$ indexes a nonempty cell with $\Delta'$ positive, and suppose further that $w\leq w'$. Since $w$ and $w'$ are comparable, we have $\sum_i\vla_i=\sum_i\vla_i'$. Hence, applying Lemma \ref{lem: Bruhat order respects area} for $j=N-1$ gives
    \[
    \min_i\vla_i\geq\min_i\vla_i'\geq0,
    \]
    In particular, the invariant set $\Delta$ corresponding to $w$
    is positive and $\area(\Delta)=\area(\Delta')$.

    Thus the cells defining $\Fl_{\gamma,\leq a}^{\mu,+}$ are stable under the Bruhat closure relations. There are only finitely many such cells, since $\vla_i\geq0$ and $\sum_i\vla_i\leq a$. Their union is therefore closed.
\end{proof}

In our main results, when we compute graded characters, we will always use the variable $q$ to denote the grading coming from $F_\bullet$. That is, we associate the polynomial
\[
\mathrm{gr}_{F}H^{BM}_*(\Fl_\gamma^{\mu,+}) = \sum_a q^a\dim\mathrm{gr}_{F}^aH^{BM}_*(\Fl_\gamma^{\mu,+})
\]

\subsection{Parahoric orbits and Hessenberg pavings} 
Let $\Pfrak:=\Pfrak_{(d^n)}$ be the standard parahoric of type $(d^n) = (d,d,\dots,d)$. Take 
\[
\ol{\gamma}:=\mathrm{ev}_{t=1}(\gamma) = \sum_{\alpha\in \Delta_{km}} u_\alpha e_\alpha + \sum_{\alpha\in \Delta_{N-km}}u_\alpha e_\alpha.
\]
For $S_d^n\cdot w \in S_d^n\setminus W^\aff$ and $S_d^n\cdot w\cdot S_\mu\in S_d^n\setminus W^\aff/S_\mu$, let 
\[
E_w:= \Fl_\gamma\cap \Pfrak w I/I,\quad E_w^\mu = \Fl_\gamma^\mu\cap \Pfrak w \Pfrak_\mu /\Pfrak_\mu
\]
\begin{prop}
\label{prop: parahoric cells nonempty}
    The cell $E_w^\mu$ is nonempty if and only if there exists an element $w' \in S_d^nwS_\mu/S_\mu$ such that $w' = w_{\cell}(\Delta,f)$ for a labeled invariant set $(\Delta,f)$ of content $\mu$.
\end{prop}
\begin{proof}
    The cell $E_w^\mu$ is a disjoint union of Iwahori orbits $C_{w'}^\mu$ where $w'$ varies among all elements in the left coset $S_d^n\cdot w\cdot S_\mu/S_\mu$, so this follows from Propositions \ref{prop: nonempty cells} and \ref{prop: parahoric cells}.
\end{proof}

In fact, the proof of Proposition \ref{prop: parahoric cells nonempty} gives an affine paving of $E_w^\mu$. Moreover, the following lemma puts nonempty cells $E_w^\mu$ in bijection with dominant labeled invariant sets $(\Delta,f)$.

\begin{lem}
\label{lem: unique dominant Delta,f}
    Let $(\Delta,f)$ be a labeled invariant set of content $\mu$. There is a unique dominant $(\Delta_+,f_+)$ such that $w_+ = w_\cell(\Delta_+,f_+)\in S_d^n\cdot w_\cell(\Delta,f)\cdot S_\mu/S_\mu$.
\end{lem}
\begin{proof}
    Suppose that $(\Delta,f)$ is a labeled invariant set of content $\mu$. It is clear that there is a unique coset $w_+S_\mu\in S_d^nw_\cell(\Delta,f)S_\mu/S_\mu$ with a representative $w_+ = w_{-\vla}w_\fin^{-1}$ satisfying
    \[
    \vla_{1+(\ell-1)d}\geq \vla_{2+(\ell-1)d}\geq\cdots \geq \vla_{d+(\ell-1)d}, \text{ for }1\leq \ell\leq n
    \]
    \[
    \text{if }\vla_{i} = \vla_{i+1}\text{ for }i\not\equiv 0\mod d,\quad \text{then }w_\fin(i)\leq w_{\fin}(i+1).
    \]
    We claim that this $S_\mu$ coset is of the form $w_+ = w_\cell(\Delta_+,f_+)$ for a labeled invariant set $(\Delta_+,f_+)$ of content $\mu$. For this, we note that $\vla_+ = s^{-1}(\vla)$ for $s\in S_d^n$ and $\vla$ satisfies the condition of Proposition \ref{prop: nonempty cells}. Therefore, as $\vla_+$ is sorted in blocks of size $d$, we have
    \[
    \max_{1\leq i\leq d}\left((\vla_+)_{i+(\ell-1)d} - (\vla_+)_{\ol{i+(\ell-1)d-M}} \right) \leq \max_{1\leq i\leq d}\left(\vla_{i+(\ell-1)d} - \vla_{\ol{i+(\ell-1)d-M}} \right)
    \]
    so that by Proposition \ref{prop: nonempty cells}, $C_{-\vla}\neq\emptyset$ and $\Delta_+ = w_{\vla_+}(\Z_{>0})$ is an invariant set. If $\vla_+= s^{-1}(\vla)$, then the function $f_+\colon \Ngen(\Delta_+)\to \Z$ is determined by $w_\fin s^{-1}$ via Lemma \ref{lem: labelings to finite weyl group}. The labeling condition for $f_+$ is then an easy exercise. We conclude that $w_+ = w_\cell(\Delta_+,f_+)$, as claimed.
\end{proof}

The dominant element $w_+ = w_\cell(\Delta_+,f_+)$ also minimizes $\dim(C_w^\mu)$, as we prove below.

\begin{lem}
\label{lem: dominant has minimal codinv}
    Let $(\Delta,f)$ be a labeled invariant set of content $\mu$, and let $(\Delta_+,f_+)$ be the unique dominant labeled invariant set of content $\mu$ such that $w_+ = w_\cell(\Delta_+,f_+)\in S_d^n\cdot w_\cell(\Delta,f)\cdot S_\mu/S_\mu$. (See Lemma \ref{lem: unique dominant Delta,f}.) Then,
    \[
    \codinv(\Delta_+,f_+) = \min\left\{\codinv(\Delta',f')\mid w_\cell(\Delta',f')\in S_d^n\cdot w_\cell(\Delta,f)\cdot S_\mu\right\}.
    \]
\end{lem}
\begin{proof}
    For an invariant set $(\Delta,f)$ and for $z\in \Z$, we assign the auxiliary ``state'' value
\[
\st(z) = \begin{cases}
    0 & \text{if }z\not\in\Delta \\
    f(z) & \text{if }z\in \Ngen(\Delta) \\
    \infty & \text{else}
\end{cases}
\]
    It is immediate from the definition of a labeling that, for any labeled invariant set, $\st(z)\leq \st(z+M)$. Moreover, if $(\Delta,f) = (\Delta_+,f_+)$ is dominant, then it is easy to check that
    \begin{equation}
    \label{eqn: states for dominant}
    \st(1+(\ell-1)d)\leq \st(2+(\ell-1)d)\leq \cdots\leq \st(\ell d)
    \end{equation}
    for every $\ell\in\Z$.. An inspection of the definition of $\codinv(\Delta,f)=\codinv(\Delta)+\lcodinv(\Delta,f)$ shows that this quantity is the sum of contributions
    \[
    C(x) = \#\{ z\in [x,x+M-1]\mid \st(z)<f(x) \}
    \]
    where $x\in \Ngen(\Delta)$. For $x = i+(\ell-1)d$ with $1\leq i\leq d$ and $\ell\in\Z$, we can rewrite this contribution as
    \begin{align*}
        C(x) = &\#\left\{ j\mid i< j\leq d,\; \st\big(j+(\ell-1)d\big)<f(x) \right\}\\ &+\sum_{r = 1}^{m-1}\#\left\{ j\mid 1\leq j\leq d,\; \st\big(j+(\ell+r-1)d\big)<f(x) \right\} \\
        &+ \#\left\{ j\mid 1\leq j< i,\; \st\big(j+(\ell+m-1)d\big)<f(x) \right\}
    \end{align*}
    After summing over all generators $x$, the central terms depend only on the multiset of generator labels in each starting block and the multiset of states in each block being counted. These multisets are preserved by the left $S_d^n$ action, so the total contribution of the central terms is independent of ordering. Let us focus on the outer two terms:
    \begin{align*}
    C'(x) ={}& \#\left\{ j\mid i< j\leq d,\; \st\big(j+(\ell-1)d\big)<f(x) \right\}\\ &+\#\left\{ j\mid 1\leq j< i,\; \st\big(j+(\ell+m-1)d\big)<f(x) \right\}. 
    \end{align*}
    Note that if $j>i$ and $\st\big(j+(\ell+m-1)d\big)<f(x)$, then we have
    \[
\st\big(j+(\ell-1)d\big)\leq\st\big(j+(\ell+m-1)d\big)<f(x).
    \]
    Moreover, the term $j=i$ in the second count contributes
    zero since $\st(x+M)\geq f(x)$. Hence, we may absorb these terms into the second set and rewrite
    \begin{align*}
C'(x) = \#\big\{ j\mid & i< j\leq d,\; \st\big(j+(\ell-1)d\big)< f(x)\leq \st\big(j+(\ell+m-1)d\big) \big\} \\
&+\#\left\{ j\mid 1\leq j\leq d,\; \st\big(j+(\ell+m-1)d\big)<f(x) \right\}. 
    \end{align*}
    The second summand above is again independent of ordering after summing over all $N$ gens and hence is constant on the left $S_d^n$ coset.

    Finally, we now specialize to the case of dominant $(\Delta_+,f_+)$. In this case, \eqref{eqn: states for dominant} implies that for all $i<j\leq d$ we have
    \[
f(x) = \st\big(i+(\ell-1)d\big)\leq \st\big(j+(\ell-1)d\big),
    \]
    so that
    \[
\#\big\{ j\mid i< j\leq d,\; \st\big(j+(\ell-1)d\big)< f(x)\leq \st\big(j+(\ell+m-1)d\big) \big\} = 0
    \]
    for all $x = i+(\ell-1)d\in \Ngen(\Delta)$. Since all other terms have constant total contribution on the $S_d^n$ coset, and the remaining terms are nonnegative, we conclude that $\codinv(\Delta_+,f_+)$ is minimal in $S_d^n\cdot w_\cell(\Delta,f)\cdot S_\mu/S_\mu$.
\end{proof}

\begin{defn}
    Suppose that $(\Delta,f)$ is a labeled invariant set with content $\mu = (1^N)$. 
    We associate a parabolic $P = P_{(\Delta,f)}$ in $H = \GL_d^n$ by the condition that the root space corresponding to the entry $(i+(\ell-1)d,j+(\ell-1)d)$ for $1\leq i,j\leq d$ and $1\leq \ell\leq n$ is contained in $\mathrm{Lie}(P)$ if and only if
\[
w_\hess(i+(\ell-1)d) - w_\hess(j+(\ell-1)d)<0.
\]
    For $1\leq i,j\leq d$ and $1\leq \ell\leq n$, let $1\leq i_\ell,j_\ell\leq N$ be the residues mod $N$ of
    \[
    i_\ell = (\ell-1)d+i\mod N,\quad \quad j_\ell = (b+\ell-1)d+j\mod N
    \]
    We associate a Hessenberg space $\fm = \fm_{(\Delta,f)}\subset \fg_\xi$, stable under $P_{(\Delta,f)}$, to $(\Delta,f)$ by the condition that the root space corresponding to the entry $(i_\ell,j_\ell)$ is contained in $\fm$ if and only if
    \begin{equation}
    \label{eqn: hessenberg definition}
w_\hess(i_\ell) - w_\hess(j_\ell)< \begin{cases}
    aN  & \text{for }\ell\leq n-b \\
    (a+1)N & \text{for }\ell> n-b
\end{cases}
    \end{equation}
\end{defn}

\begin{lem}
\label{lem: dominant implies standard hessenberg}
    If $(\Delta,f)$ is a \emph{dominant} labeled invariant set with content $\mu = (1^N)$ (see Definition \ref{def: dominance}), then $\fm_{(\Delta,f)}$ is standard in the sense that:
    \begin{itemize}
        \item $\fm_{(\Delta,f)}$ contains the root space for entry $(i_\ell,j_\ell)$ whenever $i\leq j$. (That is, each $d\times d$ block of $\fm_{(\Delta,f)}$ contains a standard Borel.)
        \item $P_{(\Delta,f)}$ is the standard, upper triangular Borel of $H =\GL_d^n$.
    \end{itemize}
    In particular, for such $(\Delta, f)$, $\fm_{(\Delta,f)}$ is equivalent to the data of an $n$ tuple of Hessenberg functions $\mathbf{m} = (\mathbf{m}^{(1)},\dots, \mathbf{m}^{(n)})$.
\end{lem}
\begin{proof}
    Since $\vla$ is non-increasing in blocks and $w_\fin$ is increasing whenever consecutive entries of $\vla$ are equal, we have
\[
w_\hess(i)-w_\hess(i+1)
=
w_\fin(i)-w_\fin(i+1)
+N(\vla_{i+1}-\vla_i)<0
\]
whenever $i$ and $i+1$ lie in the same block of size $d$. Indeed, if $\vla_i>\vla_{i+1}$, the last summand is at most $-N$, whereas $w_\fin(i)-w_\fin(i+1)\leq N-1$. If $\vla_i=\vla_{i+1}$, the inequality follows from dominance and the distinctness of the labels.

Thus $w_\hess$ is strictly increasing within each block. Consequently, the defining condition for $P$ selects precisely the upper triangular root spaces, so $P$ is the standard upper triangular Borel.

To prove the assertion about $\fm_{(\Delta,f)}$, it suffices to verify \eqref{eqn: hessenberg definition} for the diagonal positions within each $d\times d$ block.

Fix $1\leq i\leq d$ and $1\leq\ell\leq n$, and put
\[
r=i+(\ell-1)d,\qquad
\delta_\ell=
\begin{cases}
0,&\ell\leq n-b,\\
1,&\ell>n-b,
\end{cases}
\qquad
s=r+bd-\delta_\ell N.
\]
Thus $r,s\in\{1,\ldots,N\}$ index a diagonal entry of the $\ell$-th $d\times d$ block. Write $x_r=r+N\vla_r$ and $x_s=s+N\vla_s$ for the corresponding $N$-generators of $\Delta$.

Since $\Delta$ is $M$-invariant, $x_r+M\in\Delta$. This element has residue $s$ modulo $N$, so $x_s\leq x_r+M$. Using $M=aN+bd$, we obtain
\[
\vla_s-\vla_r\leq a+\delta_\ell.
\]
Moreover,
\[
w_{\mathrm{Hess}}(r)-w_{\mathrm{Hess}}(s)
=
w_\fin(r)-w_\fin(s)+N(\vla_s-\vla_r).
\]
If $\vla_s-\vla_r<a+\delta_\ell$, then integrality and $w_\fin(r)-w_\fin(s)\leq N-1$ imply
\[
w_{\mathrm{Hess}}(r)-w_{\mathrm{Hess}}(s)
<(a+\delta_\ell)N.
\]
If equality holds, then $x_s=x_r+M$. The labeling condition gives $f(x_r)\leq f(x_s)$, and this inequality is strict because $f$ has content $(1^N)$. Hence $w_\fin(r)=f(x_r)<f(x_s)=w_\fin(s)$, giving the same strict inequality.
\end{proof}

 We have the following combinatorial variant of Theorem \ref{thm: GKM} for $\Fl_\gamma$.
\begin{thm}
\label{thm: main for Fl gamma}
    Suppose that $(\Delta,f)$ is a labeled invariant set of content $\mu = (1^N)$. Recall that $\ol{\gamma} = \mathrm{ev}_{t=1}(\gamma)$ is the evaluation of $\gamma$ at $t = 1$. Then, $E_{w_\cell(\Delta,f)}$ is an affine bundle over the generalized Hessenberg variety $\cH(\ol{\gamma},\mfrak_{(\Delta,f)})$. Moreover, the rank of this affine bundle is exactly $\codinv(\Delta_+,f_+)$ where $w_+=w_\cell(\Delta_+,f_+)$ is the unique dominant representative of the coset $S_d^nw_\cell(\Delta,f)S_\mu/S_\mu$. (See Lemma \ref{lem: unique dominant Delta,f}.)
\end{thm}
\begin{proof}
    We must line up the spaces above with Theorem \ref{thm: GKM}. In this case, the parabolic of $H$ is determined by $w\cdot (\varepsilon\check\rho_{(1^N)}) = w\cdot (\varepsilon\check\rho)$, and since $\varepsilon$ is sufficiently small, this element in $\mathbb{X}_*(T)$ is regular (i.e. $\langle  w\cdot (\varepsilon\check\rho),\alpha_i \rangle \neq 0$ for all simple roots $\alpha_i$). Therefore, the parabolic $P$ of $H$ is Borel.

    We must now line up the Hessenberg space $\mfrak_{(\Delta,f)}$ with the one induced by the grading given by $w\cdot (\varepsilon\check\rho) - \frac{1}{n}\check\rho_{(d^n)}$. For $1\leq i,j\leq d$, we denote $(i_\ell,j_\ell)$ for the $(i,j)$ entry in the $\ell$-th block; that is,
 \[
 i_\ell = i+(\ell-1)d\mod N,\quad j_\ell = j+(b+\ell-1)d\mod N.
 \]   
    Write $w = w_{-\vla}\cdot w_{\fin}^{-1}$; then we have
    \[
    w\cdot (\varepsilon\check\rho) = \vla+\varepsilon( w_\fin^{-1}\cdot \check\rho).
    \]
    Therefore, an entry in position $(i_\ell,j_\ell)$ has grading under $w\cdot (\varepsilon\check\rho) - \frac{1}{n}\check\rho_{(d^n)}$ given by 
    \[
a+\delta_{\{\ell\leq n-b\}}-\frac{M}{N} + (\vla_{i_\ell}-\vla_{j_\ell}) + \varepsilon(\check\rho_{w_\fin(i_\ell)} - \check\rho_{w_\fin(j_\ell)})
    \]
    where $\delta_{\{\ell\leq n-b\}}$ is the indicator function for $\ell\leq n-b$. In particular, this is greater than or equal to $-\frac{M}{N}$ if and only if
    \[
\vla_{i_\ell}-\vla_{j_\ell}\geq\begin{cases}
     -a & \text{if }\ell\leq n-b \\
     -a-1 & \text{if }\ell>n-b
\end{cases}\quad \text{and \textit{if equality holds} in the previous, }w_\fin(i_\ell) < w_\fin(j_\ell),
    \]
    where we recall that $a = \lfloor\frac{M}{N}\rfloor$. It is not hard to see the above is equivalent to the following single inequality:
    \[
-N(\vla_{i_\ell}-\vla_{j_\ell}) + w_\fin(i_\ell) - w_\fin(j_\ell)<\begin{cases}
     aN & \text{if }\ell\leq n-b \\
     (a+1)N & \text{if }\ell>n-b
\end{cases}
    \]
    As we defined $w_\hess = w_\fin\cdot w_{-\vla}$, the left hand side is exactly $w_\hess(i_\ell) - w_\hess(j_\ell)$.

    To conclude, we must compute the rank of the affine space bundle $E_{w_\cell(\Delta,f)}$. We can decompose 
    \[
    E_{w_\cell(\Delta,f)} = \bigcup_{\substack{(\Delta',f')\colon \\w_\cell(\Delta',f')\in S_d^n\cdot w_\cell(\Delta,f)\cdot S_\mu/S_\mu}} C_{w_\cell(\Delta',f')}
    \]
    into left Iwahori cells, with $C_{w_\cell(\Delta',f')}$ being affine spaces of dimension $\codinv(\Delta',f')$ by Proposition \ref{prop: parahoric dimension}. This stratification is compatible with the affine paving of the Hessenberg variety $\cH(\ol\gamma,\fm_{(\Delta,f)})$ in the sense that $C_{w_\cell(\Delta',f')}$ is an affine space bundle of the same rank as $E_{w_\cell(\Delta',f')}$ over an affine cell of $\cH(\ol\gamma,\fm_{(\Delta,f)})$. Since any affine paving of $\cH(\ol\gamma,\fm_{(\Delta,f)})$ has at least one element which is a single point (corresponding to the trivial Bruhat cell), it follows that the rank of the affine space bundle $E_{w_\cell(\Delta',f')}$ is
    \[
    \min \big\{\codinv(\Delta',f')\mid w_\cell(\Delta',f')\in S_d^n\cdot w_\cell(\Delta,f) \big\}
    \]
    The result now follows from Proposition \ref{prop: affine Springer monodromy} combined with Lemmas \ref{lem: unique dominant Delta,f} and \ref{lem: dominant has minimal codinv}.
\end{proof}

\begin{cor}
\label{cor: Sd frob char for iwahori case}
    The $(q,t)$-bigraded $S_d$ Frobenius character $\Frob_{S_d}(\mathrm{gr}_{F}H_{*}^{BM}(\Fl_\gamma^+))$ is
    \[
\sum_{\substack{
     (\Delta,f)\colon \Delta\in I_{M,N}^+\text{ dominant,} \\
     \text{content }(1^N)}} q^{\area(\Delta)}t^{\codinv(\Delta,f)}\Xi_{\fm_{(\Delta,f)}}(Y,t)
    \]
\end{cor}

We now consider the case of $\Fl_\gamma^\mu$. There is a surjective map $\Fl_\gamma\to \Fl_\gamma^\mu$ induced by the projection $G_K/I\to G_K/\Pfrak_\mu$. Under this map, the preimage of the cell $E_w^\mu$ is the union $\cup_{w'\in S_d^n\backslash w\cdot S_\mu} E_{w'}$. Recall that $x_i\in \Ngen(\Delta)$ denotes the unique $N$ gen such that $x_i = i\mod N$. 
\begin{defn}
    For a labeled invariant set $(\Delta,f)$ of content $\mu$, we consider the parabolic $P_{(\Delta,f)}$ in $H = \GL_d^n$ whose Lie algebra is a sum of root spaces, and contains the root space for entry $(i+(\ell-1)d,j+(\ell-1)d)$ if and only if
\[
f(x_{i+(\ell-1)d}) - f(x_{j+(\ell-1)d})+N(\vla_{j+(\ell-1)d} - \vla_{i+(\ell-1)d})\leq 0.
\]
\end{defn}
The following is easy to verify.
\begin{lem}
    \label{lem: dominant implies P is standard}
    If $(\Delta,f)$ is a dominant labeled invariant set, then $P_{(\Delta,f)}$ is a standard parabolic of $H = \GL_d^n$ (which in particular contains an upper triangular Borel).
\end{lem}

\begin{defn}
\label{def: hessenberg for content mu}
    If $(\Delta,f)$ is a labeled invariant set of content $\mu$, let $w_\cell(\Delta',f')\in [w_\cell(\Delta,f)]$ be a representative of the coset, with $(\Delta,f')$ being a labeled invariant set of content $(1^N)$. We define the Hessenberg space $\fm_{(\Delta,f)}$ associated to $(\Delta,f)$ to be the orbit
    \[
   \fm_{(\Delta,f)} :=  P_{(\Delta,f)}\cdot \fm_{(\Delta,f')}
    \]
    That is, we set
    \[
\fm_{(\Delta,f)} := \sum_{w_\cell(\Delta',f')\in [w_\cell(\Delta,f)]} \fm_{(\Delta',f')}
    \]
    where the sum is over all elements of the coset $W^\aff/S_\mu$ which are of the form $w_\cell(\Delta',f')$ for $(\Delta',f')$ of content $(1^N)$. More explicitly, $\fm_{(\Delta,f)}$ is the sum of root spaces indexed by $(i_\ell,j_\ell)$ satisfying
    \[
    f(x_{i_\ell})-f(x_{j_\ell})+N(\vla_{j_\ell}-\vla_{i_\ell})\leq \begin{cases}
        aN, & \ell\leq n-b\\
        (a+1)N, & \ell>n-b
    \end{cases}.
    \]
\end{defn}

\begin{thm}
    \label{thm: GKM for content mu}
    Let $(\Delta,f)$ be a labeled invariant set of content $\mu$. The cell $E_{w_\cell(\Delta,f)}^\mu$ is an affine space bundle over the Hessenberg variety $\cH_{P_{(\Delta,f)}}(\ol\gamma,\fm_{(\Delta,f)})$. 
    
    Moreover, when $(\Delta,f)$ is chosen to be dominant, the graded $S_d$ Frobenius characteristic is 
    \[
\Frob_{S_d}(H^{BM}_*(E_{w_\cell(\Delta,f)}^\mu))=t^{\codinv(\Delta,f)}\cdot\Xi_{\fm_{(\Delta,f)}}(Y,t)/\Aut_t(\alpha_{(\Delta,f)})
    \]
    where $\alpha_{(\Delta,f)}$ is the tuple of compositions of $d$ which corresponds to $P_{(\Delta,f)}$.
\end{thm}
\begin{proof}
    That $\fm_{(\Delta,f)}$ is the relevant Hessenberg space follows from an identical argument to the proof of Theorem \ref{thm: main for Fl gamma}. We must justify why $P_{(\Delta,f)}$ agrees with the parabolic of Theorem \ref{thm: GKM}. The parabolic $P$ is determined by the grading on $\Lie(H)$ induced by $w_\cell(\Delta,f)\cdot(\varepsilon\check\rho_\mu) = \vla+\varepsilon (w_\fin^{-1}\cdot \check\rho_\mu)$, so that $\Lie(P)$ contains the root space associated to the entry $(i+(\ell-1)d,j+(\ell-1)d)$ whenever 
    \[
    \vla_{i+(\ell-1)d} - \vla_{j+(\ell-1)d} + \varepsilon((\check\rho_\mu)_{w_\fin(i+(\ell-1)d)} - (\check\rho_\mu)_{w_\fin(j+(\ell-1)d)})\geq 0
    \]
    which is equivalent to the condition for $P_{(\Delta,f)}$ since one can check
    \[
    (w_\fin^{-1}\cdot\check\rho_\mu)_i = (\check\rho_\mu)_{w_\fin(i)} = \left(\frac{1}{N}\sum_{j\geq 1}\mu_j\cdot j\right) - f(x_i).
    \]
    To conclude, we must compute the rank of the affine space bundle $E_{w_\cell(\Delta,f)}^\mu$, which we do with an identical argument to Theorem \ref{thm: main for Fl gamma}. Namely, the rank of the affine bundle is equal to the minimal dimension of an affine cell 
    \[
    \min \big\{\codinv(\Delta',f')\mid w_\cell(\Delta',f')\in S_d^n w_\cell(\Delta,f) S_\mu/S_\mu \big\},
    \]
    and the result follows from lemmas \ref{lem: unique dominant Delta,f} and \ref{lem: dominant has minimal codinv}.
\end{proof}

\begin{cor}
\label{cor: Sd frob char for parahoric case}
    The $(q,t)$-bigraded $S_d$ Frobenius character of $\gr_F H^{BM}_*(\Fl_\gamma^{\mu,+})$ is
    \[
\Frob_{S_d}(\mathrm{gr}_{F}\;H^{BM}_*(\Fl_\gamma^{\mu,+})) = \sum_{\substack{
     (\Delta,f)\colon \Delta\in I_{M,N}^+ \\
     \text{dominant,}\text{ content }\mu}}q^{\area(\Delta)}t^{\codinv(\Delta,f)}\frac{\Xi_{\fm_{(\Delta,f)}}(Y,t)}{\Aut_t(\alpha_{(\Delta,f)})}
    \]
\end{cor}

\subsection{Main results for affine Springer fibers}

We now finish the proof of the $S_N\times S_d$ Frobenius characteristic of $\gr_F H^{BM}_*(\Fl_\gamma^+)$.

\begin{thm}
    \label{thm: ASF main}
    The $(q,t)$-bigraded $S_N\times S_d$ Frobenius character of $H^{BM}_*(\Fl_\gamma^+)$ is
    \[
\Frob_{S_N\times S_d}\big(\mathrm{gr}_F H^{BM}_*(\Fl_\gamma^+) \big) = \sum_{\mu\vdash N}\sum_{\substack{
     (\Delta,f)\colon \Delta\in I_{M,N}^+ \\
     \text{dominant,}\text{ content }\mu}}q^{\area(\Delta)}t^{\codinv(\Delta,f)}\frac{\Xi_{\fm_{(\Delta,f)}}(Y,t)}{\Aut_t(\alpha_{(\Delta,f)})}m_\mu(X)
    \]
    where the variables $Y$ give $S_d$ Frobenius character and the variables $X$ give $S_N$ Frobenius character.
\end{thm}
\begin{proof}
    We pair both sides with $h_\mu(X)$ for $\mu$ a partition of $N$. On the left hand side, we have
    \begin{align*}
\left\langle \Frob_{S_N\times S_d}\big(\gr_F H^{BM}_*(\Fl_\gamma^+) \big),h_\mu(X)\right\rangle &= \Frob_{S_d}\big(\gr_F H^{BM}_*(\Fl_\gamma^{\mu,+}) \big) \\
&= \sum_{\substack{
     (\Delta,f)\colon \Delta\in I_{M,N}^+ \\
     \text{dominant,}\text{ content }\mu}} q^{\area(\Delta)}t^{\codinv(\Delta,f)}\frac{\Xi_{\fm_{(\Delta,f)}}(Y,t)}{\Aut_t(\alpha_{(\Delta,f)})}
    \end{align*}
    with the last equality following from Corollary \ref{cor: Sd frob char for parahoric case}. Since $\langle m_\lambda,h_\mu\rangle = \delta_{\lambda,\mu}$, this agrees with the right hand side, and so the result follows.
\end{proof}

We note that the formulas \eqref{eqn:equivariant BM of hessenberg} apply to compute equivariant Borel-Moore homology of $\Fl_\gamma^+$. Hence, we obtain the following corollary.
\begin{cor}
    \label{cor: ASF equivariant main}
    Let $T = C_H(\ol\gamma)$ be the diagonal torus in $H = \GL_d^n$. The $T$-equivariant $(q,t)$-bigraded $S_N\times S_d$ Frobenius character of $H^{BM}_*(\Fl_\gamma^+)$ is
    \[
\Frob_{S_N\times S_d}\big(\mathrm{gr}_F H^{T,BM}_*(\Fl_\gamma^+) \big) = \sum_{\mu\vdash N}\sum_{\substack{
     (\Delta,f)\colon \Delta\in I_{M,N}^+ \\
     \text{dominant,}\text{ content }\mu}}q^{\area(\Delta)}t^{\codinv(\Delta,f)}\frac{\Xi_{\fm_{(\Delta,f)}}\left[\frac{Y}{1-t^{-1}},t\right]}{\Aut_t(\alpha_{(\Delta,f)})}m_\mu(X).
    \]
\end{cor}

\part{Coherent side}
	Keep $M=dm$ and $N=dn\geq2$, where $d,m,n$ are positive integers and $\gcd(m,n)=1$, and set $c=M/N=m/n$. In this part, write $G=\SL_N$ and $\gfrak=\mathfrak{sl}_N$. Let $T\subset G$ be the diagonal torus, with Lie algebra 
	$\hfrak=\left\{(x_1,\ldots,x_N)\in\C^N\mid\sum_{i=1}^N x_i=0\right\}.$
	\section{Cherednik algebras and character $D$-modules}\label{sec:basics}
	\subsection{Rational Cherednik algebras}
	For $b\in \C$, the rational Cherednik algebra $\Hrm_b:=\Hrm_b({S_N},\hfrak)$ is the $\C$-algebra generated by $\hfrak$, $\hfrak^*$ and $S_N$ with relations
	\begin{align*}
		&[x,x']=[y,y']=0, \quad wxw^{-1}=w(x), \quad wyw^{-1}=w(y),\\
		&[y,x]=x(y)-\sum_{s\in S}b\langle \alpha_s,y\rangle\langle x,\alpha_s^\vee\rangle s,
	\end{align*}
	where $x,x'\in \hfrak^*$, $y,y'\in \hfrak$, $w\in {S_N}$ and $\langle-,-\rangle$ is the natural pairing. Here $S\subset S_N$ is the set of reflections, $\alpha_s$ and $\alpha_s^\vee$ are the root and coroot associated with $s$, normalized by $\langle\alpha_s,\alpha_s^\vee\rangle=2$.
	
	For an irreducible $S_N$-representation $\mu$, extend $\mu$ to $S(\hfrak)\rtimes S_N$ by letting $\hfrak$ act trivially. The Verma module $M_b(\mu)=\Hrm_b\otimes_{S(\hfrak)\rtimes {S_N}}\mu$ has a unique simple quotient $\Lrm_b(\mu)$.
	\subsection{Modules of minimal support}
	For $b\in\C$, let $\Ocal(\Hrm_b)$ denote the category of finitely generated $\Hrm_b$-modules on which $S(\hfrak)_+$ acts locally nilpotently \cite{ggor}. Here $S(\hfrak)_+$ is the ideal of positive-degree elements in the symmetric algebra $S(\hfrak)$. Its simple objects are $\Lrm_b(\mu)$, indexed by the irreducible $S_N$-representations.
	
	\begin{thm}[{\cite{Wilcox}, \cite[Theorem 2.2]{egl}}]
		\label{thm:wilcox}\label{thm:beg}
		The minimal support of a nonzero object of $\Ocal(\Hrm_c)$, as a $\C[\hfrak]$-module, is
		\begin{align*}
			X_{d,n}:=S_N\cdot\left\{
			(\underbrace{z_1,\ldots,z_1}_{n},\ldots,
			\underbrace{z_d,\ldots,z_d}_{n})
			\ \middle|\ \sum_{i=1}^d z_i=0\right\}.
		\end{align*}
		The simple objects with this support are $\Lrm_c(n\lambda)$, where $\lambda\vdash d$ and $n\lambda=(n\lambda_1,n\lambda_2,\ldots)$. The full subcategory of $\Ocal(\Hrm_c)$ supported on $X_{d,n}$ is semisimple and equivalent to $\Rep_\C(S_d)$, with $\Lrm_c(n\lambda)$ corresponding to the irreducible representation indexed by $\lambda$.
	\end{thm}
	
	\begin{rem}
		Since $\dim X_{d,n}=d-1$, the minimal-support simples are infinite-dimensional when $d>1$. When $d=1$, $\Lrm_c((N))$ is the unique finite-dimensional simple.
	\end{rem}
	
	\subsection{$D$-module lifts of minimal supported modules}
	\subsubsection{Functors on $D$-modules}
	For a smooth algebraic variety $X$, let $\Dcal_X$ be its sheaf of differential operators and put $\Dcal(X)=\Gamma(X,\Dcal_X)$. Recall that every coherent $\Dcal_X$-module $M$ admits a good filtration: an exhaustive increasing filtration $F_jM$ compatible with the order filtration on $\Dcal_X$ such that $\gr^F M$ is coherent over $\gr\Dcal_X\cong\Sym T_X$. We write $\wt{\gr}(M)$ for the corresponding coherent sheaf on $T^*X$; see \cite[Section 2.1]{htt}. The support $SS(M)$ of $\wt{\gr}(M)$ in $T^*X$ is  called the singular support of $M$, which is independent of the choice of a good filtration.
	
	For a morphism $f:X\to Y$ of smooth algebraic varieties, write $f_\dagger$, $f_!$, and $f^\dagger$ for direct image, direct image with proper support, and inverse image in the bounded derived categories of holonomic $D$-modules. In terms of the transfer modules of \cite[Section 1.3]{htt}, our conventions are
	$f_\dagger(-)=Rf_*\bigl(\Dcal_{Y\leftarrow X}\otimes_{\Dcal_X}^{\mathrm L}(-)\bigr),
	f_!(-) =\Dbb_Y f_\dagger\Dbb_X(-),
	f^\dagger(-)=\Dcal_{X\to Y}\otimes_{f^{-1}\Dcal_Y}^{\mathrm L}f^{-1}(-)$ 
	where $\Dbb_X$ denotes Verdier duality. 
	
	For a locally closed embedding $i:U\hookrightarrow X$ of smooth varieties and a holonomic $\Dcal_U$-module $\Mcal$, its minimal extension is 
	$i_{!*}\Mcal=\operatorname{im}\bigl(\mathcal H^0(i_!\Mcal)\rightarrow \mathcal H^0(i_\dagger\Mcal)\bigr).$ 
	An integrable connection on $U$ (i.e. {an $\Ocal_U$-coherent} $\Dcal_U$-module) is called clean if its minimal extension coincides with the extensions using $i_\dagger$ and $i_!$.
	
	\subsubsection{Induced cuspidal $D$-modules}\label{sec:cee}\label{sec:loc-sys}
	\providecommand{\Ebb}{\mathbb{E}}
	\providecommand{\Kbb}{\mathbb{K}}
	
	Let $\Ecal_c$ be the rank-one integrable connection on the regular nilpotent orbit $O_{(n)}\subset\mathfrak{sl}_n$ with monodromy $e^{2\pi i m/n}$, as in \cite{lusztig4}, \cite[Theorem 5.5]{finkelberg-ginzburg_mirabolic} and \cite[Section 2.5.2]{ma2}. Since $\gcd(m,n)=1$, $\Ecal_c$ is clean for the inclusion $j:O_{(n)}\hookrightarrow\mathfrak{sl}_n$, as no smaller nilpotent orbit supports its order-$n$ central character. Let $\Nbb_{m,n}$ be its minimal extension to $\mathfrak{sl}_n$, equipped with the Hodge module structure induced by the weight-zero variation of type $(0,0)$ on $\Ecal_c$. 
	
	Choose the standard parabolic $Q\subset G$ with Levi $L=S(\GL_n^d)$. Then
	\begin{align*}
		\lfrak:=&\operatorname{Lie}(L)
		=\mathfrak{sl}_n^{\oplus d}\oplus\zfrak,\\
		\zfrak:=&\mathfrak z(\lfrak)
		=\{(z_1I_n,\ldots,z_dI_n)\mid z_1+\cdots+z_d=0\}.
	\end{align*}
	We use the parabolic induction functor $\Ind_Q^G$ for equivariant Hodge modules from \cite[Sections 5.1-5.2]{BLM}. Concretely, for $\mathfrak q=\operatorname{Lie}(Q)$, it pulls an $L$-equivariant module back along the Levi projection $\mathfrak q\to\lfrak$, forms the associated $G$-equivariant module on $G\times^Q\mathfrak q$, and takes direct image along $[g,x]\mapsto\operatorname{Ad}(g)x$.
	
	\begin{defn}\label{def:induced-cuspidal}
		The character $D$-module of parameter $(M,N)$ is
		\begin{align*}
			\Nbb_{M,N}:=\Ind_Q^G\bigl(\Nbb_{m,n}^{\boxtimes d}\boxtimes\Ocal_{\zfrak}\bigr).
		\end{align*}
	\end{defn}
	
	Put
	\begin{align*}
		\zfrak^\circ
		&:=\{(z_1I_n,\ldots,z_dI_n)\in\zfrak
		\mid z_i\neq z_j\text{ for }i\neq j\}, \quad 
		\gfrak_d^\circ
		:=G\cdot(\zfrak^\circ+O_{(n)}^d),
	\end{align*}
	where $O_{(n)}^d$ is embedded block-diagonally in $\lfrak$. The map
	\begin{align*}
		q^\circ:\wt{\gfrak}_d^\circ
		:=G\times^L(\zfrak^\circ+O_{(n)}^d)
		&\rightarrow \gfrak_d^\circ,\quad 
		[g,x]\mapsto\Ad(g)x
	\end{align*}
	is a finite \'etale $S_d$-cover recording an ordering of the $d$ generalized eigenspaces.
	
	Let $\wt{\Ebb}_m^\circ$ be the connection on $\wt{\gfrak}_d^\circ$ induced from $\Ecal_c^{\boxtimes d}$, with the trivial connection along $\zfrak^\circ$. Permuting the blocks extends the equivariant structure to the normalizer of $L$ in $G$, giving an $S_d$-equivariant structure on $\wt{\Ebb}_m^\circ$ for the deck action. This is the construction in the proof of \cite[Theorem 6.1]{BLM}, applied to the Lie algebra as in \cite[Section 6.2]{BLM}. Set
	\begin{align*}
		\Ebb_m^\circ&:=q^\circ_*\wt{\Ebb}_m^\circ,\quad 
		\Ebb_{m,\tau}:=\Hom_{S_d}(V_\tau,\Ebb_m^\circ)
		\quad \forall \tau\vdash d,
	\end{align*}
	where $V_\tau$ is the irreducible $S_d$-module indexed by $\tau$. For a connection on a locally closed stratum of $\gfrak$, $\IC(\gfrak,-)$ denotes its minimal extension to the closure of the stratum, followed by closed pushforward to $\gfrak$.
	
	\begin{thm}[{\cite[Section 6.2]{BLM}}]
		\label{prop_Sd_action}
		The deck action induces an action of $S_d$ on $\Nbb_{M,N}$ as Hodge module automorphisms, and there is an $S_d$-equivariant decomposition
		\begin{align*}
			\Nbb_{M,N}
			\cong\bigoplus_{\tau\vdash d}
			V_\tau\otimes\IC(\gfrak,\Ebb_{m,\tau}).
		\end{align*}
		The IC summands carry their induced Hodge structures, and the multiplicity spaces $V_\tau$ have Hodge type $(0,0)$.
	\end{thm}
	
	\begin{proof}
		By the Lie algebra analogue of \cite[Proposition 6.5]{BLM}, $\Nbb_{M,N}$ is the intermediate extension of $\Ebb_m^\circ$. The deck action extends uniquely to the intermediate extension, and the finite-group isotypic projectors give the displayed decomposition, as in the proof of \cite[Theorem 6.1]{BLM}. These maps preserve the type-$(0,0)$ variations and hence are morphisms of Hodge modules.
	\end{proof}
	
	\begin{cor}\label{coro: gr(convolution)}
		Under the $S_d$-action from Theorem \ref{prop_Sd_action}, 
		\begin{align*}
			\gr(\Nbb_{M,N})
			\cong\bigoplus_{\tau\vdash d}
			V_\tau\otimes
			\gr\bigl(\IC(\gfrak,\Ebb_{m,\tau})\bigr).
		\end{align*}
	\end{cor}
	\begin{proof}
		The isotypic projectors are morphisms of Hodge modules and hence are strict for the Hodge filtration. Apply associated graded to the decomposition in Theorem \ref{prop_Sd_action}.
	\end{proof}
	
	For $\lambda=(\lambda_i)\vdash d$, let $O_{n\lambda}\subset\gfrak$ be the nilpotent orbit of Jordan type $n\lambda=(n\lambda_1,n\lambda_2,\ldots)$. Writing $b=\gcd(\lambda_i)$, one has $\pi_1(O_{n\lambda})\cong\Z/(nb)\Z$. Let $\Kbb_{m,\lambda}$ be the rank-one connection on $O_{n\lambda}$ corresponding to the character sending the standard generator to $e^{2\pi i m/n}$. 
	
    Throughout, we identify $\gfrak$ with $\gfrak^*$ by the Killing form and write $\Ffrak$ for Fourier transform on $\Dcal_\gfrak$-modules.
	
	\begin{prop}[{\cite[Proposition 6.7]{BLM}}]
		\label{prop:fourier-induced-cuspidal}
		One has
		\begin{align*}
			\Ffrak\bigl(\IC(\gfrak,\Ebb_{m,\tau})\bigr)
			\cong\IC(\gfrak,\Kbb_{m,\tau^t}),
		\end{align*}
		where $\tau^t$ denotes the transpose partition. Consequently, with the $S_d$-action transported by Fourier transform,
		\begin{align*}
			\Ffrak(\Nbb_{M,N})
			\cong\bigoplus_{\tau\vdash d}
			V_\tau\otimes\IC(\gfrak,\Kbb_{m,\tau^t}).
		\end{align*}
		These are isomorphisms of underlying $D$-modules.
	\end{prop}
	\begin{lem}\label{lem:cuspidal-support}
		The center $Z(G)$ acts on $\Nbb_{M,N}$ by $\zeta\mapsto\zeta^M$. Let $\gfrak_d:=\ol{\gfrak_d^\circ}$. For every $\tau\vdash d$ one has
		\begin{align}
			\label{formula:u}
			SS\bigl(\IC(\gfrak,\Ebb_{m,\tau})\bigr)
			\subset
			\bigl\{(x,y)\in\gfrak_d\times\ol{O}_{n\tau^t}\mid[x,y]=0\bigr\}_\red
		\end{align}
        where $_\red$ indicates taking reduced structure. 
	\end{lem}
	\begin{proof}
		On each $\Kbb_{m,\lambda}$, a generator of $Z(G)$ acts by $e^{2\pi i m/n}=e^{2\pi i M/N}$. Fourier transform preserves the $G$-action, so Proposition \ref{prop:fourier-induced-cuspidal} gives the asserted central character of $\Nbb_{M,N}$.
		
		The module $\IC(\gfrak,\Ebb_{m,\tau})$ is supported on $\gfrak_d$, and Proposition \ref{prop:fourier-induced-cuspidal} identifies its Fourier transform with $\IC(\gfrak,\Kbb_{m,\tau^t})$, supported on $\ol{O}_{n\tau^t}$. Since the original module is monodromic for scalar dilation (the Euler vector field acts locally finitely on its global sections), Fourier transform acts on its singular support by  $(x,y)\mapsto(y,-x)$. Thus its singular support lies in $\gfrak_d\times\ol{O}_{n\tau^t}$. Finally, $G$-equivariance imposes the moment-map equation $[x,y]=0$.
	\end{proof}
	\begin{rem}\label{rem:GLn-action}
		By Lemma \ref{lem:cuspidal-support}, the $G$-action on $\Nbb_{M,N}$ extends to $\GL_N$ by letting $zI_N$ act as $z^M$. Indeed, let $(g,z)\in G\times\C^*$ act by $z^M$ times the given action of $g$. The kernel of the map $G\times\C^*\to\GL_N$, $(g,z)\mapsto zg$, consists of $(\zeta I_N,\zeta^{-1})$ with $\zeta\in Z(G)$. Each such element acts trivially by the central-character formula, so the action descends to $\GL_N$. Replacing $m$ by $m+kn$ for an integer $k$ leaves the underlying $G$-equivariant Hodge module unchanged and replaces this $\GL_N$-action by its tensor product with $\det^k$, since $M$ changes to $M+kN$. 
	\end{rem}
	
	\subsection{The Calaque-Enriquez-Etingof construction}
	\label{sec:cee-construction}
	Write $\mathsf{N}_{M,N}:=\Gamma(\gfrak,\Nbb_{M,N})$. We use superscripts to indicate the rank: $\Hrm_b^{(r)}:=\Hrm_b(S_r,\hfrak_r)$, $\Arm_b^{(r)}:=\e^{(r)}\Hrm_b^{(r)}\e^{(r)}$, and $\Lrm_b^{(r)}(\mu)$ denotes the simple module with lowest $S_r$-type $\mu$. Here $\hfrak_r$ is the reflection representation of $S_r$ and $\e^{(r)}$ is its symmetrizing idempotent. With this notation, the module introduced in the introduction is
	\begin{align*}
		\Lrm_{M,N}:=\bigoplus_{\tau\vdash d}V_{\tau^t}\otimes\Lrm_{M/N}^{(N)}(n\tau).
	\end{align*}
	For a $G$-equivariant $\Dcal_\gfrak$-module $\Mcal$, put
	\begin{align*}
		\mathsf F_M^{(N)}(\Mcal)
		:=\bigl(\Gamma(\gfrak,\Mcal)\otimes(V^*)^{\otimes M}\bigr)^G,
		\quad V=\C^N.
	\end{align*}
	The Calaque-Enriquez-Etingof (CEE) construction equips this space with a Cherednik-algebra action \cite{cee}. We follow the formulation in \cite[Section 7.6]{egl}, which uses $V^*$ as above. In particular, the number of tensor factors is the rank of the resulting Cherednik algebra, and its parameter is $N/M=1/c$. 
	\begin{thm}[{\cite[Proposition 7.5]{egl}}]
		\label{thm:cee}
		The functor $\mathsf F_M^{(N)}$ takes values in $\Hrm_{N/M}^{(M)}\Mod$, with $S_M$ acting by permutation of the tensor factors. For each partition $\lambda\vdash d$,
		\begin{align*}
			\mathsf F_M^{(N)}
			\bigl(\Ffrak(\IC(\gfrak,\Kbb_{m,\lambda}))\bigr)
			\cong\Lrm_{N/M}^{(M)}(m\lambda).
		\end{align*}
		Consequently, there is an $S_d$-equivariant isomorphism of $\Hrm_{N/M}^{(M)}$-modules
		\begin{align}
			\label{iso: cee iso}
			\bigl(\mathsf{N}_{M,N}\otimes(V^*)^{\otimes M}\bigr)^G
			\cong\bigoplus_{\tau\vdash d}
			V_{\tau^t}\otimes\Lrm_{N/M}^{(M)}(m\tau)=\Lrm_{N,M}.
		\end{align}
	\end{thm}
	\begin{proof}
		In \cite[Section 7.6]{egl}, replace the matrix size by $N$ and the number of tensor factors by $M$. The orbital module $\IC(\gfrak,\Kbb_{m,\lambda})$ has support $\ol{O}_{n\lambda}$ and central character $M\bmod N$. The CEE correspondence therefore sends its Fourier transform to $\Lrm_{N/M}^{(M)}((M/d)\lambda)$. Moreover, $\Ffrak^2$ is pullback by $x\mapsto-x$, which preserves the isomorphism class of each of these orbital modules. Thus Proposition \ref{prop:fourier-induced-cuspidal} gives
		\begin{align*}
			\mathsf F_M^{(N)}\bigl(\IC(\gfrak,\Ebb_{m,\sigma})\bigr)
			\cong\Lrm_{N/M}^{(M)}(m\sigma^t).
		\end{align*}
		Apply $\mathsf F_M^{(N)}$ to the decomposition in Theorem \ref{prop_Sd_action} and reindex by $\tau=\sigma^t$.
	\end{proof}

	For $0\leq i<M$, the Koszul differential sends $f\otimes(v_1\wedge\cdots\wedge v_i)$ to $\sum_{j=1}^i(-1)^{j-1}fv_j\otimes(v_1\wedge\cdots\widehat{v_j}\cdots\wedge v_i)$, where $f\in\Sym^{M-i}(V^*)$ and $v_j\in V^*$. We take its target to be zero when $i=0$. It gives
	\begin{align*}
		\Hom_{S_M}\bigl(\wedge^i\hfrak_M,(V^*)^{\otimes M}\bigr)
		\cong\ker\left[
		\begin{aligned}
			&\Sym^{M-i}(V^*)\otimes\wedge^i(V^*)\\
			&\quad\rightarrow 
			\Sym^{M-i+1}(V^*)\otimes\wedge^{i-1}(V^*)
		\end{aligned}\right].
	\end{align*}
	Tensoring with $\mathsf{N}_{M,N}$, taking $G$-invariants, and using \eqref{iso: cee iso} yields
	\begin{align}
		\label{eqn: hook kernel}
		&\bigoplus_{\tau\vdash d}V_{\tau^t}\otimes
		\Hom_{S_M}\bigl(\wedge^i\hfrak_M,
		\Lrm_{1/c}^{(M)}(m\tau)\bigr)
		\cong\ker\left[
		\begin{aligned}
			&\bigl(\mathsf{N}_{M,N}\otimes\Sym^{M-i}(V^*)
			\otimes\wedge^i(V^*)\bigr)^G\\
			&\hspace{1cm}\rightarrow 
			\bigl(\mathsf{N}_{M,N}\otimes\Sym^{M-i+1}(V^*)
			\otimes\wedge^{i-1}(V^*)\bigr)^G
		\end{aligned}\right].
	\end{align}
	\begin{rem}
		Mirabolic $D$-modules provide a related approach to Cherednik algebras through quantum Hamiltonian reduction \cite{finkelberg-ginzburg_mirabolic,ganginzburg}, as discussed in \cite{cee}. 
	\end{rem}
	
	\section{Hodge associated graded of character $D$-modules}
	\label{sec:coherent-commuting}
	Following \cite[Sections 3-4]{ma2}, we construct two Koszul complexes whose derived pushforwards give the Hodge associated graded sheaves of $\Nbb_{M,N}$ and its trivial $S_d$-summand $\IC(\gfrak,\Ebb_{m,(d)})$. Throughout this section, $\wt\gr$ denotes the coherent sheaf associated with the Hodge filtration, and grading shifts are defined by $(\Fcal\langle s\rangle)_k=\Fcal_{k-s}$.
	For a vector space $W$ with linear coordinates $x_j$, our Euler operator on $D$-modules is
	\begin{align}
		H_W:=\sum_jx_j\partial_{x_j}+\frac{\dim W}{2}.
        \label{eqn:euler}
	\end{align}
	This convention is compatible with direct sums: $H_{W_1\oplus W_2}=H_{W_1}+H_{W_2}$. We record Euler and filtration degrees by $q,h$, respectively; Euler characters may involve half-integral powers of $q$.
	
	\subsection{The induced $DG$-module}
	\label{sec:cuspidal-dg-modules}
	We adapt the construction of \cite[Section 3.7]{ma2} to the Levi $L=S(\GL_n^d)$ of Definition \ref{def:induced-cuspidal}. Let $B\subset Q$ be the upper-triangular Borel, $\Bcal=G/B$, $\bfrak=\operatorname{Lie}(B)$, $\nfrak=[\bfrak,\bfrak]$, and $\hfrak_L=\hfrak\cap[\lfrak,\lfrak]$. With $\zfrak=\mathfrak z(\lfrak)$, set $\bfrak_L:=\hfrak_L\oplus\nfrak$, 
	\begin{align*}
\mfrak_L:=\bfrak_L\oplus\bigoplus_{\substack{1\leq i<N\\n\nmid i}}\C E_{i+1,i},\quad 
\Xfrak_{M,N}:=G\times^B\bigl((\zfrak+\nfrak)\times\mfrak_L\bigr).
	\end{align*}
	Here $E_{ij}$ is the matrix unit with its only nonzero entry equal to $1$ in position $(i,j)$.  In other words, the negative simple-root spaces in $\mfrak_L$ lie within the size-$n$ Levi blocks. Since $[\zfrak+\nfrak,\mfrak_L]\subset\bfrak_L$, the commutator defines the first of the maps
	\begin{align*}
		q_L:\Xfrak_{M,N}&\rightarrow  G\times^B\bfrak_L,
		&[g,x,y]&\mapsto[g,[x,y]],\\
		a_L:\Xfrak_{M,N}&\rightarrow  T^*\gfrak,
		&[g,x,y]&\mapsto(\Ad(g)x,\Ad(g)y).
	\end{align*}
	Let $pr:\zfrak+\nfrak\to\nfrak\cap\lfrak$ be the projection taking the nilpotent parts of the diagonal blocks, and put
$U:=pr^{-1}(O_{(n)}^d).$ 
The $B$-equivariant pullback $pr^*(\Ecal_c^{\boxtimes d})$ on $U$ induces a connection on $G\times^B U$. Let $\Mbb_{M,N}$ be its minimal extension to $G\times^B\gfrak$, equipped with the Hodge module structure induced by the variation of type $(0,0)$. By construction, the direct image of $\Mbb_{M,N}$ along $[g,x]\mapsto\Ad(g)x$ is $\Nbb_{M,N}$.

	Let $\pi:\Xfrak_{M,N}\to\Bcal$ be the projection and $\ul{\bfrak_L^*}=G\times^B\bfrak_L^*$ the associated vector bundle. Put
	\[\lceil\mu_c\rceil:=\bigl(\lceil ic\rceil-\lceil(i-1)c\rceil\bigr)_{i=1}^N.\]
	and define the DG module 
	\begin{align}
		\label{eqn: Ae}\Acal_{M,N}&:=\left(\Ocal_{\Xfrak_{M,N}}\otimes\pi^*\bigl(\wedge^\bullet\ul{\bfrak_L^*}\otimes\Lcal_{\lceil\mu_c\rceil}\bigr),\iota_{[x,y]}\right).
	\end{align}
	Here $\iota_{[x,y]}$ is contraction with the commutator, and $\Lcal_\lambda=G\times^B\C_\lambda$, where $\C_\lambda$ is the one-dimensional $B$-module of weight $\lambda$ with trivial unipotent action. The exterior term $\wedge^j$ has cohomological degree $-j$.
	\subsection{The rectangular $DG$-module}
	Let $P\subset G$ be the standard parabolic with Levi $S(\GL_d^n)$ and nilradical $\ufrak_{P}$. It stabilizes a flag of dimensions $d,2d,\ldots,(n-1)d$. In $n\times n$ block notation, with each block of size $d$, put
	\begin{align*}
		\mfrak_{P}&:=\{x=(x_{ij})\in\gfrak\mid x_{ij}=0\text{ if }i>j+1\}.
	\end{align*}
	The action map $G\times^{P}\ufrak_{P}\to\ol O_{(n^d)}$ resolves the orbit closure and is an isomorphism over $O_{(n^d)}$. Let $\Mbb_P$ be the $P$-equivariant minimal extension to $\ufrak_{P}$ of the restriction of $\Kbb_{m,(1^d)}$ to $O_{(n^d)}\cap\ufrak_{P}$, with the Hodge structure induced by its variation of type $(0,0)$. Take the induced $G$-equivariant Hodge module on $G\times^P\ufrak_P$, and take its pushforward $\Mbb_{M,N}^{(d)}$ on $G\times^P\gfrak$ along the closed embedding. By construction, the direct image of $\Mbb_{M,N}^{(d)}$ along $[g,x]\mapsto\Ad(g)x$ is $\IC(\gfrak,\Kbb_{m,(1^d)})$.
	
	The module $\Mbb_P$ has an explicit description in blocks. For $0\leq i\leq n$, put
	\begin{align*}
		k_i:=\lceil ic\rceil ,
		\quad\alpha_i:=k_i-ic.
	\end{align*}
	For $1\leq i<n$, set $A_i=x_{i,i+1}$; then $0<\alpha_i<1$, whereas $\alpha_0=\alpha_n=0$. Each adjacent block $A_i$ carries the determinant connection $\Ocal_{\mathfrak{gl}_d}[\det^{-1}]\det^{-\alpha_i}$. Here $\Ocal[\det^{-1}]\det^\alpha$ denotes the rank-one connection with a formal generator $e$ satisfying $\partial e=\alpha\,\partial(\det)/\det\cdot e$ for every vector field $\partial$ on $\{\det\ne0\}$.
	Each block $x_{ij}$, with $j>i+1$, carries the trivial connection. 
	For $0<\alpha_i<1$, the determinant connection on $\{\det\neq0\}$ is clean for the inclusion $\GL_d\hookrightarrow\mathfrak{gl}_d$ by \cite[Remark 5.3]{lorincz-yang}. On the open set where every $A_i$ is invertible, their tensor product is generated by $\bigotimes_{i=1}^{n-1}\det(A_i)^{-\alpha_i}$. The equivariant structure inherited from $\Kbb_{m,(1^d)}$ gives this generator the Levi character $\prod_{i=1}^n\det(g_i)^{k_i-k_{i-1}}$ for $(g_1,\ldots,g_n)\in S(\GL_d^n)$.
	
	Let $i_{P}:\ufrak_{P}\hookrightarrow\gfrak$ be the closed embedding. Using the monodromic Hodge-module Fourier transform of \cite[Theorem 1.4]{chen-dirks}, composed with pullback by $x\mapsto-x$ on the target to obtain our sign convention for $\Ffrak$, define $\Cscr_{P}^\vee:=\wt\gr(\Mbb_P^\vee)$ for 
	\begin{align}
		\label{eqn:rectangular-fourier-coefficient}
		\Mbb_P^\vee
		&:=\Ffrak\bigl((i_{P})_\dagger\Mbb_P\bigr)
		\left(R-\dim\gfrak+\left\lfloor\frac{N-d}{2}\right\rfloor\right).
	\end{align}
	Parentheses $(-)$ denote the Tate twist, with $F_p(\Mcal(t))=F_{p-t}\Mcal$ and $R:=\dim(G/P)=d^2\binom n2.$ Let $(\Mbb_{M,N}^{(d)})^\vee$ be the $G$-equivariant Hodge module on $G\times^P\gfrak$ associated with $\Mbb_P^\vee$; thus it is the Fourier transform along the fibers of $\Mbb_{M,N}^{(d)}$, with the displayed Tate twist. By Proposition \ref{prop:fourier-induced-cuspidal}, its proper direct image along $[g,x]\mapsto\Ad(g)x$ is $\IC(\gfrak,\Ebb_{m,(d)})$, with the twist normalizing the Hodge type of $\Ebb_{m,(d)}$ to $(0,0)$.
	
	The embedding $\ufrak_P\to \gfrak$ induces the maps 
	\[T^*\ufrak_P \xleftarrow{\rho} \ufrak_P\times_\gfrak T^*\gfrak\xrightarrow{\varpi}T^*\gfrak\]
	where $\rho$ is the co-differential. The filtered closed-embedding formula of \cite[Section 1.3]{htt} gives
	\begin{align*}
		\wt\gr\bigl((i_{P})_\dagger\Mbb_P\bigr)
		\cong \varpi_*\left(\rho^*\wt\gr(\Mbb_P)\otimes
		\omega_{\ufrak_{P}/\gfrak}\right)\langle\dim\gfrak-R\rangle,
	\end{align*}
	where $\omega_{\ufrak_P/\gfrak}:=\omega_{\ufrak_P}\otimes i_P^*\omega_\gfrak^{-1}$ is the relative canonical bundle. 
    
	Let $\sigma:T^*\gfrak\to T^*\gfrak$ be the map  $\sigma(x,y)=(y,-x)$. By \cite[Theorem 1.4]{chen-dirks}, there is an isomorphism of $P$-equivariant coherent sheaves, with the regrading specified in \eqref{eqn:rectangular-fourier-grading} below,
	\begin{align*}
		\Cscr_P^\vee
		\cong \sigma_*\varpi_*
		\left(\rho^*\wt\gr(\Mbb_P)\otimes\omega_{\ufrak_P/\gfrak}\right).
	\end{align*}
	Take $W=\gfrak$ in \eqref{eqn:euler}. Let $(\beta,k)$ denote the $H_\gfrak$-eigenvalue and Hodge degree of a homogeneous global section of $\wt\gr((i_P)_\dagger\Mbb_P)$. Fourier transform sends $H_\gfrak$ to $-H_\gfrak$. In these degrees, the regrading above is
	\begin{align}
		(\beta,k)\longmapsto
		\left(-\beta,\ k+R+\left\lfloor\frac{N-d}{2}\right\rfloor+\left\lceil\beta-\frac{\dim\gfrak}{2}\right\rceil\right).
		\label{eqn:rectangular-fourier-grading}
	\end{align}
Here $k$ includes the shift $\dim\gfrak-R$ from the closed embedding, and the target grading includes the Tate twist in \eqref{eqn:rectangular-fourier-coefficient}.

	The sheaf $\Cscr_{P}^\vee$ is supported on the scheme
	\begin{align*}
		Z_{P}^0:=\{(x,y)\in\mfrak_{P}\times\ufrak_{P}\mid
		x_{i+1,i}y_{i,i+1}=y_{i,i+1}x_{i+1,i}=0\ (1\leq i<n)\}.
	\end{align*}
	Indeed, these equations are the symbols of the infinitesimal left and right $\gl_d$ actions on the adjacent size-$d$ blocks, and hence they annihilate the associated graded. 
    They also imply $[x,y]\in\ufrak_{P}$, so contraction with $[x,y]$ defines a Koszul differential with values in $\ufrak_{P}$.
	
	Set $Z_{P}=G\times^{P}Z_{P}^0$, let $\pi_{P}:Z_{P}\to G/P$ be the projection, and write $\wt{\Cscr}_{P}^\vee$ for the induction of $\Cscr_{P}^\vee$ to $Z_{P}$. With $\ul{\ufrak_{P}^*}=G\times^{P}\ufrak_{P}^*$, define
	\begin{align}
		\label{eqn:rectangular-dg-module}
		\Acal_{M,N}^{(d)}:=\left(
		\wt{\Cscr}_{P}^\vee\otimes
		\pi_{P}^*\bigl(\omega_{G/P}\otimes\wedge^\bullet\ul{\ufrak_{P}^*}\bigr),
		\iota_{[x,y]}\right).
	\end{align}
	Put 
	  $a_{P}:Z_{P}\rightarrow  T^*\gfrak,
		\quad [g,x,y]\mapsto(\Ad(g)x,\Ad(g)y).$
	
	\subsection{Associated graded of character $D$-modules}
	\begin{thm}\label{thm:gr N Ac}
		There are $\GL_N\times\C^*$-equivariant isomorphisms of coherent sheaves on $T^*\gfrak$, where the $\GL_N$-actions extend the $G$-actions by letting $zI_N$ act as $z^M$, as in Remark \ref{rem:GLn-action}, and $\C^*$ records the Hodge grading:
		\begin{enumerate}
			\item $				\wt\gr(\Nbb_{M,N})\cong R(a_L)_*\Acal_{M,N}\langle N-d\rangle.$
			\item $\wt\gr\bigl(\IC(\gfrak,\Ebb_{m,(d)})\bigr)
			\cong R(a_{P})_*\Acal_{M,N}^{(d)}\langle-R\rangle.$ 
		\end{enumerate}
		Both derived pushforwards are hence concentrated in cohomological degree zero.
	\end{thm}
	\begin{proof}
		The arguments of \cite[Theorem 3.4, Proposition 3.5, Theorem 3.8]{ma2} apply to $\Mbb_{M,N}$ and $(\Mbb_{M,N}^{(d)})^\vee$, whose proper direct images are $\Nbb_{M,N}$ and $\IC(\gfrak,\Ebb_{m,(d)})$, respectively.
	\end{proof}

	\section{Elliptic Hall algebras and Hilbert schemes}\label{sec:coherent-characters}
	\subsection{Bigradings}
	Put $A=(\C^*)^2$. On $\gfrak\times\gfrak$ we use the action
	\begin{align}
		\label{eq: A action}
		(u,v)\cdot(x,y)=(u^{-1}x,v^{-1}y),
	\end{align}
	so that the coordinate functions on the two factors have weights $q$ and $t$, respectively. Write $F^{\Hrm}$ for the Hodge filtration and normalize it by
	\begin{align*}
		F_j^0\mathsf{N}_{M,N}:=F_{j+N-d}^{\Hrm}\mathsf{N}_{M,N},
		\quad \gr^0\mathsf{N}_{M,N}:=\gr^{F^0}\mathsf{N}_{M,N}.
	\end{align*}
	
	By Theorem \ref{thm:gr N Ac}(1), $\gr^0\mathsf{N}_{M,N}$ corresponds to the coherent sheaf $R(a_L)_*\Acal_{M,N}$ with its geometric $A$-action. The CEE construction gives
	\begin{align*}
		\Lrm_{N,M}
		&\cong\mathsf F_M^{(N)}(\Nbb_{M,N})
		\cong\bigoplus_{\tau\vdash d}
		V_{\tau^t}\otimes\Lrm_{N/M}^{(M)}(m\tau),\\
		F_j^0\Lrm_{N,M}
		&:=\bigl(F_j^0\mathsf{N}_{M,N}\otimes(V^*)^{\otimes M}\bigr)^G.
	\end{align*}
	Each simple $\Lrm_{N/M}^{(M)}(m\tau)$ inherits a filtration by applying $\Hom_{S_d}(V_{\tau^t},-)$ to $F^0\Lrm_{N,M}$. Exchanging $M,N$ defines $F^0\Lrm_{M,N}$ from $\Nbb_{N,M}$, with Hodge shift $M-d$. 
	
	The dilation vector field $\sum_jx_j\partial_{x_j}$ preserves the Hodge filtration on the monodromic module $\mathsf{N}_{M,N}$, so $H_\gfrak$ acts on $\gr^0\mathsf{N}_{M,N}$. Under the CEE isomorphism \eqref{iso: cee iso}, it induces the  Euler element in $\Hrm_{N/M}^{(M)}$
	\begin{align*}
		\hrm_{N/M}^{(M)}:=\frac12\sum_{a=1}^{M-1}(x_ay_a+y_ax_a)
	\end{align*}
	where $x_a\in\hfrak_M^*$ and $y_a\in\hfrak_M$ are dual bases. Restriction to the subgroup $\{(u,u^{-1})\}\subset A$ sends a geometric weight $q^it^j$ to $q^{i-j}$, while its filtration degree remains $j$. Thus characters in the $H_\gfrak$ and $F^0$ gradings are obtained by substituting $t=h/q$ and applying a uniform Euler shift. Theorem \ref{thm:ch Lc Q} determines this shift; we record the resulting convention in \eqref{eq:geometric-Euler-conversion}.
	
	\subsection{Bigraded characters}
	Let $k=\C(q,t)$ and let $\Lambda$ be the ring of symmetric functions over $k$. The polynomial representation below takes values in symmetric functions in $X$. We use the positive half $\Ecal^+$ of the elliptic Hall algebra, its elements $Q_{u,v}$, and its polynomial representation in the conventions of \cite[Section 4.4]{ma2}. 

    For symmetric-group representations, $\ch$ denotes the same Frobenius characteristic as $\Frob$ in Part I, with the torus gradings included; for $\GL_r$-representations, $\ch$ denotes the ordinary torus character. For both $S_N$- and $S_M$-characters, the Frobenius alphabet is $X$. Write $\omega_X$ for the involution sending $s_\lambda(X)$ to $s_{\lambda^t}(X)$, where $\lambda^t$ is the transpose partition.
	\begin{thm}\label{thm:ch Lc Q}
		With $q$ recording the $\hrm_{M/N}^{(N)}$-grading and $h$ the filtration degree, we have
		\begin{align*}
			q^{-(d-1)/2}(1-q)^{d-1}
			\ch_{S_N\times A}(\gr^0\Lrm_{M,N})
			=\omega_X(Q_{m,n}^{\,d}\cdot1)(X;q,h/q).
		\end{align*}
	\end{thm}
	\begin{proof}
		We compute $\Lrm_{N,M}$ from $\Nbb_{M,N}$ and then exchange $M,N$. An analogous proof of \cite[Theorem 4.6]{ma2} applies to show 
		\begin{align*}
			\ch_{S_M\times A}(\gr^0\Lrm_{N,M})
			=\frac{q^b}{(1-q)^{d-1}}
			\omega_X(Q_{n,m}^{\,d}\cdot1)(X;q,h/q)
		\end{align*}
		for a scalar shift $b$. Here the trivial extension to $\zfrak(\lfrak)$ contribute $(1-q)^{-(d-1)}$.

		To determine $b$, set $h=1$. Write $\kappa(\lambda)=\sum_{(i,j)\in\lambda}(j-i)$ and let $\chi^\lambda((m^d))$ denote the character of $V_\lambda$ on cycle type $(m^d)$. The operator formulas of \cite[Theorems I.2 and I.6]{glwx} gives 
		\begin{align}
			\left.\omega_X(Q_{n,m}^{\,d}\cdot1)\right|_{t=q^{-1}}
			=(1-q)^d q^{d(m-1)/2}
			\sum_{\lambda\vdash M}
			\chi^\lambda((m^d))q^{-(n/m)\kappa(\lambda)}
			s_\lambda\left[\frac{X}{1-q}\right]. 
			\label{formula:pmn-t-inverse-q}
		\end{align}
		  By \cite[Theorem 1.4]{egl} and $\sum_{\tau\vdash d}\dim V_\tau\,s_\tau[p_m]=p_m^d$, the Euler-graded character of $\Lrm_{N,M}$ is the same sum as \eqref{formula:pmn-t-inverse-q} but multiplied $(1-q)q^{(M-1)/2}$. Hence $b+d(m-1)/2=(M-1)/2$, giving $b=(d-1)/2$. Exchanging $M,N$ gives the stated formula.
	\end{proof}
	
	\begin{rem}\label{rem:grading-conversion}
		Theorem \ref{thm:ch Lc Q}, with $M,N$ exchanged, shows that a vector in $\gr^0\mathsf{N}_{M,N}$ of geometric bidegree $(i,j)$ has $(H_\gfrak,F^0)$-bidegree
		\begin{align}
			(i,j)\longmapsto\left(i-j-(N-d)+\frac{d-1}{2},j\right).
			\label{eq:geometric-Euler-conversion}
		\end{align}
		Thus the character conversion is $t=h/q$ followed by multiplication by $q^{-(N-d)+(d-1)/2}$. After tensoring with $\Ocal_\C$ to pass from $\mathfrak{sl}_N$ to $\mathfrak{gl}_N$, use $H_{\mathfrak{gl}_N}=H_{\mathfrak{sl}_N}+H_\C$; the scalar factor has geometric character $(1-q)^{-1}$ and Euler character $q^{1/2}(1-q)^{-1}$. The corresponding multiplier for the full matrix space is therefore $q^{-(N-d)+d/2}$. We use these conversions throughout.
	\end{rem}
	
	\subsection{Coherent sheaves on Hilbert schemes of points}
    We give a realization of $Q_{n,m}^d$ in terms of coherent sheaves on the Hilbert schemes that enjoy cohomology vanishing properties. 
	
	Let $\ol{\Hilb}^N$ be the moduli of colength $N$ ideals in $\C[x,y]$, and let $\xi:\ol{\Hilb}^N\to\Sym^N(\C^2)$ be the Hilbert-Chow map. Take 
	\begin{align*}
		\Hilb^N:=\xi^{-1}\bigl((\hfrak_N\oplus\hfrak_N)/S_N\bigr).
	\end{align*}
	
	Let $\Cfrak\subset\gfrak^2$ be the reduced commuting scheme. Set
	\begin{align*}
		\wt{\Hilb}^N:=\bigl\{(x,y,v)\in\Cfrak\times V\mid\C[x,y]v=V\bigr\}.
	\end{align*}
	on which the diagonal action of $\GL_N$ is free and the resulting quotient is $\Hilb^N$. 
	The action of $A$ in \eqref{eq: A action}, acting trivially on $v$, descends to $\Hilb^N$.
	
	\subsubsection{The BKR equivalence}
	The isospectral Hilbert scheme is the reduced fiber product
	\begin{align*}
		\IHilb^N:=\bigl(\Hilb^N\times_{(\hfrak_N\oplus\hfrak_N)/S_N}(\hfrak_N\oplus\hfrak_N)\bigr)_\red.
	\end{align*}
	Let $\alpha$ and $\beta$ be its projections to $\Hilb^N$ and $\hfrak_N\oplus\hfrak_N$, respectively. The Procesi bundle $\Pcal:=\alpha_*\Ocal_{\IHilb^N}$ is an $S_N\times A$-equivariant vector bundle of rank $N!$. The Bridgeland-King-Reid equivalence \cite{BKR,haiman_macdonald} gives
	\begin{align*}
		\Phi_N:D_A^b\Coh(\Hilb^N)
		&\xrightarrow{\sim}D_{S_N\times A}^b\Coh(\hfrak_N\oplus\hfrak_N),\notag\\
		\Fcal&\mapsto R\beta_*L\alpha^*\Fcal.
	\end{align*}
	By the projection formula,
	\begin{align*}
		R\Gamma(\hfrak_N\oplus\hfrak_N,\Phi_N(\Fcal))
		\cong R\Gamma(\Hilb^N,\Pcal\otimes^L\Fcal).
	\end{align*}
	Composing the $K$-theory of $\Phi_N$ and the Frobenius characteristic yields an isomorphism
	\begin{align}
		\kappa:\bigoplus_{r\geq0}K^A(\Hilb^r)_\loc
		&\xrightarrow{\sim}\Lambda,\quad
		[\Fcal]\mapsto\ch_{S_r\times A}R\Gamma(\Hilb^r,\Pcal\otimes^L\Fcal).
		\label{iso: K Lambda}
	\end{align}
	\subsubsection{The induced coherent sheaf}
	\begin{prop}\label{prop: support gr Nc}
		The coherent sheaf $\wt\gr^0(\Nbb_{M,N})$ is scheme-theoretically supported on $\Cfrak$, and $F_{-1}^0\mathsf{N}_{M,N}=0$. 
	\end{prop}
\begin{proof}
Put $R=\C[T^*\gfrak]$ and $E=\gr^0\mathsf{N}_{M,N}$, and let $I\subset R$ be the ideal generated by the entries of $[x,y]$. The infinitesimal $G$-action preserves the Hodge filtration, so its principal symbols annihilate $E$; hence $IE=0$. Since $\Nbb_{M,N}$ underlies a Hodge module, $E$ is Cohen--Macaulay by \cite[Lemme 5.1.13]{saito-hodge}. In particular, its associated primes are the generic points of components in its supports, with no embedded associated primes.

By Lemma \ref{lem:cuspidal-support}, $\Supp(E)$ is contained in the locus of commuting pairs $(x,y)$ with $y$ nilpotent. For each nilpotent orbit $\Ocal\subset\gfrak$, the variety
$\bigl\{(x,y)\in\gfrak\times \Ocal \mid[x,y]=0\bigr\}$ 
is isomorphic to the conormal bundle to $\Ocal$, and is irreducible of dimension $\dim\gfrak$. There are finitely many nilpotent orbits. Since $E$ is the associated graded of a holonomic $D$-module, every component of $\Supp(E)$ has dimension $\dim\gfrak$ and therefore equals the closure of one of these conormal bundles.

Every such component meets the open set $U=\gfrak^{\mathrm{reg}}\times\gfrak$, where $\gfrak^{\mathrm{reg}}$ consists of matrices whose centralizer in $\gfrak$ has dimension $N-1$. Indeed, write a nilpotent matrix $y$ in Jordan blocks $y_1,\ldots,y_s$. The block-diagonal matrix with blocks $a_iI+y_i$, for pairwise distinct $a_i$, commutes with $y$ and is regular.

On $U$, the scheme defined by $I$ is the total space of the rank-$(N-1)$ vector bundle $\ker(\operatorname{ad}x)$ over $\gfrak^{\mathrm{reg}}$, so it is reduced. Thus $I\mathcal O_U=\sqrt I\,\mathcal O_U$, and $IE=0$ implies $(\sqrt I\,E)|_U=0$. But $\sqrt I\,E$ is a submodule of $E$, so all its associated primes belong to $\operatorname{Ass}_R(E)$. Every such prime is the generic point of a support component meeting $U$. Consequently $\sqrt I\,E=0$, proving scheme-theoretic support on the reduced scheme $\Cfrak$.

By Theorem \ref{thm:gr N Ac}(1), we have $\wt\gr^0(\Nbb_{M,N})\cong R(a_L)_*\Acal_{M,N}$. The complex $\Acal_{M,N}$ has no negative internal degrees: the base coordinates have degree zero, the cotangent coordinates and exterior generators have degree one, and the twisting line has degree zero. Hence $\gr_j^0\mathsf{N}_{M,N}=0$ for $j<0$. Since the filtration is bounded below, $F_{-1}^0\mathsf{N}_{M,N}=0$.
\end{proof}
	 Take the natural maps
	$p:\wt{\Hilb}^N\rightarrow \Cfrak$ and 
		$\varpi:\wt{\Hilb}^N\rightarrow \Hilb^N$
	that forget the cyclic vector and take the quotient, respectively. Proposition \ref{prop: support gr Nc} justifies the following definition. 
	\begin{defn}\label{def:Fc}
		Let $\Scal_{M,N}$ be the coherent sheaf on $\Cfrak$ whose pushforward to $T^*\gfrak$ is $\wt\gr^0(\Nbb_{M,N})$. Define its descent to the Hilbert scheme by
		\begin{align*}
			\Fcal_{M,N}:=(\varpi_*p^*\Scal_{M,N})^{\GL_N}.
		\end{align*}
		We use the $\GL_N$-lift of Remark \ref{rem:GLn-action}. The $A$-action on the Koszul complex $\Acal_{M,N}$ \eqref{eqn: Ae} commutes with this lift and hence descends to $\Fcal_{M,N}$.
	\end{defn}
	
	Define
	\begin{align*}
		\Sigma_{n,d}:=\left\{\sum_{i=1}^d n(a_i,0)\ \middle|\ a_i\in\C,\ \sum_i a_i=0\right\}
		\subset(\hfrak_N\oplus\hfrak_N)/S_N.
	\end{align*}
    Then $\Supp(\Fcal_{M,N})\subset\xi^{-1}(\Sigma_{n,d}).$ 
	\begin{prop}\label{prop:Fc and P}
		One has
		\begin{align*}
			\kappa([\Fcal_{M,N}])
			=\frac{q^{N-d}}{(1-q)^{d-1}}(Q_{m+n,n}^{\,d}\cdot1).
		\end{align*}
	\end{prop}
	\begin{proof}
		Apply the proof of \cite[Proposition 5.4]{ma2} using Theorem \ref{thm:ch Lc Q}.
	\end{proof}

	\subsubsection{Vanishing of higher cohomology}
	\label{sec:cohomology-vanishing}
	Let $\Pcal^\vee$ be the dual of the Procesi bundle. For a simple $S_N$-representation $W$, write
	\begin{align*}
		\Pcal^\vee_W:=\Hom_{S_N}(W,\Pcal^\vee).
	\end{align*}
	The proof of \cite[Theorem 6.4]{ma2} can be adapted verbatim to imply the following result. 
	\begin{thm}\label{thm: vanishing}
		For $0\leq i\leq N-1$ and $j>0$, one has
		\begin{align*}
			H^j\bigl(\Hilb^N,\Fcal_{M,N}\otimes\Pcal^\vee_{\wedge^i\hfrak_N}\bigr)=0.
		\end{align*}
		Moreover, there is an $S_d$-equivariant isomorphism
		\begin{align*}
			H^0\bigl(\Hilb^N,\Fcal_{M,N}\otimes\Pcal^\vee_{\wedge^i\hfrak_N}\bigr)
			\cong\Hom_{S_M}\bigl(\wedge^i\hfrak_M,\gr^0\Lrm_{N,M}\bigr).
		\end{align*}
	\end{thm}
    More generally, we expect that 
    \begin{conj}For $j>0$,  
        $H^j\bigl(\Hilb^N,\Fcal_{M,N}\otimes\Pcal^\vee\bigr)=0$  
        and there is an $S_N\times S_d$-equivariant isomorphism
        \[H^0\bigl(\Hilb^N,\Fcal_{M,N}\otimes\Pcal^\vee\bigr)\cong 
			\gr^0\Lrm_{M,N}.\]
    \end{conj}
	\section{The $S_d$-action and refined characters}
	\label{sec:Sd-refinement}
	We now go back to the $S_d$-action of Theorem \ref{prop_Sd_action} and refine the character formula in Theorem \ref{thm:ch Lc Q}. The approach is as follows: for the fixed primitive pair $(m,n)$, we let the number of blocks vary. Parabolic induction defines an algebra homomorphism from symmetric functions to the localized $K$ theory, as shown in Proposition \ref{prop:Psi}. This reduces the calculation to the trivial summands in Proposition \ref{prop:generator}.  We compute the character of the trivial summand in Proposition \ref{prop: char triv part} using Theorem \ref{thm:gr N Ac} and compare it with $h_r[(t-1)X^{m,n}]$ in Theorem \ref{thm:Sd char}. This computes the $S_N\times S_d\times A$ character of $\Lrm_{M,N}$.
	
	Throughout this section, coherent sheaves carry the geometric $q,t$-grading of \eqref{eq: A action}, while $D$-modules and Cherednik modules carry the $H_{\mathfrak{gl}_{N_r}}$ and $F^0$ gradings tracked by $q,h$. The conversion is fixed by Theorem \ref{thm:ch Lc Q} and Remark \ref{rem:grading-conversion}.
	\subsection{Hall convolution}
	For each $r\geq0$, let $\ol{\Cfrak}_r\subset\mathfrak{gl}_r^2$ be the commuting scheme. We equip
	\begin{align}
		\bigoplus_{r\geq0}K^{\GL_r\times A}(\ol{\Cfrak}_r)_\loc
		\label{def: K loc}
	\end{align}
	with the Hall convolution product $*$ of \cite[Sections 7.2--7.5]{SV_Hilbert} and \cite[Section 5.1(i)]{BLM}. 

	For $r\geq1$, write $N_r:=rn$. With respect to the direct sum $\mathfrak{gl}_{N_r}=\mathfrak{sl}_{N_r}\oplus\C I_{N_r}$, we take 
	\begin{align*}
		\wt{\Nbb}_r&:=\Nbb_{rm,rn}\boxtimes\Ocal_\C,\\
		\wt{\Nbb}_{r,\tau}&:=\Hom_{S_r}(V_\tau,\wt{\Nbb}_r)
		\quad \text{ for }\tau\vdash r.
	\end{align*}
	Put $K_r:=r(n-1)$ and use the normalized filtration $F_j^0=F_{j+K_r}^{\Hrm}$, with the factor $\Ocal_\C$ generated in filtration degree zero. For the closed embedding $i_r:\ol{\Cfrak}_{N_r}\hookrightarrow T^*\mathfrak{gl}_{N_r}$, Proposition \ref{prop: support gr Nc} gives a coherent sheaf $\wt{\Scal}_{r,\tau}$ characterized by
$(i_r)_*\wt{\Scal}_{r,\tau}=\wt\gr^{F^0}\wt{\Nbb}_{r,\tau}$
	such that 
	\begin{align*}
		\ch_{\GL_{N_r}\times A}\wt{\Scal}_{r,\tau}
		=\frac{1}{1-q}\,
		\ch_{\GL_{N_r}\times A}\Hom_{S_r}(V_\tau,\Scal_{rm,rn}). 
	\end{align*}
	
	Define the $k$-linear map
	\begin{align*}
		\Psi:\Lambda&\rightarrow  \bigoplus_{r\geq0}K^{\GL_r\times A}(\ol{\Cfrak}_r)_\loc,
		&\Psi(s_\tau)&:=q^{-K_r}[\wt{\Scal}_{r,\tau}]
		\quad \forall \tau\vdash r,
		&\Psi(1)&:=1.
	\end{align*}
	Thus symmetric functions of degree $r$ map to the rank-$N_r$ summand. 
	By Remark \ref{rem:grading-conversion}, the $(H_{\mathfrak{gl}_{N_r}},F^0)$-bigraded character of $\wt{\Nbb}_{r,\tau}$ becomes $q^{r/2}\ch\Psi(s_\tau)$ after setting $h=qt$.
	
	\begin{prop}\label{prop:Psi}
		The map $\Psi$ is an algebra homomorphism.
	\end{prop}
	\begin{proof}
		Let $\alpha\vdash a$ and $\beta\vdash b$. We induce from the Levi $\GL_{N_a}\times\GL_{N_b}\subset\GL_{N_{a+b}}$. By transitivity of induction, $\Ind(\wt{\Nbb}_a\boxtimes\wt{\Nbb}_b)\cong\wt{\Nbb}_{a+b}$. On the stratum of matrices with $a+b$ distinct eigenvalues, each having a single Jordan block of size $n$, the deck cover records an ordering of these blocks. Permutations within the first $a$ and last $b$ blocks give the restriction of the $S_{a+b}$-action to $S_a\times S_b$. This compatibility extends to the minimal extensions by Theorem \ref{prop_Sd_action}. Taking multiplicity spaces and applying Frobenius reciprocity gives
		\begin{align}
			\Ind(\wt{\Nbb}_{a,\alpha}\boxtimes\wt{\Nbb}_{b,\beta})
			\cong\bigoplus_{\tau\vdash a+b}
			\wt{\Nbb}_{a+b,\tau}^{\oplus c_{\alpha\beta}^{\tau}},
			\label{eqn:decompose_induction}
		\end{align}
		where $c_{\alpha\beta}^{\tau}$ are the Littlewood-Richardson coefficients, defined by $s_\alpha s_\beta=\sum_\tau c_{\alpha\beta}^{\tau}s_\tau$. 	
		By \cite[Lemma 5.2]{BLM}, the relation \eqref{eqn:decompose_induction} passes to associated graded sheaves on $\ol{\Cfrak}_{N_{a+b}}$ and gives the proposition. 
	\end{proof}
	
	\subsection{Reduction to the trivial component}
	\label{section: EHA}
	Let $K$ be the subalgebra of \eqref{def: K loc} generated by $z^a[\Ocal_{\ol{\Cfrak}_1}]$, $a\in\Z$, where $z$ is the standard character of $\GL_1$. We now use the geometric realization of $\Ecal^+$ from \cite{SV_Hilbert}:
	\begin{align*}
		C:\Ecal^{+,\mathrm{op}}\xrightarrow{\sim}K,
		\qquad C(Q_{1,a})=(1-q)(1-t)z^a[\Ocal_{\ol{\Cfrak}_1}].
	\end{align*}

	Write $h_r,e_r$ for the homogeneous and elementary symmetric functions. For coprime $u>0$ and $v\in\Z$, the formal alphabet $X^{u,v}$ identifies $\Lambda$ with the commutative subalgebra $\Lambda(X^{u,v})\subset\Ecal^+$ of slope $v/u$, where $p_r[X^{u,v}]$ has bidegree $(ru,rv)$ and $Q_{u,v}=p_1[-(1-q)(1-t)X^{u,v}]$.  We make the following reduction.
	\begin{prop}
		\label{prop:generator}
		For a fixed $d$, each of the following two families of identities  
		\begin{align}
			q^{-K_r}[\wt{\Scal}_{r,(1^r)}]
			&=C(e_r[(t-1)X^{n,m}])\quad \text{ for all }1\leq r\leq d ,
			\label{eq:antisym-class}\\
			q^{-K_r}[\wt{\Scal}_{r,(r)}]
			&=C(h_r[(t-1)X^{n,m}])\quad \text{ for all }1\leq r\leq d 
			\label{eq:symmetric-class}
		\end{align}
		is equivalent to the identities
		\begin{align}
			\Psi(s_\tau)=C(s_\tau[(t-1)X^{n,m}])\quad \text{ for all }|\tau|\leq d.
			\label{eq:all-Schur-classes}
		\end{align}
		When \eqref{eq:all-Schur-classes} hold, we have 
		\begin{align}
			\sum_{\tau\vdash d}\Psi(s_\tau)s_\tau[Y]=C(h_d[(t-1)X^{n,m}Y]).
			\label{eq:Sd-class-Cauchy}
		\end{align}
	\end{prop}
	\begin{proof}
		For the first statement, use Proposition \ref{prop:Psi} and the Jacobi-Trudi identities,
		\begin{align*}
			s_\tau
			=\det(h_{\tau_i-i+j})_{1\leq i,j\leq\ell(\tau)}
			=\det(e_{\tau_i^t-i+j})_{1\leq i,j\leq\tau_1}.
		\end{align*}
		 For \eqref{eq:Sd-class-Cauchy}, apply the Cauchy identity
		\begin{align*}
			C(h_d[(t-1)X^{n,m}Y])=\sum_{\tau\vdash d}C(s_\tau[(t-1)X^{n,m}])s_\tau[Y].   
		\end{align*}
	\end{proof}

	\subsection{Character of the trivial component}
	For the alphabet $Z=(z_1,\ldots,z_{N_r})$ for characters of $\GL_{N_r}$, set
	\begin{align*}
		A_r(q,h;Z)&:=\ch_{\GL_{N_r}\times A}
		\Gamma(\mathfrak{gl}_{N_r},\wt{\Nbb}_{r,(r)}).
	\end{align*}
	By Remark \ref{rem:grading-conversion}, $q^{-r/2}A_r(q,qt;Z)=\ch\Psi(h_r)$.  
    
	For two size-$r$ alphabets $U,V$, define
	\begin{align}
		D_r(U,V):=
		\frac{t^{-\binom r2}}{\delta(U)\delta(V)}
		\det_{1\leq i,j\leq r}\left(
		\frac{1-qt}{(1-tU_iV_j)(1-q/(U_iV_j))}\right),
		\quad \delta(U):=\prod_{i<j}(U_i-U_j).
		\label{eq:Sd-determinant}
	\end{align}
	Write
	\begin{align*}
		b_i:=\lceil im/n\rceil-\lceil(i-1)m/n\rceil,
		\quad g(a):=\frac{(1-qta)(1-q/a)}{1-a^{-1}},\quad  \Delta_{q,N_r}(Z):=\prod_{i,j=1}^{N_r}(1-qz_i/z_j).
	\end{align*}
	
    	\begin{prop}\label{prop: char triv part}
		For every $r\geq1$, the character $A_r$ satisfies
		\begin{align}
			\Delta_{q,N_r}(Z)A_r(q,qt;Z)
			={}&q^{r/2+n\binom r2}t^{-K_r}
			\sum_{\substack{Z_1\sqcup\cdots\sqcup Z_n=Z\\|Z_1|=\cdots=|Z_n|=r}}
			\prod_{i=1}^n(\det Z_i)^{b_i}
			\prod_{i=1}^{n-1}D_r(Z_i^{-1},Z_{i+1})\notag\\
			&\quad\cdot
			\prod_{j\geq i+2}\prod_{\substack{z\in Z_i\\z'\in Z_j}}
			(1-tz'/z)^{-1}
			\prod_{i<j}\prod_{\substack{z\in Z_i\\z'\in Z_j}}
			g(z'/z).
			\label{eq:rectangle-F}
		\end{align}
		The sum is over ordered decompositions of $Z$ into size-$r$ alphabets. We write $\det Z_i=\prod_{z\in Z_i}z$ and $Z_i^{-1}=\{z^{-1}:z\in Z_i\}$. All block indices lie between $1$ and $n$.
	\end{prop}
	\begin{proof}
		Apply Theorem \ref{thm:gr N Ac}(2) with $d=r$, $\gfrak=\mathfrak{sl}_{N_r}$, and $R=r^2\binom n2$. We compute the alternating character of $R\Gamma(\Acal_{rm,rn}^{(r)})$ from its fiber over the standard flag in $G/P$, tensor with  $\Ocal_\C$. By \eqref{eqn:rectangular-dg-module}, the restriction is  $\Cscr_{P}^\vee$ tensored with the relative canonical line bundle and $\wedge^\bullet\ufrak_{P}^*$. 
		
		\vspace{5pt}
		\emph{The $D_r$-term.} 	Let $U,V$ record the left and right torus weights. Then the torus character of $\Ocal_{\operatorname{Mat}_r}[\det^{-1}]\det^{-\alpha}$ is
		$\sum_{\lambda_1\geq\cdots\geq\lambda_r}
			s_\lambda(U)s_\lambda(V).$ 
		The term indexed by $\lambda$ has relative Euler degree $|\lambda|=\sum_a\lambda_a$. For rational $0<\alpha<1$, \cite[Theorem 1.4 and Proposition 5.1]{lorincz-yang}  give the first Hodge degree in which this summand occurs as
		\begin{align*}
			d_\lambda
			:=\sum_{a=1}^r\max\{\lceil\alpha-(1+\lambda_a+r-a)\rceil,0\}
			=\sum_{a=1}^r\max\{-\lambda_a-r+a,0\}.
		\end{align*}
		Since the rational Schur decomposition is multiplicity-free, this is its degree in the associated graded. The Fourier regrading \eqref{eqn:rectangular-fourier-grading} changes its relative $H$ and Hodge degrees from $(|\lambda|,d_\lambda)$ to $(-|\lambda|,d_\lambda+|\lambda|)$. Indeed, the Euler degree changes sign, and $\lceil\beta+|\lambda|\rceil-\lceil\beta\rceil=|\lambda|$. Consequently, setting $h=qt$ the block $(i,i+1)$ contributes
		\begin{align*}
			\sum_{\lambda_1\geq\cdots\geq\lambda_r}
			q^{d_\lambda}t^{|\lambda|+d_\lambda}s_\lambda(Z_i^{-1})s_\lambda(Z_{i+1})
			=D_r(Z_i^{-1},Z_{i+1}).
		\end{align*}
		Indeed, expand
		\begin{align*}
			\frac{1-qt}{(1-tw)(1-q/w)}
			=\sum_{k\geq0}t^kw^k+\sum_{k<0}q^{-k}w^k.
		\end{align*}
		Cauchy-Binet and the Weyl character formula index the determinant expansion by $k_1>\cdots>k_r$, where $k_a=\lambda_a+r-a$. The $q$-exponent is $\sum_a\max\{-k_a,0\}=d_\lambda$, and the $t$-exponent, including $t^{-\binom r2}$, is $\sum_a\max\{k_a,0\}-\binom r2=|\lambda|+d_\lambda$. 
        
		\vspace{5pt}
		\emph{The $(1-tz'/z)^{-1}$ and $g(z'/z)$ terms.}
		After restriction to the Levi, the normal derivatives, the polynomial coordinates on the blocks $x_{ij}$ with $j\geq i+2$, and the Koszul exterior factor contribute $\Sym(\mathfrak{gl}_{N_r}/\ufrak_P)$, $\Sym([\ufrak_P,\ufrak_P]^*)$, and $\wedge^\bullet\ufrak_P^*$, respectively. Their characters are respectively 
		\begin{align*}
			&\prod_i\Delta_{q,r}(Z_i)^{-1}\prod_{i<j}\prod_{\substack{z\in Z_i\\z'\in Z_j}}(1-qz'/z)^{-1},\quad 
			\prod_{j\geq i+2}\prod_{\substack{z\in Z_i\\z'\in Z_j}}(1-tz'/z)^{-1},
			\quad
			\prod_{i<j}\prod_{\substack{z\in Z_i\\z'\in Z_j}}(1-qtz'/z).
		\end{align*}
		
		\vspace{5pt}
		\emph{The line bundle and grading shift.}
		At the standard flag, $\omega_{\ufrak_P/\gfrak}$ has character $\det(\ufrak_P^*)$, and $\omega_{G/P}$ has character $\det(\ufrak_P^*)^{-1}$. They cancel, leaving the generator of character $\prod_i(\det Z_i)^{b_i}$. Since $m$ and $n$ are coprime, the numbers $\alpha_i=\lceil im/n\rceil-im/n$ for $1\leq i<n$ permute $1/n,\ldots,(n-1)/n$. Thus the common generator $\bigotimes_{i=1}^{n-1}\det(A_i)^{-\alpha_i}$ has $H_{\ufrak_P}$-degree $R/2-r\sum_i\alpha_i=(R-K_r)/2$ and Hodge degree zero. Pushforward along the closed embedding subtracts $(\dim\gfrak-R)/2$ from its $H$-degree and adds $\dim\gfrak-R$ to its Hodge degree. Consequently, \eqref{eqn:rectangular-fourier-grading} gives
		\begin{align*}
			\left(-\frac{K_r}{2}-\frac{\dim\gfrak}{2}+R,\ \dim\gfrak-R\right)
			\longmapsto\left(\frac{K_r}{2}+\frac{\dim\gfrak}{2}-R,\ R\right),
		\end{align*}
		where the Hodge degree uses $\lceil-K_r/2\rceil=-\lfloor K_r/2\rfloor$. Including the scalar constant of $H_\C$-degree $1/2$, the shift $\langle-R\rangle$ in Theorem \ref{thm:gr N Ac}(2), and the normalization $F_j^0=F_{j+K_r}^{\Hrm}$ gives the common factor
		\begin{align*}
			q^{K_r/2+N_r^2/2-R}(qt)^{-K_r}
			=q^{N_r^2/2-K_r/2-R}t^{-K_r}
			=q^{r/2+n\binom r2}t^{-K_r}.
		\end{align*}
		
		\vspace{5pt}
		\emph{Induction from $P$ to $G$.}
		Localization on $G/P$ produces the term $\prod_{i<j}\prod_{z\in Z_i,z'\in Z_j}(1-z/z')^{-1}$. Multiply the resulting character by $\Delta_{q,N_r}(Z)$. The diagonal-block factors $\Delta_{q,r}(Z_i)$ cancel, and for each pair of distinct blocks the remaining factors combine as
		\begin{align*}
			(1-qa)(1-q/a)\frac{1-qta}{(1-qa)(1-a^{-1})}=g(a),
			\quad a=z'/z.
		\end{align*}
    		Altogether, we obtain \eqref{eq:rectangle-F}. 
	\end{proof}

	\subsection{Comparison with the complete shuffle generators}
	
	\begin{thm}\label{thm:Sd char}
		Write $\ch:=\ch_{\GL_{N_r}\times A}$. For every $r\geq1$,
		\begin{align*}
			q^{-r/2}A_r(q,qt;Z)=\ch C(h_r[(t-1)X^{n,m}])(Z).
		\end{align*}
		Consequently \eqref{eq:symmetric-class} and \eqref{eq:all-Schur-classes} hold.
	\end{thm}
	\begin{proof}
			Recall that $N_r=rn$ and $K_r=r(n-1)$. In this proof, write $\mathsf{A}_r:=\Psi(h_r)$,  $\mathsf{B}_r:=C(h_r[(t-1)X^{n,m}])$ and $\mathsf{A}_0=\mathsf{B}_0=1$. Since  $\ch(\mathsf{A}_r)(Z)=q^{-r/2}A_r(q,qt;Z)$, it suffices to prove $\mathsf{A}_r=\mathsf{B}_r$. We argue in steps, following \cite[Proposition 5.2, Lemma 5.5 and Section 6.11]{negutshuffle}.
			
			\smallskip\noindent\emph{Step 1: We show that $\mathsf{A}_r,\mathsf{B}_r\in C(\Lambda(X^{n,m}))$.}
			This is immediate for $\mathsf{B}_r$. For $\mathsf{A}_r$, nilpotent cotangent support \eqref{formula:u} and \cite[Section 1.2 and Theorem 1.4]{negut-integral} implies that $\mathsf{A}_r$ lies in the image of $C$. Its character has degree $rm$ in $N_r=rn$ variables, so the slope criterion of \cite[Sections 2.6 and 6.11]{negutshuffle} reduces the claim to
			\begin{align}
				\deg_\xi\left(q^{-\binom{N_r}{2}}\Delta_{q,N_r}(\xi U,V)\ch(\mathsf{A}_r)(\xi U,V)\right)\leq\frac mn|U|,
				\label{eqn: slope}
			\end{align}
			where $Z=U\sqcup V$.
			
			In a summand of \eqref{eq:rectangle-F}, write $Z_i=\xi U_i\sqcup V_i$, with $|U_i|=r_i$. Put $\ell=|r_i-r_{i+1}|$. In the determinant defining $D_r(Z_i^{-1},Z_{i+1})$, expansion by minors forces a minor containing both $U$ and $V$ to be of size at least $\ell$. Its entries have expansions in positive powers of $\xi^{-1}$ whose coefficients factor into row and column factors. Nonzero minors therefore require distinct exponents, with sum at least $\ell(\ell+1)/2$. Including the Vandermonde factors gives
			\begin{align*}
				\deg_\xi D_r(Z_i^{-1},Z_{i+1})
				&\leq-\frac{\ell(\ell+1)}2+\binom{r_i}{2}-\binom{r_{i+1}}2-r_{i+1}(r-r_{i+1})\\
				&=-r_{i+1}(r-r_i)-(r_i-r_{i+1})_+,
			\end{align*}
			where $a_+:=\max\{a,0\}$. Since $g(a)\sim-qt\,a$ at infinity and $g(a)\sim q$ at zero, the $g$-factors for $i<j$ contribute degree $r_j(r-r_i)$. For nonadjacent blocks this is canceled by $(1-tz'/z)^{-1}$. For adjacent blocks the bound above leaves $-(r_i-r_{i+1})_+$. Hence
			\begin{align*}
				\deg_\xi(\text{each summand})\leq\sum_i b_i r_i-\sum_{i<n}(r_i-r_{i+1})_+.
			\end{align*}
			Decompose each level set $\{i:r_i\geq h\}$ into maximal intervals $[l,j]$. Such an interval contributes
			\begin{align*}
				\sum_{i=l}^j b_i-\mathbf1_{j<n}
				=\left\lceil\frac{jm}{n}\right\rceil-\left\lceil\frac{(l-1)m}{n}\right\rceil-\mathbf1_{j<n}
				\leq\frac{(j-l+1)m}{n},
			\end{align*}
			with equality only for $[l,j]=[1,n]$: use $\lceil jm/n\rceil-1<jm/n$ when $j<n$, and coprimality when $j=n$ and $l>1$. Summing proves \eqref{eqn: slope}, and equality holds only when $r_1=\cdots=r_n$.
			
			\smallskip\noindent\emph{Step 2: We show that $\mathsf{A}_r$ and $\mathsf{B}_r$ satisfy the coproduct identities.}
			Transport the symmetric-function coproduct $p_j\mapsto p_j\otimes1+1\otimes p_j$ through $C$, and let $\Delta^{\mathrm{cop}}_{ns,nr}$ denote its component with output matrix sizes $ns,nr$. For $s,r\geq1$, we claim
			\begin{align*}
				\Delta^{\mathrm{cop}}_{ns,nr}(\mathsf{A}_{s+r})=\mathsf{A}_s\otimes \mathsf{A}_r,
				\quad\Delta^{\mathrm{cop}}_{ns,nr}(\mathsf{B}_{s+r})=\mathsf{B}_s\otimes \mathsf{B}_r.
			\end{align*}
			The second identity follows from the coproduct of homogeneous symmetric functions. By \cite[Proposition 5.2 and Section 6.11]{negutshuffle}, the first is equivalent to
			\begin{align}
				&\lim_{\xi\to\infty}\xi^{-sm}q^{-\binom{N_{s+r}}2}\Delta_{q,N_{s+r}}(\xi U,V)\ch(\mathsf{A}_{s+r})(\xi U,V)\notag\\
				=&\left(q^{-\binom{N_s}2}\Delta_{q,N_s}(U)\ch(\mathsf{A}_s)(U)\right)
				\left(q^{-\binom{N_r}2}\Delta_{q,N_r}(V)\ch(\mathsf{A}_r)(V)\right),
				\label{eq:rectangle-grouplike}
			\end{align}
			where $|U|=ns$, $|V|=nr$.
			
			By Step 1, only decompositions $Z_i=\xi U_i\sqcup V_i$ with $|U_i|=s$, $|V_i|=r$ survive. The off-diagonal determinant entries tend to zero, giving
			\begin{align}
            \label{eqn: D}
				D_{s+r}(Z_i^{-1},Z_{i+1})
				\sim(-t)^{-sr}\xi^{-sr}\frac{(\det V_i)^s}{(\det U_{i+1})^r}
				D_s(U_i^{-1},U_{i+1})D_r(V_i^{-1},V_{i+1}).
			\end{align}
			For every $i<j$, the $g$ terms containing both $U$ and $V$ variables contribute
\begin{align*}
    (-qt)^{sr}q^{sr}\xi^{sr}\frac{(\det U_j)^r}{(\det V_i)^s}.
\end{align*}
The factors involving $t$, $\xi$, and the determinants cancel against the prefactor in \eqref{eqn: D} for adjacent blocks, or against the $(1-tz'/z)^{-1}$ factors containing both $U$ and $V$ variables for nonadjacent blocks, leaving $q^{2sr}$ in either case. The determinant monomials contribute $\xi^{sm}$, while
			\begin{align*}
				q^{-(s+r)^2\binom n2}t^{-K_{s+r}}q^{2sr\binom n2}
				=q^{-s^2\binom n2}t^{-K_s}q^{-r^2\binom n2}t^{-K_r}.
			\end{align*}
			All remaining factors belong separately to the $U$- and $V$-summands. Their ordered decompositions are independent, proving \eqref{eq:rectangle-grouplike}.
			
			\smallskip\noindent\emph{Step 3: Induction on $r$.}
			The coprime identity from \cite[proof of Proposition 5.4]{ma2}, after tensoring with $\Ocal_\C$ and multiplying by $q^{-K_1}=q^{-(n-1)}$, gives
			\begin{align*}
				\mathsf{A}_1=(1-q)^{-1}C(Q_{n,m})=C(h_1[(t-1)X^{n,m}])=\mathsf{B}_1.
			\end{align*}
			If $\mathsf{A}_j=\mathsf{B}_j$ for $j<r$, Step 2 implies that $C^{-1}(\mathsf{A}_r-\mathsf{B}_r)$ is primitive. Since the degree-$r$ primitive subspace of $\Lambda$ is spanned by $p_r$, we have
			\begin{align*}
				\mathsf{A}_r-\mathsf{B}_r=\gamma_r C(p_r[X^{n,m}]),\quad\gamma_r\in k.
			\end{align*}
			
			\smallskip\noindent\emph{Step 4: For $r>1$, we prove $\gamma_r=0$ by specific evaluation.}
			For an alphabet $Z=(z_1,\ldots,z_k)$, put
			\begin{align*}
				\Delta_{t,k}(Z):=\prod_{i,j=1}^k(1-tz_i/z_j),\quad
				\Delta_k:=\Delta_{q,k}\Delta_{t,k}.
			\end{align*}
			
			Expand the adjacent determinants in \eqref{eq:rectangle-F}. The denominators $1-tz'/z$ cancel against $\Delta_{t,N_r}$, and $1-qz/z'$ against the $g$-factors; the remaining denominators are regular at the specialization. For each $a$, canceling the zero factor $1-t(t^ax)/(t^{a+1}x)$ requires $t^{a+1}x$ to lie in a block above $t^ax$. A nonzero term would therefore require $N_r>n$ strictly decreasing block indices. Hence
			\begin{align*}
				(\Delta_{N_r}\ch(\mathsf{A}_r))(xt^{N_r-1},\ldots,x)=0.
			\end{align*}
			
			For $\mathsf{B}_r$ for $r>1$, the shuffle formula \cite[(46)]{gn_commuting_idempotent}, gives
\begin{align*}
    (\Delta_{N_r}\ch(\mathsf{B}_r))(xt^{N_r-1},\ldots,x)=0.
\end{align*}
			By \cite[Section 6.6 and Proposition 6.9]{negutshuffle},  we have
\begin{align*}
    \bigl(\Delta_{N_r}\ch C(p_r[X^{n,m}])\bigr)(xt^{N_r-1},\ldots,x)\ne0.
\end{align*}
			Applying these specializations to $\mathsf{A}_r-\mathsf{B}_r=\gamma_r C(p_r[X^{n,m}])$ gives $\gamma_r=0$ and completes the induction.
		\end{proof}

	\subsection{$S_d$-characters}
	Let $\alpha:\Ecal\to\End_k(\Lambda)$ be the polynomial representation of the full elliptic Hall algebra, let $\tau:\Ecal^{+,\mathrm{op}}\to\Ecal$ exchange slopes by $\tau(Q_{u,v})=Q_{v,u}$, and let $\rho:K\to\End_k(\Lambda)$ be the geometric Hall action on the full Hilbert schemes, transported by the Frobenius characteristic. Define $\nabla$ on the modified Macdonald basis in $X$ by $\nabla\wt H_\lambda=q^{n(\lambda^t)}t^{n(\lambda)}\wt H_\lambda$, where $n(\lambda)=\sum_i(i-1)\lambda_i$. The action comparison \cite[Proposition 5.1]{ma2} is
	\begin{align}
		\rho(C(f))=\nabla\,\alpha(\tau(f))\,\nabla^{-1},\qquad f\in\Ecal^+.
		\label{eq:Hall-action-comparison}
	\end{align}
	Return to $M=dm$ and $N=dn$. For $\tau\vdash d$, put
	\begin{align*}
		\Fcal_{M,N,\tau}&:=\Hom_{S_d}(V_\tau,\Fcal_{M,N}).
	\end{align*}
	
	\begin{thm}\label{thm:Sd-joint-character}
		With respect to the Euler and filtration gradings,
		\begin{align}
			\ch_{S_N\times S_d\times A}(\gr^0\Lrm_{M,N})(X,Y;q,h)
			=q^{(d-1)/2}(1-q)
			\left.\omega_X\bigl(h_d[(t-1)X^{m,n}Y]\cdot1\bigr)\right|_{t=h/q}.
			\label{eq:joint-candidate}
		\end{align}
		Moreover, 
		\begin{align}
			\kappa([\Fcal_{M,N,\tau}])
			=q^{N-d}(1-q)\nabla\bigl(s_\tau[(t-1)X^{m,n}]\cdot1\bigr).
			\label{eq:Sd-individual-kappa}
		\end{align}
	\end{thm}
	\begin{proof}
		By the definition of $\Psi$ and \eqref{eq:all-Schur-classes},
		$[\wt{\Scal}_{d,\tau}]=q^{N-d}C(s_\tau[(t-1)X^{n,m}]).$ 
		Setting the scalar base coordinate to zero multiplies this class by $1-q$ and gives the class of $\Hom_{S_d}(V_\tau,\Scal_{M,N})$, extended by zero to $\ol{\Cfrak}_N$.

		By \cite[Proposition 5.2]{ma2} and \cite[Proposition 3.5.2]{BHMPS},
		\begin{align*}
			\left(\ch_{\GL_N\times A}C(s_\tau[(t-1)X^{n,m}])\right)_{\pol}
			=\omega_X\bigl(s_\tau[(t-1)X^{n,m}]\cdot1\bigr).
		\end{align*}
		Summing over $s_\tau[Y]$, using the Cauchy identity, and noting the factor $q^{(d-1)/2}$ from Theorem \ref{thm:ch Lc Q} gives
		\begin{align*}
			\ch_{S_M\times S_d\times A}(\gr^0\Lrm_{N,M})(X,Y;q,h)
			=q^{(d-1)/2}(1-q)
			\left.\omega_X\bigl(h_d[(t-1)X^{n,m}Y]\cdot1\bigr)\right|_{t=h/q}.
		\end{align*}
		Exchanging $M,N$ gives \eqref{eq:joint-candidate}.

	For \eqref{eq:Sd-individual-kappa}, $\kappa([\Fcal_{M,N,\tau}])$ equals the Hall action of the class $\wt{\Scal}_{r,\tau}$ on $1$, as in the proof of \cite[Proposition 5.4]{ma2}. Applying \eqref{eq:Hall-action-comparison} and $\nabla^{-1}1=1$ gives the result.
	\end{proof}
	
	\subsection{Comparison of the two joint characters}\label{sec:character-comparison}

Recall that $\delta_{M,N}:=\frac{(M-1)(N-1)+d-1}{2}.$

\begin{thm}\cite{KO2} \label{thm: ko}
\begin{align*}
t^{\delta_{M,N}}
 \left.\omega_X\omega_Y\bigl(h_d[(t-1)X^{m,n}Y]\cdot1\bigr)\right|_{t\mapsto t^{-1}}=\sum_{\mu\vdash N}\ \sum_{\substack{(\Delta,f)\colon\Delta\in I_{M,N}^+\\ \text{ dominant, with content }\mu}}
 q^{\area(\Delta)}t^{\codinv(\Delta,f)}
 \frac{\Xi_{\fm_{(\Delta,f)}}(Y,t)}{\Aut_t(\alpha_{(\Delta,f)})}m_\mu(X).
\end{align*}
\end{thm}
Keep the alphabets of Theorem \ref{thm: ASF main} with $X$ and $Y$ recording the $S_N$- and $S_d$-actions, respectively. The variable $q$ records the Hikita filtration, and $t^i$ records $H^{BM}_{2i}$. The involutions $\omega_X$ and $\omega_Y$ transpose Schur functions in the indicated alphabet.
\begin{thm}\label{thm:character-comparison}
We have the following identities 
\begin{align*}
&t^{-\delta_{M,N}}\Frob_{S_N\times S_d}\bigl(\gr_F H_*^{BM}(\Fl_\gamma^+)\bigr)(X,Y;q,t)\notag\\
=&\left.\omega_X\omega_Y
 \bigl(h_d[(t-1)X^{m,n}Y]\cdot1\bigr)\right|_{t\mapsto t^{-1}}\\
=&\frac{q^{-(d-1)/2}}{1-q}\,
 \omega_Y\ch_{S_N\times S_d\times A}
 \bigl(\gr^0\Lrm_{M,N}\bigr)(X,Y;q,q/t)
\label{eq:ASF-Cherednik-normalized}\\
=&\frac{q^{d-N}}{1-q}\,
 \left.\omega_Y
 \ch_{S_N\times S_d\times A}R\Gamma
 \bigl(\Hilb^N,\Pcal^\vee\otimes\Fcal_{M,N}\bigr)
 \right|_{t\mapsto t^{-1}}.
\end{align*}

\end{thm}
\begin{proof}
Theorem \ref{thm: ASF main} and Theorem \ref{thm: ko} give the first equality. The second equality is Theorem \ref{thm:Sd-joint-character}. 

For the third equality, multiply \eqref{eq:Sd-individual-kappa} by $s_\tau[Y]$ and sum over $\tau\vdash d$. By \eqref{iso: K Lambda} and the Cauchy identity, this gives
\begin{align*}
\ch_{S_N\times S_d\times A}R\Gamma\bigl(\Hilb^N,\Pcal\otimes\Fcal_{M,N}\bigr)
=q^{N-d}(1-q)\nabla\bigl(h_d[(t-1)X^{m,n}Y]\cdot1\bigr).
\end{align*}
The equivariant perfect pairing on the Procesi bundle gives
\begin{align*}
\Pcal^\vee
\cong \Pcal\otimes \Ocal(-1)\otimes\operatorname{sgn}_{S_N}.
\end{align*}
Under $\kappa$, tensoring by $\Ocal(-1)$ corresponds to $\nabla^{-1}$, while tensoring by $\operatorname{sgn}$ corresponds to $\omega_X$. Hence the third equality follows.
\end{proof}

\begin{rem}
   Similar to Corollary \ref{cor: ASF equivariant main}, we have also have an analogue of Theorem \ref{thm:character-comparison} in terms of $h_d[X^{m,n}Y]$. Accordingly, on the coherent side, one takes the minimal extension of $\wt{\Ebb}_m^\circ$ from $\wt{\gfrak}_d^\circ$ to $\gfrak\times \tfrak_d$. The pushforward of its associated graded to $T^*\gfrak$ would have singular support 
   $\{(x,y)\in\gfrak_d\times\gfrak_d\mid[x,y]=0\bigr\}_\red$. In particular, when $m=n=1$, the resulting sheaf on $\Hilb^d$ is the Procesi bundle according to \cite{isospec}. 
\end{rem}

\end{document}